\documentclass[11pt]{article}
\usepackage[utf8]{inputenc}
\usepackage{fullpage}
\usepackage{algorithm}
\usepackage{algorithmic}
\usepackage{boxedminipage}
\usepackage{comment}
\usepackage{eqparbox}
\usepackage{multirow}
\usepackage{caption}
\usepackage{subcaption}

\usepackage{amsmath}
\usepackage{amsfonts}
\usepackage{graphicx}

\newlength{\commentindent}
\usepackage{array}
\newcolumntype{R}[1]{>{\raggedleft\let\newline\\\arraybackslash\hspace{0pt}}m{#1}}

\newcommand\BibTeX{{\rmfamily B\kern-.05em \textsc{i\kern-.025em b}\kern-.08em
T\kern-.1667em\lower.7ex\hbox{E}\kern-.125emX}}

\newcommand{\EN}{\mathcal{E}}
\newcommand{\dE}{\partial \EN}

\newcommand{\jump}[1]{[\![#1]\!]}
\newcommand{\avg}[1]{\{\!\!\{#1\}\!\!\}}

\usepackage[ backend=bibtex, style=authoryear, url=false, isbn=false, maxcitenames=2, maxbibnames=10,giveninits=true,uniquename=init, natbib]{biblatex}
\renewbibmacro{in:}{\ifentrytype{article}{}{\printtext{\bibstring{in}\intitlepunct}}}

\renewbibmacro*{volume+number+eid}{%
  \setunit{\addcomma\space}%
  \printfield{volume}%
  \setunit*{\addnbthinspace}
  \printfield{number}%
  \setunit{\addcomma\space}%
  \printfield{eid}}
\DeclareFieldFormat[article]{number}{\mkbibparens{#1}}

\begin{document}

\title{Discontinuous Galerkin Discretizations of the Boltzmann Equations in 2D: semi-analytic time stepping and absorbing boundary layers} 

\author{A. Karakus\thanks{Department of Mathematics, Virginia Tech, 225 Stanger Street, Blacksburg, VA 24061, USA, akarakus@vt.edu},  N. Chalmers\textsuperscript{*}, J.S. Hesthaven \thanks{Computational Mathematics and Simulation Science, EPFL SB, MATH MCSS, MA C2 652 (Bâtiment MA), Station 8, CH-1015 Lausanne, Switzerland} \& T.~Warburton\textsuperscript{*}}

\maketitle

\begin{abstract}
We present an efficient nodal discontinuous Galerkin method for approximating nearly incompressible flows using the Boltzmann equations. The equations are discretized with Hermite polynomials in velocity space yielding a first order conservation law. A stabilized unsplit perfectly matching layer (PML) formulation is introduced for the resulting nonlinear flow equations. The proposed PML equations exponentially absorb the difference between the nonlinear fluctuation and the prescribed mean flow. We introduce semi-analytic time discretization methods to improve the time step restrictions in small relaxation times. We also introduce a multirate semi-analytic Adams-Bashforth method which preserves efficiency in stiff regimes. Accuracy and performance of the method are tested using distinct cases including isothermal vortex, flow around square cylinder, and wall mounted square cylinder test cases.     
\end{abstract}

\tableofcontents

\section{Introduction}
The Boltzmann equations, based on kinetic theory, describe fluids at the microscopic level. It has been shown that the Boltzmann equations recover the Navier-Stokes equations in the low Mach limit \citep{chapman_mathematical_1970,cercignani_boltzmann_1988}. The Boltzmann equations are also used in describing rarefied flows \citep{yang_rarefied_1995}. The main difficulty encountered in studying the equations is the complex non-linear integral nature of the collision term, which is often replaced with statistical or relaxation models. In this work we adopt the Bhatnaggar-Gross-Krook (BGK) \citep{bhatnagar_model_1954} single rate relaxation approximation.  

Lattice Boltzmann methods (LBM) are widely used to discretize the simplified Boltzmann equations. The classical LBM is a first order, explicit, upwind finite difference scheme for the discrete Boltzmann equation where the continuous velocity space is reduced a finite number of velocities \citep{yu_viscous_2003}. Although the LBM has several advantages including simplicity, easy parallelization, and relatively low floating point operations per lattice node, it is limited to structured meshes, suffers instabilities at high Reynolds numbers, and is difficult to accurately enforce boundary conditions. Several approaches have been proposed to address these limitations through replacing the lattice with finite volume \citep{nannelli_lattice_1992,peng_lattice_1998}, finite element \citep{lee_characteristic_2001} and discontinuous Galerkin \citep{shi_discontinuous_2003,duster_high-order_2006,min_spectral-element_2011} methods. 

The time discretization of the Boltzmann equations in stiff regimes presents a computational challenge. In the small relaxation time regime, the time-scale of the collision operator dominates the transport of particles and forces the numerical methods to operate with small time discretization steps. Fully implicit time integration techniques are limited for most applications because of the cost of inversion of the nonlinear collision operator. Semi-analytic or exponential time discretization methods allow us to overcome this problem. Semi-analytic methods \citep{cox_exponential_2002, kassam_fourth-order_2005} are special numerical time discretization schemes in which the traditional linear system solves for implicit schemes are replaced with computing an action of a matrix exponential. These methods have also been developed for the discrete non-homogeneous \citep{dimarco_exponential_2011} and homogeneous \citep{li_exponential_2014} Boltzmann equations. 

Implicit-Explicit (IMEX) methods are also popular schemes to relax the time step restriction in stiff ODEs such as the discrete Boltzmann equations \citep{dimarco_implicit-explicit_2017}. We refer to \citep{dimarco_numerical_2014} for a survey of semi-analytic and implicit-explicit techniques for discrete velocity Boltzmann equations. In this study, we explore the performance of time discretization methods by employing semi-analytic and low-storage IMEX methods for the Boltzmann equations, which fully exploit the specific structure of the non-linear collision operator. 

Accurate representation of the underlying domain geometry or complex flow field often requires the use of unstructured grids clustered at some specific locations. If a classical explicit scheme is used, varying length and time scales in the model introduce Courant-Friedrichs-Lewy (CFL) type time step restriction, which must be enforced globally. Multirate time discretizations address this restriction by using different time steps for each grid partition, using only local CFL stability conditions. A coherent flux transport between the partitions preserves the order of accuracy and the stability of the method. Due to their inherent efficiency, various multirate methods have been developed based on Runge-Kutta  \citep{seny_multirate_2013,schlegel_multirate_2009,constantinescu_multirate_2007} and multistep \citep{godel2010gpu,sandu_multirate_2009}
schemes for purely hyperbolic equations. For the Boltzmann equations, time-scales of the collision operator may dominate the advective scales depending on the flow regime and local grid resolution. This phenomenon makes the classical multirate methods as the problem becomes globally stiff limiting the number of possible multirate partitions. In this study, we extend the semi-analytic Adams-Bashforth approach to multirate time discretization which preserves efficiency and accuracy in stiff regimes.      

Perfectly matching layers (PML) were first introduced by Berenger \citep{berenger_perfectly_1994} for the Maxwell equations and is one of the preferred techniques for the computation of wave problems in open domains. PMLs rely on the fact that absorbing material zones surrounding the computational domain are theoretically non-reflecting, irrespective of the frequency and angle of outgoing multi-dimensional linear waves. Due to its simplicity and performance, PMLs are used extensively for modelling many physical phenomena such as the linearized Euler equations \citep{hesthaven_analysis_1998,hu_perfectly_2005}, wave equations \citep{collino_application_2001, becache_mixed_2005,appelo_new_2006}, Schr\"{o}dinger equations \citep{zheng_perfectly_2007}, Boltzmann equations \citep{najafi-yazdi_absorbing_2012,sutti_analysis_2015},  nonlinear Euler equations, and Navier-Stokes equations \citep{hagstrom_experiments_2007,hu_absorbing_2008}. 
In the original PML formulation \citep{berenger_perfectly_1994}, the field variables are split into nonphysical components to incorporate the mathematical formulation for the desired absorption. This approach is therefore referred to as {\it split-field PML}. It has been shown that the classical split model is dynamically stable but only weakly well-posed \citep{abarbanel_mathematical_1997,becache_analysis_2002}. As a result unsplit formulations, based on the causal frequency dependent PML (referred  to as  {\it convolutional PML or  C-PML}), have been proposed and analyzed \citep{abarbanel_construction_1998,appelo_new_2006,komatitsch_unsplit_2007}. However, C-PML also manifests slowly growing instabilities especially in anisotropic media \citep{matzen_efficient_2011} and loses its absorption capability at low frequencies \citep{meza-fajardo_nonconvolutional_2008}.  The so-called M-PML method was proposed in \citep{meza-fajardo_nonconvolutional_2008}. This method is a multiaxial version of the classical split-field PML formulation but it uses a more general coordinate stretching with anisotropic damping. M-PML has been shown to be more stable and efficient than classical PML in the long term simulation of wave propagation in elastic and anisotropic media yet it has similar reflection properties with increased efficiency when compared with C-PML.   \citep{meza-fajardo_stability_2010,mezafajardo_study_2012}.

In this study, the Boltzmann equations under the BGK relaxation approximation are discretized in velocity space  using Galerkin approach with Hermite polynomials \citep{grad_kinetic_1949,tolke_discretization_2000}. A nodal discontinuous Galerkin method is used for the resulting first order system in terms of Hermite polynomial coefficients. We propose an M-PML formulation for the resulting DG method and demonstrate its accuracy in truncated domains. We also introduced single and multirate semi-analytic time discretization methods to increase the performance of temporal integration in stiff regimes. The remainder of this paper is structured as follows: Section \ref{Sec:Formulation} introduces the derivation of the continuous Boltzmann equations, the design of perfectly matching layers, and the discontinuous Galerkin discretization for the resulting first order system. Section \ref{Sec:TimeDiscretization} is devoted to time discretization strategies including semi-analytic, implicit-explicit and multirate methods. Section \ref{Sec:Implementation} briefly describes the implementation of the numerical scheme leveraging GPU acceleration. Finally, in Section \ref{Sec:NumericalTests} we show numerical results which validate the formulation and demonstrate the applicability and performance of the approach, before giving some concluding remarks in Section \ref{Sec:Conclusion}.

\section{Formulation}
\label{Sec:Formulation}
In this section we begin by describing the Galerkin-Boltzmann equations. We then explain our proposed unsplit M-PML formulation for the discrete system. Finally, we detail the high-order nodal DG discretization and definitions of the resulting discrete operators. 

\subsection{Galerkin-Boltzmann Equations}
The Boltzmann equations describe the time evolution of a phase-space distribution function, $f(\mathbf{x},\mathbf{v},t)$ which is a function of the spatial variable $\mathbf{x}$, microscopic particle velocity $\mathbf{v}$, and time $t$. Neglecting external particle acceleration and under the BGK single-rate relaxation approximation \citep{bhatnagar_model_1954}, the continuous Boltzmann-BGK equation reads 
\begin{equation}
\label{eq:Boltzmann-BGK}
    \frac{\partial f}{\partial t} + \mathbf{v}\cdot\nabla_\mathbf{x} f =\frac{\left(f^{eq}-f\right)}{\tau},
\end{equation}
where $\tau$ is the relaxation time and $f^{eq}$ is the equilibrium phase space density which attains the macroscopic density, $\rho$, through the Maxwell velocity distribution as follows,
\begin{eqnarray*}
\label{eq:Equilibrium-Distribution}
     f^{eq} = \frac{\rho}{2\pi RT} \exp\left(-\frac{\left(\mathbf{v}-\mathbf{u}\right)^2}{2RT}\right),
\end{eqnarray*}
where $R$, $T$, and $\mathbf{u}$ are the gas constant, temperature, and macroscopic vector velocity field, respectively. 

Following the work of T\"{o}lke et. al. \citep{tolke_discretization_2000} the velocity field $\mathbf{v}$ is approximated by a discrete set of Hermite polynomials in the velocity space. The order of the polynomials needs to be sufficiently large to recover macroscopic flow properties which are the moments of the phase space distribution function. To model isothermal and nearly incompressible flows, second or higher order Hermite polynomials are required. The application of Galerkin formalism and analytic integration of the weak form of \eqref{eq:Boltzmann-BGK} in velocity space leads to the following first order semi-discrete PDE
\begin{equation}
\label{eq:BoltzmannSystem}
\frac{\partial q}{\partial t}= A_x\frac{\partial q}{\partial x} + A_y\frac{\partial q}{\partial y} + \mathcal{N}(\mathbf{q}),
\end{equation}
where $\mathbf{q} = \mathbf{q}(\mathbf{x},t)$ is the vector of Hermite coefficients to be solved, $A_x,$ and $A_y$ are matrices giving the directional coefficients, and $\mathcal{N}$ is the non-linear collision term. 

For the remainder of this paper, we assume a spatial of dimension two for simplicity but note that generalization to the third dimension is straightforward. In two dimensions, $\mathbf{x} = [x,y]$ and we assume a second order velocity approximation so that $\mathbf{q}(\mathbf{x},t) = [ q_1(x,y,t), \ldots,q_6(x,y,t)]^T$, $A_x$, $A_y$, and  $\mathcal{N}$ are given by
\begin{eqnarray*}
A_x = -\sqrt{RT} \left(\begin{array}{cccccc}
0 & 1 & 0 & 0 & 0 & 0\\
1 & 0 & 0 & 0 & \sqrt{2} & 0\\
0 & 0 & 0 & 1 & 0 & 0\\
0 & 0 & 1 & 0 & 0 & 0\\
0 & \sqrt{2} & 0 & 0 & 0 & 0\\
0 & 0 & 0 & 0 & 0 & 0
\end{array}
\right),\;
A_y = -\sqrt{RT} \left(\begin{array}{cccccc}
0 & 0 & 1 & 0 & 0 & 0\\
0 & 0 & 0 & 1 & 0 & 0\\
1 & 0 & 0 & 0 & 0 & \sqrt{2} \\
0 & 1 & 0 & 0 & 0 & 0\\
0 & 0 & 0 & 0 & 0 & 0\\
0 & 0 & \sqrt{2} & 0 & 0 & 0
\end{array}
\right) 
\end{eqnarray*}

\begin{eqnarray}
\label{eq:nonlinearDef}
\mathcal{N} = -\frac{1}{\tau}\left(\begin{array}{cccccc}
0 &
0 &
0 &\
\left(q_4 - \frac{q_2q_3}{q_1}\right) &
\left(q_5 - \frac{q_2^2}{q_1\sqrt{2}}\right) &
\left(q_6 - \frac{q_3^2}{q_1\sqrt{2}}\right)
\end{array}
\right)^T,
\end{eqnarray}
where $c=\sqrt{RT}$ represents the speed of sound of the fluid. This particular form of the system has a symmetric advection operator coupled with a nonlinear collision source term.  

The system \eqref{eq:BoltzmannSystem} recovers the Navier-Stokes equations for low Mach number, nearly incompressible flows with kinematic viscosity, $\nu = \tau RT$. Macroscopic flow properties are computed using the moment of the distribution function as follows,
\begin{eqnarray*}
\rho = q_1, \; \rho u = \sqrt{RT}q_2, \;\rho v = \sqrt{RT}q_3.
\end{eqnarray*}
Similarly, components of deviatoric stress tensor is given by
\begin{eqnarray*}
\sigma_{11} =  -RT\left(\sqrt{2} q_5 - \frac{q_2^2}{q_1}\right),\;
\sigma_{22} =  -RT\left(\sqrt{2} q_6 - \frac{q_3^2}{q_1}\right),\;
\sigma_{12} =  -RT\left(q_4 - \frac{q_2q_3}{q_1}\right).  
\end{eqnarray*}
Finally, the pressure is recovered through equation of state for ideal gases $p = \rho RT$. 

In the Galerkin-Boltzmann system of \eqref{eq:BoltzmannSystem}, physical quantities can be connected to the unknown numerical parameters through choosing the reference Mach number, \emph{Ma}, a free parameter that determines the compressibility of the fluid, and Reynolds number $Re$, a parameter which determines the ratio of inertial effects to viscous dissipation. These parameters are connected to the physical quantities via the relations
\begin{equation}
    \label{eq:SpeedofSound}
    \mathrm{\emph{Ma}} = \frac{U_r}{c}, \quad Re = \frac{U_r L_r}{\nu}, 
\end{equation}
where $U_r$ and $L_r$ are characteristic velocity and length, respectively. The value of $\tau$ follows directly via $\tau = \nu/RT$ from the choice of \emph{Ma} and $Re$. 

\subsection{An Unsplit Perfectly Matched Layer for the Galerkin-Boltzmann Equations}
\label{Sec:unplit PML}
The perfectly matched layer (PML) method requires introducing a finite width absorbing layer, called the PML, which surrounds the physical domain of interest so that waves leaving the domain and entering the PML are damped out. Suppose that the interface between the physical domain and the absorbing layer is aligned with the $x$ axis and is located at $x=0$ such that $x<0$ and $x>0$ correspond the physical domain and PML medium, respectively. The main idea of the PML approach is to construct a wave equation that admits a plane wave solution
\begin{equation*}
  \mathbf{q} = \mathbf{C}\exp\left(i\left(\mathbf{k}\cdot\mathbf{x}-wt\right) - \frac{k_x}{\omega}s(x)\right),  
\end{equation*}
where $\mathbf{C}$ is the polarization vector, $\mathbf{k} = [k_x,k_y]$ is the wave vector, and $s^x(x)$ is monotonic positive scalar function. The additional exponential term, $\exp(-\frac{k_x}{\omega}s^x(x))$ leads to an exponentially decaying wave amplitude in the increasing $x$ direction. Thus, the classical PML approach can be considered as an analytic continuation of PML medium in complex space having the following transformation
\begin{equation*}
  x \rightarrow x + \frac{i}{\omega}\int^x \sigma^x(r)\;dr,
\end{equation*}
where $s^x(x) = \int^x \sigma^x(r)\;dr$ and $\sigma^x$ is referred to as the damping profile, which is selected to be zero at interface between the physical domain and PML, and smoothly increasing through PML medium to avoid reflections. This transformation results in a new spatial differentiation operator in PML region
\begin{equation}
\label{Eq.ComplexTransform}
  \frac{\partial}{\partial \hat{x}} \rightarrow \frac{1}{1 + \frac{\sigma^x}{iw}}\frac{\partial }{\partial x},
\end{equation}
where $1 + \frac{\sigma^x}{iw}$ is called the coordinate stretching factor. If $\sigma^x = 0$ in the PML region, the transformation in \eqref{Eq.ComplexTransform} is reduced to simply the original physical coordinates.

Applying a Fourier transform in time, the Galerkin-Boltzmann equation given in (\ref{eq:BoltzmannSystem}) can be represented in the frequency domain as
\begin{equation}
\label{Eq:FourierBoltzmannSystem}
i\omega \hat{\mathbf{q}}= A_x\frac{\partial\hat{\mathbf{q}}}{\partial x} + A_y\frac{\partial\hat{\mathbf{q}}}{\partial y} + \widehat{\mathcal{N}(\mathbf{q})},
\end{equation}
where the hats are used to denote the time Fourier transformed fields. The PML equations are constructed by replacing the $x$ derivative operator via \eqref{Eq.ComplexTransform} with an analogous replacement for the $y$ derivative operator,
\begin{eqnarray*}
i\omega\hat{\mathbf{q}} = \frac{A_x}{1+\frac{\sigma^x}{i\omega}}\frac{\partial \hat{\mathbf{q}}}{\partial x} + \frac{A_y}{1+\frac{\sigma^y}{i \omega}}\frac{\partial \hat{\mathbf{q}}}{\partial y}  + \widehat{\mathcal{N}(\mathbf{q})}.
\end{eqnarray*}
Re-writing slightly, we obtain
\begin{eqnarray}
\label{Eq.FrequencyDomainEqs}
i\omega \hat{\mathbf{q}}= \left( A_x - \frac{A_x \sigma^x }{i\omega+\sigma^x}\right)\frac{\partial \hat{\mathbf{q}}}{\partial x} + \left( A_y - \frac{A_y \sigma^y }{i \omega + \sigma^y }\right)\frac{\partial \hat{\mathbf{q}}}{\partial y}  + \widehat{\mathcal{N}(\mathbf{q})}.
\end{eqnarray}
Next, we define two new variables in PML medium, i.e.,
\begin{equation}
\label{Eq.PMLvars}
\hat{\mathbf{q}}^x = \frac{1}{i\omega + \sigma^x}A_x\frac{\partial \hat{\mathbf{q}}}{\partial x}, \quad
\hat{\mathbf{q}}^y = \frac{1}{i\omega + \sigma^y}A_y\frac{\partial \hat{\mathbf{q}}}{\partial y},    
\end{equation}
and insert them into (\ref{Eq.FrequencyDomainEqs}). The PML equations in the frequency domain then take the following form,
\begin{eqnarray}
\label{Eq.PMLeqn}
i\omega \hat{\mathbf{q}}=  A_x\frac{\partial \hat{\mathbf{q}}}{\partial x} -\sigma^x \hat{\mathbf{q}}^x +  A_y \frac{\partial \hat{\mathbf{q}}}{\partial y} -\sigma^y \hat{\mathbf{q}}^y  + \widehat{\mathcal{N}(\mathbf{q})}.
\end{eqnarray}

Finally, we apply the inverse Fourier transform to both \eqref{Eq.PMLeqn} and \eqref{Eq.PMLvars} and transform back to the physical time domain which yields the unsplit equations 
\begin{eqnarray}
&\frac{\partial \mathbf{q}}{\partial t}&=  A_x \frac{\partial \mathbf{q}}{\partial x} - \sigma^x \mathbf{q}^x  +  A_y \frac{\partial \mathbf{q}}{\partial y} - \sigma^y \mathbf{q}^y+ \mathcal{N}(\mathbf{q}), \label{eq:PML_1}\\ 
&\frac{\partial \mathbf{q}^x}{\partial t} &=   - \sigma^x \mathbf{q}^x  +  A_x\frac{\partial \mathbf{q}}{\partial x},\label{eq:PML_2}\\
&\frac{\partial \mathbf{q}^y}{\partial t} &=  -  \sigma^y \mathbf{q}^y  + A_y\frac{\partial \mathbf{q}}{\partial y}. \label{eq:PML_3}
\end{eqnarray}

To avoid reflections, the damping profile is set to be zero at the interface of the physical domain and the PML, and smoothly increased across the PML width. In the corner regions, the damping profiles are taken as the superposition of the intersecting PML media. In the selection of PML profiles, we follow the M-PML formulation \citep{meza-fajardo_nonconvolutional_2008} to increase the damping performance and long term stability. M-PML introduces additional damping in the orthogonal directions as follows
\begin{equation}
    \begin{split}
     \sigma^x &= \hat{\sigma}^x(x) + \alpha^x \hat{\sigma}^y(y), \\
     \sigma^y &= \hat{\sigma}^y(y) + \alpha^y \hat{\sigma}^x(x), 
    \end{split}
\end{equation}
where $\hat{\sigma}^{x}$ and $\hat{\sigma}^y$ are the classical damping profiles for the regions having normal vectors parallel to $x$ and $y$, respectively, and $\alpha^x$ and $\alpha^y$ are constants that can be tuned for stability. With these multiaxial profiles, the M-PML applies additional damping the direction orthogonal to the usual PML dampening profile, which helps to damp shear waves generated in the PML region due to the relaxation term. One of the important advantages of using the unsplit PML formulation is that the nonlinear terms are not split in the PML. This allows us to directly implement the semi-analytic temporal discretizations, detailed in Section \ref{Sec:TimeDiscretization}, without any additional modifications.

\subsection{Nodal Discontinuous Galerkin Spatial Discretization}
\label{Sec:Spatial}
We assume that the domain $\Omega\in \mathbb{R}^{2}$ is well approximated by a computational domain, $\Omega_{h}$, which is partitioned into $K$ non-overlapping triangular elements, $\EN^e$, $e=1,\ldots,K$, such that, 
\begin{equation*}
    \Omega_h = \bigcup_{e=1}^K \EN^e.
\end{equation*}
Two elements, $\EN^{e+}$ and $\EN^{e-}$ are neighbours if they have a common face, that is $\dE^{e-} \cap \dE^{e+} \neq \emptyset$, where $\dE^{e}$ is the element boundary. We use $\mathbf{n} = \left(n_x, n_y \right)$ to denote the unit outward normal vector of $\dE$. 

We denote the approximation to $\mathbf{q}$, $\mathbf{q}^x$, and $\mathbf{q}^y$ on element $\EN^e$ as $\mathbf{q}^e$, $\mathbf{q}^{x,e}$, and $\mathbf{q}^{y,e}$, respectively. The local trace values of $\mathbf{q}^e$ along $\dE^{e}$ are denoted as $\mathbf{q}^-$ and the corresponding neighboring trace values are denoted using $\mathbf{q}^+$, omitting the $e$ superscript when it is clear which element has the local trace. We define $\avg{\mathbf{q}}$ and $\jump{\mathbf{q}}$ to be the average and jump of $\mathbf{q}^e$ along the the trace $\dE^e$, i.e.,
\begin{equation}\label{Eq.AverageJumpScalar}
\avg{\mathbf{q}} = \frac{\mathbf{q}^{+} + \mathbf{q}^{-}}{2}, \quad  \jump{\mathbf{q}} = \mathbf{q}^{+} -\mathbf{q}^{-}.
\end{equation}

Each element of triangulation $\Omega_h$ is an affine mapping of the reference bi-unit triangular element, $\hat{\EN} =  \left\{ -1\leq r,s,r+s\leq 1 \right\}$ under the map, $\Phi^e$, given by  
\begin{equation*}
    \label{eq:operators_1}
    \left(x,y\right) =\Phi^e\left(r,s\right), \quad \left(x,y\right)\in \EN^e,\ \left(r,s\right)\in\hat{\EN}.
\end{equation*}
The Jacobian of this mapping can be written
\begin{equation*}
    \label{eq:operators_2}
    G^e = \begin{bmatrix}
    x^e_r & x^e_s \\
    y^e_r & y^e_s 
    \end{bmatrix},
\end{equation*}
and we denote its determinant as $J^e = \det G^e$. We also define the surface scaling factor $J^{ef}$ as the determinant of the Jacobian $G^e$ restricted to the face, $\dE^{ef}$.

Let $\EN^e$ be an element. We select the finite element space $V_N^e$ to be $\mathcal{P}_N(\EN^e)$, the space of polynomial functions of degree $N$ on this element. For a choice of basis, we use $N_p = |V_n^e|$ Lagrange polynomials interpolating at the Warp \& Blend nodes \citep{warburton2006explicit} mapped to the element $\EN^e$, which we denote $\{\phi^e_i\}_{i=1}^{N_p}$. 

We selected the unsplit PML equation for $\mathbf{q}$ i.e., \eqref{eq:PML_1} in order to describe the spatial discretization of the equation system, because this equation includes all terms for the discretization of the PML equations and recovers the physical domain equations, \eqref{eq:BoltzmannSystem} for vanishing $\sigma^x$ and $\sigma^y$. Multiplying \eqref{eq:PML_1} by a test function $v \in V^e_N$, integrating over the element $\EN^e$, and performing integration by parts twice, we arrive to the following strong variational form to be solved, 
\begin{equation}
\label{eq:spatial_1}
\begin{split}
\int_{\EN^e} v \frac{\partial \mathbf{q}^e}{\partial t} =& \int_{\EN^e} v \left(A_x\frac{\partial \mathbf{q}^e}{\partial x}+A_y\frac{\partial \mathbf{q}^e}{\partial y}-\sigma^x \mathbf{q}^{x,e}-\sigma^y \mathbf{q}^{y,e}\right) \\ &+\int_{\dE^e} \phi F\left(\mathbf{q}^{*}-\mathbf{q}^{-}\right) + \int_{\EN^e} \phi\mathcal{N}(\mathbf{q}^e).
\end{split}
\end{equation}
where $F = n_x A_x+n_y A_y$ is the flux matrix in the direction of the element normal vector $\mathbf{n}$ and $\mathbf{q}^{e*}$ is a trace state defined using an upwind numerical flux function, which depends on the local and neighboring traces values along $\partial\EN^e$. The upwind flux is can be formulated by diagonalizing the operator $F$ as $F =  \mathcal{R} \Lambda \mathcal{R}^{-1}$. Because the transport terms of the Galerkin-Boltzmann equations are purely hyperbolic, the diagonal matrix, $\Lambda$ has only real entries of $0,0,\pm c, \pm c\sqrt{3}$. Splitting the eigenvalues that have positive signs $\Lambda^{+}$ and negative signs $\Lambda^{-}$, the upwind flux can then be written as,
\begin{equation*}
 \label{eq:spatial_2}
 F \mathbf{q}^* = \mathcal{R}\left( \Lambda^{+} \mathcal{R}^{-1} \mathbf{q}^- + \Lambda^{-} \mathcal{R}^{-1} \mathbf{q}^+ \right).
\end{equation*}

To evaluate the integrals involving $\sigma$ or the nonlinear term $\mathcal{N}(\mathbf{q}^e)$ in \eqref{eq:spatial_1}, we use a sufficiently high-order cubature rule to reduce aliasing errors. The cubature-based integration using an interpolation operator which interpolates the solution field to cubature nodes on each element. We select a nodal set of $N_c$ cubature nodes with coordinates $(r^c_i,s^c_i)$ for $i=1,\ldots,N_c$ on the reference element $\hat{\EN}$, and associated weights, $w^c_i$, for $i=1,\ldots,N_c$. We then define a set of cubature nodes $(x^{e,c}_i,y^{e,c}_i)$ for $i=1,\ldots,N_c$ on each element $\EN^e$ to be the cubature nodes on the reference element mapped to $\EN^e$ via $\Phi^e$. The interpolation operator, $\mathcal{I}^e$ can then be defined as follows
\begin{equation*}
    \mathcal{I}^e_{ij} = \phi^e_i(x_j^{e,c}, y_j^{e,c}),
\end{equation*}
for $j=1,\ldots,N_c$ and $i=1,\ldots,N_p$. We also define on each element, we define mass, surface mass, and stiffness operators as follows
\begin{align*}
    \mathcal{M}_{ij}^e = \int_{\EN^e}\phi^e_j \phi^e_i, &\quad 
    \mathcal{M}_{ij}^{ef} = \int_{\partial\EN^{e}}\phi^e_j \phi^e_i, \\
    (\mathcal{S}^e_x)_{ij} = \int_{\EN^e}\phi^e_j \frac{\partial\phi^e_i}{\partial x}, &\quad
    (\mathcal{S}^e_y)_{ij} = \int_{\EN^e}\phi^e_j \frac{\partial\phi^e_i}{\partial y}, 
\end{align*}
respectively. Then, selecting the test function to be a basis function, i.e. $v=\phi_i$ and writing the nodal values of $\mathbf{q}^e$, $\mathbf{q}^{x,e}$, and $\mathbf{q}^{y,e}$ as $\mathbf{q}^e_i$, $\mathbf{q}^{x,e}_i$, and $\mathbf{q}^{y,e}_i$, respectively, for $i=1,\ldots,N_p$ we obtain that \eqref{eq:spatial_1} can be written 
\begin{equation}
\label{eq:spatial_1_1}
\begin{split}
\mathcal{M}_{ij}^e\frac{\partial \mathbf{q}^e_j}{\partial t} =& A_x (\mathcal{S}^e_x)_{ij} \mathbf{q}^e_j +A_y(\mathcal{S}^e_x)_{ij} \mathbf{q}^e_j - J^e\mathcal{I}^e_{ki} w_k\sigma^x_k \mathcal{I}^e_{kj}\mathbf{q}^{x,e}_j - J^e\mathcal{I}^e_{ki} w_k \sigma^y_k \mathcal{I}^e_{kj}\mathbf{q}^{y,e}_j \\ &+ \mathcal{M}_{ij}^{ef} (F\left(\mathbf{q}^{*}-\mathbf{q}^{-}\right))_j + J^e\mathcal{I}^e_{ki} w_k \mathcal{N}(\mathcal{I}^e_{kj}\mathbf{q}^e_j),
\end{split}
\end{equation}
where we have made use of Einstein repeated index summation notation for $j = 1,\ldots,N_p$, $k=1,\ldots, N_c$, and $f=1,\ldots,N_f$ where $N_f$ is the number of faces per element. Here $\sigma^x_k$ and $\sigma^y_k$ are the PML damping profiles evaluated at the cubature point $(x^{e,c}_k,y^{e,c}_k)$. Upon multiplying \eqref{eq:spatial_1_1} by $(\mathcal{M}^e)^{-1}$, we define the differentiation, lift, and cubature projection operators as
\begin{align*}
    \mathcal{D}^e_x = (\mathcal{M}^e)^{-1}\mathcal{S}^e_x, &\quad 
    \mathcal{D}^e_y = (\mathcal{M}^e)^{-1}\mathcal{S}^e_y, \\
    \mathcal{L}^{ef} = (\mathcal{M}^e)^{-1}\mathcal{M}^{ef}, &\quad
    \mathcal{P}^e = (\mathcal{M}^e)^{-1}(\mathcal{I}^{e})^T \mathrm{diag}(w),
\end{align*}
respectively, where $\mathrm{diag}(w)$ is a diagonal matrix with entries $w_i$ for $i=1,\ldots,N_c$, we can write \eqref{eq:spatial_1_1} as
\begin{equation}
\label{eq:spatial_3}
\begin{split}
\frac{\partial \mathbf{q}^e_i}{\partial t} =& A_x (\mathcal{D}^e_x)_{ij} \mathbf{q}^e_j +A_y(\mathcal{D}^e_y)_{ij} \mathbf{q}^e_j - J^e\mathcal{P}^e_{ik} \sigma^x_k \mathcal{I}^e_{kj}\mathbf{q}^{x,e}_j - J^e\mathcal{P}^e_{ik} \sigma^y_k \mathcal{I}^e_{kj}\mathbf{q}^{y,e}_j \\ &+ \mathcal{L}_{ij}^{ef} (F\left(\mathbf{q}^{*}-\mathbf{q}^{-}\right))_j + J^e\mathcal{P}^e_{ik} \mathcal{N}(\mathcal{I}^e_{kj}\mathbf{q}^e_j),
\end{split}
\end{equation}

Finally, since we assume all elements are images under an affine mapping of the reference element $\hat{\EN}$, the nodal DG spatial discretization \eqref{eq:spatial_3} can be expressed simply in terms of reference differentiation matrices $\mathcal{D}_r$ and $\mathcal{D}_s$, lift matrices $\mathcal{L}^f$, interpolation $\mathcal{I}$, and projection $\mathcal{P}$ defined on the reference element $\hat{\EN}$ through the geometric factors of $\Phi^e$ via
\begin{align*}
    \mathcal{D}^e_x &= r^e_x\mathcal{D}_r + s^e_x\mathcal{D}_s, \\
    \mathcal{D}^e_y &= r^e_y\mathcal{D}_r + s^e_y\mathcal{D}_s, \\
    \mathcal{L}^{ef} &= \frac{J^{ef}}{J^e}\mathcal{L}^{f},  \\
    \mathcal{I}^e &= \mathcal{I}, \\
    \mathcal{P}^e &= \frac{1}{J^e}\mathcal{P},
\end{align*}
where $r^e_x, r^e_y, s^e_x,$ and $s^e_y$ are defined via
\begin{equation*}
    (G^e)^{-1} = \begin{bmatrix}
    r^e_x & r^e_y \\
    s^e_x & s^e_y 
    \end{bmatrix}.
\end{equation*}
Using these reference operators we write in the semi-discrete scheme \eqref{eq:spatial_3} as
\begin{equation}
\label{eq:spatial_9}
\begin{split}
\frac{\partial \mathbf{q}^e_i}{\partial t} =& A^e_r (\mathcal{D}_r)_{ij} \mathbf{q}^e_j +A^e_s(\mathcal{D}_s)_{ij} \mathbf{q}^e_j - \mathcal{P}_{ik} \sigma^x_k \mathcal{I}_{kj}\mathbf{q}^{x,e}_j - \mathcal{P}_{ik} \sigma^y_k \mathcal{I}_{kj}\mathbf{q}^{y,e}_j \\ &+ \frac{J^{ef}}{J^e}\mathcal{L}_{ij}^{f} (F\left(\mathbf{q}^{*}-\mathbf{q}^{-}\right))_j + \mathcal{P}_{ik} \mathcal{N}(\mathcal{I}_{kj}\mathbf{q}^e_j),
\end{split}
\end{equation}
where
\begin{eqnarray*}
A^e_r & = & r^e_x A_x + r^e_y A_y, \\
A^e_s & = & s^e_x A_x + s^e_y A_y.
\end{eqnarray*}

The operators present in \eqref{eq:spatial_9} describes all the actions required to solve our PML formulation in the full system \eqref{eq:PML_1}-\eqref{eq:PML_3}. Repeating the procedure above, we obtain the semi discrete forms of \eqref{eq:PML_2} and \eqref{eq:PML_3} as
\begin{align}
\frac{\partial \mathbf{q}^{x,e}_i}{\partial t} =& r^e_x A_x (\mathcal{D}_r)_{ij} \mathbf{q}^{e}_j +s^e_x A_x(\mathcal{D}_s)_{ij} \mathbf{q}^e_j \nonumber\\ & \qquad- \mathcal{P}_{ik} \sigma^x_k \mathcal{I}_{kj}\mathbf{q}^{x,e}_j + \frac{J^{ef}}{J^e}\mathcal{L}_{ij}^{f} (n_x A_x\left(\mathbf{q}^{*}-\mathbf{q}^{-}\right))_j, \label{eq:spatial_10}\\
\frac{\partial \mathbf{q}^{y,e}_i}{\partial t} =& r^e_y A_y (\mathcal{D}_r)_{ij} \mathbf{q}^{e}_j +s^e_y A_y(\mathcal{D}_s)_{ij} \mathbf{q}^e_j \nonumber\\ & \qquad - \mathcal{P}_{ik} \sigma^y_k \mathcal{I}_{kj}\mathbf{q}^{y,e}_j + \frac{J^{ef}}{J^e}\mathcal{L}_{ij}^{f} (n_y A_y\left(\mathbf{q}^{*}-\mathbf{q}^{-}\right))_j.\label{eq:spatial_11}
\end{align}
For the vanishing $\sigma^x$ and $\sigma^y$, \eqref{eq:spatial_9} also gives all the required operators for the semi-discrete form of the system in the non-PML region. In the next section, we cover the semi-analytic and implicit-explicit time discretizations using the semi-discrete equation \eqref{eq:spatial_9}.

\section{Time Discretization}
\label{Sec:TimeDiscretization}
The nonlinear collision term in the Galerkin Boltzmann equation becomes  stiff in the limit of small relaxation times $(\tau<<1)$ which introduces a severe time step restriction if a fully explicit time integrator is used.  In this section, we discuss two different temporal integration methods for the Boltzmann equation in stiff regimes: a semi-analytic time discretization method, also called exponential time discretization, and a low storage implicit-explicit Runge-Kutta method.

Assembling the semi-discrete system in \eqref{eq:spatial_9}-\eqref{eq:spatial_11} on each element $\EN^{e}$ into a global system, we arrive to the following problem,
\begin{equation}
\label{eq:time_1}
\frac{d\mathbf{q}}{dt}= \mathbf{L}(\mathbf{q}) + \mathbf{N}(\mathbf{q}),
\end{equation}
where $\mathbf{q}$ denotes assembled global vector of degrees of freedom,  $\mathbf{L}$ collects all the linear terms, and $\mathbf{N}$ includes the relaxation terms. We also use $\mathbf{N}$ and $\mathbf{L}$ for stiff and non-stiff terms depending on the coefficient $\tau$.

Due to the special structure of the nonlinear term, stiffness affects only the last three equations of the system. The time step restriction for the first three equations is thus always the advective time scale and it is independent of the $\tau$ term. We split the equation system \ref{eq:time_1} in two parts in a way that last three equations advanced with the specific time integration methods for stiff problems and first three equations always integrated with explicit time stepper that semi-analytic or implicit-explicit method reduce in the limiting case, $1/\tau \to 0$.  

To derive a semi-analytic explicit time discretization, we note that from the form of the nonlinear term in \eqref{eq:nonlinearDef} we can write \eqref{eq:time_1} as
\begin{equation}
\label{eq:time_3}
\frac{d\mathbf{q}}{dt}= -\boldsymbol\Lambda \mathbf{q} + \mathbf{L}(\mathbf{q})  + \mathbf{\tilde{N}}(\mathbf{q})
\end{equation}
where $\boldsymbol\Lambda = \mathrm{diag}\left(0,0,0,\frac{1}{\tau},\frac{1}{\tau},\frac{1}{\tau}\right)$ and   $\mathbf{\tilde{N}}(\mathbf{q})=\left(0,0,0,\frac{q_2 q_3}{\tau q_1}, \frac{q_2^2}{\tau q_1\sqrt{2}}, \frac{q_3^2}{\tau q_1\sqrt{2}}\right)^T$. Note that $\mathbf{\tilde{N}}(\mathbf{q})$ now does not depend on $q_4, q_5,$ nor $q_6$. Finally, we define $\mathbf{F}(\mathbf{q}) = \mathbf{L}(\mathbf{q})  + \mathbf{\tilde{N}}(\mathbf{q})$ to simplify the notation in the derivation of semi-analytic time discretization methods. With these modifications, we obtain 
\begin{equation}
\label{eq:time_4}
\frac{d\mathbf{q}}{dt}= -\boldsymbol\Lambda\mathbf{q} + \mathbf{F}(\mathbf{q}).
\end{equation}
Multiplying \eqref{eq:time_4} by $e^{\boldsymbol\Lambda t}$ and integrating from $t_n$ to $t_{n+1}$, we obtain the following Voltera integral equation,
\begin{equation}
\label{eq:ExponentialExact}
\mathbf{q}(t_{n+1})= \mathbf{q}(t_n) e^{-\boldsymbol\Lambda (t_{n+1}-t_n)} + \int_{t_n}^t e^{\boldsymbol\Lambda (\theta-t_{n+1}) }\mathbf{F}\left(\mathbf{q}\left(\theta\right), \theta \right)d\theta.
\end{equation}
Note that since the first three rows of $\boldsymbol\Lambda$ are zero, the first three equations of \eqref{eq:ExponentialExact} are simply the first three equations of \eqref{eq:time_1}
integrated in time. 

The derived formula is exact and we use it to derive semi-analytic time integration methods. To simplify the notation in the following sections, we denote the numerical approximation to $\mathbf{q}(t_n)$ by $\mathbf{q}_n$ and we shorten $\mathbf{F}(\left(\mathbf{q}\left(t_n\right), t_n \right)$ to $\mathbf{F}_{n}$. We denote the $i^{\text{th}}$ history field of $\mathbf{q}$ and $\mathbf{F}$ at given discrete time level $t_{n-i}$ by $\mathbf{q}_{n-i}$ and  $\mathbf{F}_{n-i}$, respectively. Similarly, we use $\mathbf{q}_{n,i}$  and  $\mathbf{F}_{n,i}$ denote the corresponding states at the same intermediate stage times.  

\subsection{Semi-analytic  Multistep Methods}
In this section, we present formal derivation of semi-analytic Adams-Bashforth (SAAB) methods and we extend the idea to multirate semi-analytic Adams-Bashforth methods (MRSAAB) with different level difference between the groups. 

We start with the basic structure of a linear multistep method, which is the polynomial extrapolation of the integration function given in \eqref{eq:ExponentialExact} from the arbitrary order $s$ with an extrapolation function, $P_s(\theta)$. This leads to a scheme of the form,
\begin{equation}
\label{eq:SAAB_1}
\mathbf{q}(t_{n+1})= \mathbf{q}_n e^{-\boldsymbol\Lambda (t_{n+1}-t_n)} + \int_{t_n}^{t_{n+1}} e^{\boldsymbol\Lambda(\theta-t_{n+1})} P_s(\theta)d\theta,
\end{equation}
where $P_s(\theta)$ is extrapolated from $s$ sampling points of $\mathbf{F}\left(\mathbf{q}( \theta), \theta \right)$. In addition, the following property holds,
\begin{equation}
\label{eq:SAAB_2}
P_{s} (t_{n-i}) = \mathbf{F}\left(\mathbf{q}( t_{n-i}), t_{n-i} \right),
\end{equation}
for $i=0,\ldots,(s-1)$. To construct $P_s$, we use the classical Lagrange interpolating polynomials as, 
\begin{eqnarray*}
\label{eq:SAAB_3}
l_j(t) = \prod_{i=0, i\neq j}^{s-1} \frac{t - t_{n-i}}{ t_{n-j}  - t_{n-i}},   
\end{eqnarray*}
After rewriting the second term of the \eqref{eq:SAAB_1} in terms of Lagrange interpolating polynomials we obtain, 
\begin{equation}
\label{eq:SAAB_5}
\int_{t_n}^{t_{n+1}} P_s(\theta)d\theta = \int_{t_n}^{t_{n+1}} e^{\boldsymbol\Lambda(\theta-t_{n+1}) }\mathbf{F}_{n-i} l_i(\theta) d\theta,
\end{equation}
for $i=0,\ldots,(s-1)$. Selecting a uniform step size of $\Delta t$ so that $t_{n+1}=t_n + \Delta t$, we arrive at the multistep semi-analytic Adams-Bashforth method
\begin{equation*}
\label{eq:SAAB_5_1}
\mathbf{q}(t_{n+1})= \mathbf{q}_n e^{-\boldsymbol\Lambda \Delta t} + \Delta t\sum_{i=0}^{s} \tilde{a}_i \mathbf{F}_{n-i},
\end{equation*}
where $\tilde{a}_i$ are coefficients that can be computed analytically using,
\begin{align}
\tilde{a}_i \Delta t &= \int_{t_n}^{t_{n}+\Delta t} e^{\boldsymbol\Lambda(\theta-t_{n+1}) } l_i(\theta) d\theta, \nonumber\\
&= \int_{0}^{\Delta t} e^{\boldsymbol\Lambda(\theta-\Delta t) } l_i(\theta-t_n) d\theta,\label{eq:SAAB_5_2}
\end{align}
for $i=0,\ldots,(s-1)$. Note that from the definition of $l_i$ in \eqref{eq:SAAB_3} and the assumption of a uniform time step size, $l_i(\theta-t_n)$ can be expressed in terms of $\Delta t$ and $\theta$ only. 

Since $ e^{-\boldsymbol\Lambda t} = \mathrm{diag}(1,1,1, e^{-\frac{t}{\tau}}, e^{-\frac{t}{\tau}}, e^{-\frac{t}{\tau}})$, the first three components of \eqref{eq:SAAB_5_1} contain no exponential terms and the coefficients $\tilde{a}_i$ simply reduce to the coefficients of the classical Adams-Bashforth methods, which we denote $a_i$. The coefficients for the last three equations are modified to include integration with the exponential factor $e^{-\frac{t}{\tau}}$. For $s=3$, these semi-analytic coefficients can be written
\begin{equation}
\label{eq:SAAB_6}
\begin{split}
\tilde{a}_0 &= \gamma^{-3}\left[\left(-1-\frac{5}{2}\gamma - 3\gamma^2\right) -  e^{\gamma}
\left( -1 - \frac{3}{2}\gamma - \gamma^2 \right) \right],  \\
\tilde{a}_1 &= \gamma^{-3}\left[\left(2+4\gamma + 3\gamma^2\right) -  e^{\gamma}
\left(2 + 2\gamma \right) \right],  \\
\tilde{a}_2 &= \gamma^{-3}\left[\left(-1-\frac{3}{2}\gamma - \gamma^2\right) -  e^{\gamma}
\left(-1 - \frac{1}{2}\gamma\right) \right]. 
\end{split}
\end{equation}
where $\gamma = -\frac{\Delta t}{\tau}$.

In the formal limiting case, $\frac{1}{\tau} \to 0$, the SAAB coefficients become the classical third order Adams-Bashforth coefficients, i.e. $\lim_{\frac{\Delta t}{\tau}\to0}\tilde{a}_i = a_i$, where $a_0 = 23/12$, $a_1 = -16/12$, and $a_2 = 5/12$. Semi-analytic schemes with arbitrary order have been derived elsewhere \citep{cox_exponential_2002}, but we include an explicit expression and derivation in preparation for a multirate version that allows elements to make different time steps. Although we discuss multirate time stepping methods in detail below, we include here a brief overview in order to include the necessary multirate SAAB coefficients in this section for completeness. We obtain the coefficients required to perform a fractional time-step in the SAAB method by setting $t_{n+1} = t_n + \Delta t/2$ in \eqref{eq:SAAB_5} and repeating the process in \eqref{eq:SAAB_5_1}-\eqref{eq:SAAB_5_2} to obtain that the multirate semi-analytic coefficients can be written as,
\begin{equation}
\label{eq:MRSAAB_7}
\tilde{b}_i \Delta t = \int_0^{\Delta t/2} e^{\boldsymbol\Lambda(\theta-\Delta t/2) }l_i(\theta-t_n) d\theta.
\end{equation}
For the order $s=3$ method gives the MRSAAB coefficients required for a fractional step of $\Delta t/2$ to be either the classical multirate Adams-Bashforth coefficients, i.e. $b_0 = 17/24$, $b_1 = -7/24$, and $b_2 = 1/12$, or the modified coefficients,
\begin{equation}
\label{eq:MRSAAB_8}
 \begin{split}
\tilde{b}_0 &=\gamma^{-3}\left[\left(-1 -2\gamma - \frac{15}{8}\gamma^2\right) -  e^{\frac{\gamma}{2}}
\left( -1 - \frac{3}{2}\gamma - \gamma^2 \right) \right],  \\
\tilde{b}_1 &=\gamma^{-3}\left[\left(2 +3\gamma + \frac{5}{4}\gamma^2\right) -  e^{\frac{\gamma}{2}}
\left( 2 + 2\gamma \right) \right],  \\
\tilde{b}_2 &=\gamma^{-3}\left[\left(-1 -1\gamma - \frac{3}{8}\gamma^2\right) -  e^{\frac{\gamma}{2}}
\left( -1 - \frac{1}{2}\gamma \right) \right], 
\end{split}   
\end{equation}

Adams-Bashforth methods are not self starting and need $s$ initial history fields to start extrapolation. Although both MRSAAB and SAAB histories can be computed with a self-starting time discretization methods, we also present the first and second order SAAB and MRSAAB coefficients which can be used to initialize the time-stepping method and prevent additional algorithmic complexity in start-up. For $s=1$ the coefficients can be obtained by a similar procedure described above to obtain
\begin{eqnarray}
\label{eq:MRSAAB_9}
\tilde{a}_0 = \gamma^{-1}\left[ e^{\gamma} -1 \right],\quad 
\tilde{b}_0 = \gamma^{-1}\left[ e^{\frac{\gamma}{2}} -1 \right], 
\end{eqnarray}
with classical Adams-Bashforth coefficients being $a_0 = 1$ and $b_0 = 1/2$. Then, for the second order $s=2$ methods,  
\begin{equation}
\label{eq:SAAB_10}
\begin{split}
\tilde{a}_0 &= \gamma^{-2}\left[\left(-1 -2\gamma\right) -  e^{\gamma}
\left( -1 - \gamma\right) \right],  \\
\tilde{a}_1 &=\gamma^{-2}\left[\left(1 +\gamma\right) +  e^{\gamma} \right],  \\
\tilde{b}_0 &=\gamma^{-2}\left[\left(-1 -\frac{3}{2}\gamma\right) -  e^{\frac{\gamma}{2}}
\left( -1 - \gamma \right) \right],  \\
\tilde{b}_1 &=\gamma^{-2}\left[\left(1 +\frac{1}{2}\gamma\right) -  e^{\frac{\gamma}{2}}\right], 
\end{split}
\end{equation}
and the single rate and multirate classical AB coefficients are $a_0 = 3/2$, $a_1 = -1/2$ and $b_0 = 5/8$, $b_1 = -1/8$, respectively. It is worthwhile to mention that exponential and non-exponential integrated parts of the Galerkin-Boltzmann equation are consistent in time when all coefficients converge to a classical multirate or a single rate Adams-Bashforth methods in the non-stiff limit. 

\subsection{Semi-analytic Runge-Kutta Methods}
Runge-Kutta methods can be constructed analogously to multistep methods. Let us begin by integrating the equation \eqref{eq:ExponentialExact} from $t=t_n$ to some intermediate time level $t=t_n+\Delta t_i$ which leads to a variation-of-constants formula,
\begin{equation*}
\mathbf{q}_{ni}= \mathbf{q}_n e^{-\boldsymbol\Lambda \Delta t_i} + \int_0^{\Delta t_i} e^{\boldsymbol\Lambda(\theta-\Delta t_i) }\mathbf{F}\left(\mathbf{q}( t_n+\theta) , t_n + \theta\right)d\theta.
\end{equation*}
For general one-step methods, the internal and the final stages are approximated in the following way,
\begin{equation*}
\label{Eq:ERK_1}
    \begin{split}
    \mathbf{q}_{ni}&= \mathbf{q}_n e^{-\boldsymbol\Lambda \Delta t_i} + \Delta t\sum_{j=0}^{s-1}\tilde{a}_{ij}\mathbf{F}\left(\mathbf{q}( t_n+\Delta t_j) , t_n + \Delta t_j\right) = \mathbf{q}_n e^{-\boldsymbol\Lambda \Delta t_i} + \Delta t\sum_{j=0}^{s-1}\tilde{a}_{ij}\mathbf{F}_{nj}, \\
\mathbf{q}_{n+1}&= \mathbf{q}_n e^{-\boldsymbol\Lambda \Delta t} + \Delta t\sum_{i=0}^{s-1}\tilde{b}_{i}\mathbf{F}\left(\mathbf{q}( t_n+\Delta t_i) , t_n + \Delta t_i\right) = \mathbf{q}_n e^{-\boldsymbol\Lambda \Delta t} + \Delta t\sum_{j=0}^{s}\tilde{b}_{i}\mathbf{F}_{ni},
    \end{split}
\end{equation*}
where $s$ is the number of stages, $\tilde{a}$ and $\tilde{b}$ are the semi analytic Runge-Kutta (SARK) method coefficients computed using exponential functions or some approximation of exponential functions. We assume that all methods satisfy $\Delta t_1 = c_1 \Delta t=0$ leading to $e^{\boldsymbol\Lambda \Delta t_1} = I$ for consistency reasons. Similar to the derivation of the semi-analytic Adams-Bashforth methods above, a semi-analytic method reduces to the base Runge-Kutta method in the limit $\frac{1}{\tau} \to 0$ which makes exponential and non-exponential parts of the integrated equation consistent. We also assume that base Runge-Kutta method satisfies,
\begin{equation}
\label{Eq:RK_general_constraint}
 \sum_{j=0}^{s-1}b_j = 1, \quad \sum_{j=0}^{s-1}a_{ij} = c_i. 
\end{equation}
The semi-analytic Runge-Kutta time discretization satisfies an analogous constraint,
\begin{equation}
\label{Eq:SARK_general_constraint}
 \sum_{j=0}^{s-1}\tilde{b}_j    = \gamma^{-1}\left(e^{\gamma}-1\right), \quad \sum_{j=0}^{s-1}\tilde{a}_{ij} = \frac{1}{c_i}\gamma^{-1}\left(e^{c_i\gamma}-1\right),
\end{equation}
for $i=0,\ldots,(s-1)$. We introduce a class of third order SARK schemes by modifying the base method coefficients. The internal stages  are computed using following relation,
\begin{equation*}
\tilde{a}_{ij} = \frac{1}{c_i}\gamma^{-1}\left(e^{c_i\gamma-1}\right) {a}_{ij},
\end{equation*}
which directly satisfies \eqref{Eq:SARK_general_constraint} if the base method satisfies \eqref{Eq:RK_general_constraint}. The final update stage is then computed using Lagrange interpolation of the function values at internal stages assuming the non-repeating internal stage time levels as given below,
\begin{equation}
\label{Eq:ERK_2}
\mathbf{q}_{n+1}= \mathbf{q}_n e^{-\boldsymbol\Lambda \Delta t} + \Delta t\sum_{j=0}^{s-1}\left(\int_0^{\Delta t} e^{\boldsymbol\Lambda( \theta-\Delta t)}l_i(\theta-t_n) d\theta\right) \mathbf{F}_{n,i},
\end{equation}
where $l_i$ are again the Lagrange interpolating polynomials, this time constructed as interpolating at the intermediate stage times, i.e.  
\begin{equation}
\label{Eq:ERK_3}
l_i(t) = \prod_{j=0, i\neq j}^{s-1} \frac{t - t_n - c_j \Delta t}{ \Delta t(c_i  - c_j)},  
\end{equation}
for $i=0,\ldots,(s-1)$.
\bgroup
\def\arraystretch{1.5}%
\begin{table}[t]
\caption{Butcher tableaus for the classical third-order  and adapted (b) method, based on RK2a, with coefficients }
    \begin{subtable}{.45\linewidth}
    \centering
    \setlength{\tabcolsep}{5pt}
        \caption{}
        \begin{tabular}{c| c c c}
            0             &                &  & \\
            $\frac{1}{2}$ & $\frac{1}{2}$  &  & \\
            $1$           & $-1$  & $2$ & \\ \hline
              & $\frac{1}{6}$ & $\frac{2}{3}$ &$\frac{1}{6}$
        \end{tabular}
    \end{subtable}
     \begin{subtable}{.45\linewidth}
    \setlength{\tabcolsep}{5pt}
      \centering
        \caption{}
        \begin{tabular}{c| c c c}
            0             &                &  & \\
            $\frac{1}{3}$ & $\frac{1}{3}$  &  & \\
            $\frac{3}{4}$ & $-\frac{3}{16}$  & $\frac{15}{16}$ & \\ \hline
              & $\frac{1}{6}$ & $\frac{3}{10}$ &$\frac{8}{15}$
        \end{tabular}
    \end{subtable}%
    \label{table:Butcher}
\end{table}
\egroup
If we start with the classical third-order RK method with the Butcher tableau given in Table \ref{table:Butcher}(a), the coefficients of the SARK method are
\begin{equation}
    \label{eq:ERK_4}
    \begin{split}
        \tilde{a}_{10} &= \gamma^{-1} \left[-1 + e^{\frac{\gamma}{2}}\right], \\
        \tilde{a}_{20} &= \gamma^{-1} \left[ 1 - e^{\gamma}\right], \\
        \tilde{a}_{21} &= \gamma^{-1} \left[-2 + 2e^{\gamma}\right], \\
        \tilde{b}_{0}  &= \gamma^{-3}\left[ -4 - \gamma - e^{\gamma}\left(-4 + 3\gamma - \gamma^2 \right) \right],\\
        \tilde{b}_{1}  &= \gamma^{-3}\left[ 8 + 4\gamma - e^{\gamma}\left(8 - 4\gamma\right) \right],\\
        \tilde{b}_{2}  &= \gamma^{-3}\left[ -4 - 3\gamma - \gamma^2 - e^{\gamma}\left(-4 + \gamma\right) \right].
    \end{split}
\end{equation} 
This SARK method reproduces the results for the third-order exponential RK scheme reported in \citep{cox_exponential_2002} with more straightforward derivation. By following the same procedure, we construct a third-order SARK method with better truncation errors by using the base RK method  where its Butcher tableau given in Table \ref{table:Butcher}(b). The coefficients of the SARK scheme used in this study are given below,
%
\begin{equation}
    \label{eq:ERK_5}
    \begin{split}
        \tilde{a}_{10} &= \gamma^{-1} \left[-1 + e^{\frac{\gamma}{3}}\right], \\
        \tilde{a}_{20} &= \frac{1}{4}\gamma^{-1} \left[ 1 - e^{\frac{3\gamma}{4}}\right], \\
        \tilde{a}_{21} &= \frac{1}{4}\gamma^{-1} \left[-5 + 5e^{\frac{3\gamma}{4}}\right], \\
        \tilde{b}_{0}  &= \gamma^{-3}\left[ -24 - 11\gamma - 2\gamma^2 - e^{\gamma}\left(-24 + 13\gamma - 3\gamma^2 \right) \right],\\
        \tilde{b}_{1}  &= \frac{36}{5}\gamma^{-3}\left[ 2 + \frac{5}{4}\gamma + \frac{1}{4}\gamma^2 - e^{\gamma}\left(2 - \frac{3}{4}\gamma\right) \right],\\
        \tilde{b}_{2}  &= \frac{16}{5}\gamma^{-3}\left[ -2 - \frac{5}{3}\gamma - 2\gamma^2 - e^{\gamma}\left(-2 + \frac{1}{3}\gamma\right) \right].
    \end{split}
\end{equation} 
All semi-analytic time stepping coefficients include a terms similar to 
\[
f(z) = \frac{\exp(z) -1}{z},
\]
or some high-order polynomial variant of this expression. The accuracy of SARK and SAAB methods depend strongly on the accurate computation of this function. For small $z$, direct computation of such expressions can encounter large cancellation errors. Therefore, in the small $z$ limit a Taylor series approximation can be better option. On the other hand, Taylor series approximation is inaccurate if $z$ is large. To obtain a more robust approach we compute the coefficients using a complex contour integral \citep{kassam_fourth-order_2005}. For example, the evaluation of $f(z)$ is done by integrating over the contour $\Gamma$ in the complex plane enclosing $z$ as follows
\begin{equation*}
    f(z) = \frac{1}{2\pi i} \int_{\Gamma} \frac{f(\theta)}{\theta-z}d\theta.
\end{equation*}
In our numerical tests, we take $\Gamma$ as the unit circle sampled with 64 equally spaced points. Due to symmetry, integration only requires $32$ points on the upper plane. We compute function values at these points and we take the mean of the real part of function values. Using $32$ points on half plane gives full accuracy of coefficients (in 14 digits) independent of the magnitude of $z$ \citep{kassam_fourth-order_2005}.  
 
\subsection{Low-Storage Implicit-Explicit Time Discretization}
In order to to avoid the time step restriction in the stiff regime $\frac{1}{\tau} >> 1$, we have adapted a low-storage implicit explicit (LSIMEX) Runge-Kutta method to the Galerkin-Boltzmann system. Because implicit explicit Runge-Kutta schemes are well documented in the literature (see \citep{kennedy_additive_2003}), we only provide a short description of its efficient application to Galerkin-Boltzmann system. For the ODE system \eqref{eq:time_1}, an LSIMEX scheme is
\begin{align}
\mathbf{q}_{ex} &= 
    \begin{cases}
        \mathbf{q}  & \text{if} \quad i = 1, \\ 
        \mathbf{q} + \left( \tilde{a}_{i, i-1} - \tilde{b}_{i-1} \right) \Delta t  \mathbf{q}_{im} +\left( a_{i, i-1} - b_{i-1} \right) \Delta t \mathbf{q}_{ex} & \text{else}, 
                                    \end{cases} \label{eq:LSIMEX_1} \\
        \mathbf{q}_{im} &=\mathbf{N}\left(\mathbf{q}_{ex} + \tilde{a}_{i,i} \Delta t \mathbf{N}(\mathbf{q}_{im}) \right) \label{eq:LSIMEX_2},\\
        \mathbf{q}_{ex}  &= \mathbf{L}(\mathbf{q}_{ex} + \tilde{a}_{i,i} \Delta t \mathbf{q}_{im}, t_n + c_i \Delta t ) \label{eq:LSIMEX_3},\\
        \mathbf{q} &=\mathbf{q} + \tilde{b}_i \Delta t  \mathbf{q}_{im} + b_i \Delta t \mathbf{q}_{ex}, \label{eq:LSIMEX_4} 
\end{align}
where $i=1\dots s$, $s$ is the stage number. Here $\mathbf{q}_{ex}$  and  $\mathbf{q}_{im}$ denote the explicit and implicit parts of the right hand side of \eqref{eq:time_1} at each stage, and $\tilde{a}$ and $\tilde{b}$ are the coefficients of implicit scheme while $a$,$b$, and $c$ are the explicit scheme coefficients. The Butcher tableau for a class of LSIMEX schemes can be found in \citep{cavaglieri_low-storage_2015}. We use the third order method presented in \citep{cavaglieri_low-storage_2015} in the numerical tests presented below.           

In each stage, the LSIMEX formulation consists of two update steps through \eqref{eq:LSIMEX_1} and \eqref{eq:LSIMEX_4}. Equation \eqref{eq:LSIMEX_3} requires one explicit function evaluation per stage. An important part of this implementation comes from the efficient implicit solve stage given in \eqref{eq:LSIMEX_2}. Since the first three entries of the nonlinear collision term $\mathbf{N}(\mathbf{q})$ are zero, $q_1$, $q_2$ and $q_3$ remain constant when solving for $\mathbf{q}_{im}$ in \eqref{eq:LSIMEX_2}. Furthermore, the remaining three entries of $\mathbf{N}(\mathbf{q})$ are linear in $q_4$, $q_5$ and $q_6$ which allows to us to solve for the entries of $\mathbf{q}_{im}$ without matrix inversion or iterative procedure and reduces the operation to a local node-wise update.

\subsection{Multirate Time Integration}
Explicit time stepping techniques offer numerous computational advantages but their stability is only guaranteed under the Courant-Friedrichs-Lewy (CFL) condition which imposes a limit on maximal time step size. Global stability is then determined by the element having the smallest mesh size. This can result in an increased computational expense, especially in realistic flow problems requiring wide spread of element sizes. 

In most of the realistic flow applications, unstructured meshes are refined around some specific regions to accurately capture the topology of the geometry or complex physical phenomena. Due to the varying resolution and physics of the problem, the CFL stabilty condition is not generally constant in space and time. Ignoring the non-linear relaxation term fo rthe moment, the Galerkin-Boltzmann equations have the following wave-transport time step restriction
\begin{equation}
\label{eq:multirate_1}
    \Delta t \leq \min_{e} C \frac{h^e }{N^2 \lambda_{max}},
\end{equation}
where $h^e$ a the characteristic length of the element, $\EN^e$, $\lambda_{max}=\sqrt{3RT}$ is the maximum wave speed for the Galerkin Boltzmann system in element $\EN^e$, and $C$ is the CFL number which depends on the stability region of the time discretization scheme. Besides the advective time step restriction, the relaxation time, $\tau$ is also an important parameter in designing the time discretization scheme. For the small relaxation times, the Boltzmann equations becomes very stiff. The stiff term, $1/\tau$ is of the order of $Re/$\emph{Ma}$^2$ which is independent of the mesh resolution and polynomial order of approximation. For weakly incompressible and high Reynolds number flows, the time step size required to explicitly step the relaxation term becomes very restricted. On the other hand, semi-analytic time discretization avoids stiffness through analytic integration of some stiff linear terms and recovers the efficiency using only the advective time scale independent of the flow regime. 

The key idea of multirate methods is to achieve a reduced computational expense per time step by partitioning the mesh into groups wherein we advance time using a locally stable time step choice. To construct the semi-analytic multirate groups, we first compute the stable time step of each element using  \eqref{eq:multirate_1}. The global minimum and maximum time step sizes  are denoted by $\Delta t_{\min}$ and $\Delta t_{\max}$, respectively. As noted above, for the sake of simplicity we follow an approach similar to \citep{seny_multirate_2013}, i.e.,  we assume that the successive multirate groups have the time step ratio of $2$. Then, selecting a reference time step, $\Delta t_r$ as the power of two of $\Delta t_{\min}$, the maximum exponent  of the multirate groups is defined as,
\begin{equation}
    \label{eq:multirate_2}
    l^{*} = \log_{2} \frac{\Delta t_{r}}{\Delta t_{\min}},
\end{equation}
which gives $N_l = l^{*} + 1$ groups. We partition the mesh into multirate groups so that elements in the same group have stable time steps in the range $[2^l \Delta t_{\min}, 2^{l+1}\Delta t_{\min}]$. We call $l$ the level of the multirate group.  

We also introduced $N_l -1$ buffer groups to connect the bulk groups. Because elements are weakly connected with their immediate neighbors in DG spatial discretization, we select buffer groups to be the single element layer along the interface of groups of levels $l$ and $l+1$. We store the element numbers in the buffer region for each level and update only the trace values when required for the efficient implementation. 

We adopt the fastest-first approach \citep{gear_multirate_1984} which requires the integration starting from smallest levels i.e., groups with smaller time step sizes. Integration of these levels requires the trace values of one level larger groups at buffer region which is not available and must be interpolated to the sub-step levels. We summarize our implementation of MRSAAB in Table \ref{Table:MRSAAB_Algorithm} for a sample multirate system with $N_l = 3$. In the table, $R$, $U$ and $T$ denote all the required computations i.e., evaluation of the right hand side of \eqref{eq:time_1}, time step update and trace update, respectively. We assume that all history is known for all levels at the beginning of stage $0$ and $t=0$ at this synchronization level. 

\begin{table}[ht]
	\caption{Illustration of MRSAAB algorithm for $N_l =3$ in terms the operations performed in each stage where $R$ is the right hand side evaluation, $U$ is the temporal update, and $T $ is the trace update operation.}
	\centering
	\label{Table:MRSAAB_Algorithm}
	\small\addtolength{\tabcolsep}{5pt}
	\begin{tabular}{|c|c|c|c|c|c|c|c|c|c|}
		\hline
		       & \multicolumn{3}{|c|}{$l=0$} &\multicolumn{3}{|c|}{$l=1$} & \multicolumn{3}{|c|}{$l=2$} \\ \cline{2-10} 
		 Stage & $R$ & $U$ & $T$ & $R$ & $U$ & $T$ & $R$ & $U$ & $T$ \\
		\hline \hline
		 0 & \checkmark & \checkmark & $\times$  & \checkmark & $\times$ & \checkmark & \checkmark & $\times$ &$\times$ \\ \hline
		 1 & \checkmark & \checkmark & $\times$  & $\times$ & \checkmark & $\times$ & $\times$ & $\times$ &\checkmark \\ \hline
		 2 & \checkmark & \checkmark & $\times$  & \checkmark & $\times$ & \checkmark & $\times$  & $\times$ &$\times$ \\ \hline
		 3 & \checkmark & \checkmark & $\times$  & $\times$  & \checkmark & $\times$ & $\times$ & \checkmark &$\times$ \\ \hline
	\end{tabular}
\end{table}
In the first stage, all the levels compute the RHS contributions first. Then, level-$0$ is updated to the time level $\Delta t_{\min}$ using the SAAB coefficients given in \eqref{eq:SAAB_6}. Subsequently, required trace values at the buffer region of level $1$ are extrapolated to the time level $\Delta t_{\min}$ using $\Delta t=2\Delta t_{\min}$ and the half step coefficients listed in \eqref{eq:MRSAAB_8} to evolve level $0$ in the next stage. In the second stage, level $0$  computes the RHS using the extrapolated trace values between level-$0$ and level-$1$ bulk groups. Then, level-$0$ and level-$1$ are updated to $2 \Delta t_{\min}$ with their stable time step sizes. The stage ends with the trace update of level-$2$ for $\Delta t= 4 \Delta t_{\min}$ and half step coefficients that will be used to evolve level-$1$. Next stage starts with RHS evaluation of level-$0$ and level-$1$ and continues with advancing level-$0$ to $3 \Delta t_{\min}$ and extrapolating  the trace values of level-$1$ to the same time level. The final stage brings all levels to the same time with RHS evaluation level-$0$ and update operations for all levels.

\section{GPU Implementation}
\label{Sec:Implementation}
In this section, we give a brief overview of the implementation used in the numerical tests below, conducted using Graphics Processing Unit (GPU) acceleration. The Galerkin-Boltzmann solver described here has been implemented in C/C++ using the Open Concurrent Compute Abstraction (OCCA) API and OKL kernel language \citep{medina_occa:_2014}. OCCA is an abstracted programming model designed to encapsulate native languages for parallel devices such as CUDA, OpenCL, Threads, and OpenMP. OCCA thus offers flexibility in choosing hardware architectures and programming models at run-time by allowing customized implementations of algorithms for several computing devices with a single code. Parallelization on distributed multi-GPU/CPU platforms is achieved using MPI. The source code was compiled using the GNU GCC $5.2.0$  and the Nvidia CUDA V$8.0.61$ NVCC compilers. All the tests presented in the next section were run on a Nvidia Tesla P100 GPU paired with a Xeon E5-$2680$v4 processor. 

The solution process consists of four major computational kernels: (1) evaluation of volume integrals, (2) evaluation of surface integrals, (3) cubature-based integration of the non-linear relaxation contributions, and (4) time-step updates. We refer to each of these processes as the volume, surface, cubature, and update kernels, respectively. In all the kernel implementations, the work load is partitioned in such a way that each thread in a thread block performs computations related to a single node while a thread-block processes multiple elements.

\begin{itemize}

\item {\bf Volume Kernel:} The volume integral terms in the semi-discrete form given in \eqref{eq:spatial_9} are computed in this kernel. The kernel first loads the solution fields from global device memory and loads these fields to shared memory arrays of size $N_p$. Two differentiation matrices are re-used within the kernel taking advantage of L1 or L2 caches depending on the size of matrices. Each thread computes the derivative at a single node by calculating the inner product of a row of each differentiation matrix with the nodal element solution vectors stored in the shared memory. Resulting values and all other necessary data i.e., geometric factors and the Jacobian of local to global transformation, are stored in register memory. This kernel requests $N_p$ threads per element per thread block to perform all required computations. 
\item {\bf Surface Kernel:} The surface kernel computes the contributions of the surface integral term in \eqref{eq:spatial_9}. The structure of the surface kernel is similar to the structure of the volume kernel. The kernel loads the trace data of an element and all its neighbors to registers, computes numerical flux and scales the result with geometric data. Results are then stored in the shared memory array of the size of the total number of face nodes i.e., $N_{f}\times N_{fp}$, where $N_f$ is the number of faces per element and $N_{fp}$ is the number of nodes along each element face. The computed surface fluxes are then lifted to the interpolation nodes by performing matrix-vector multiplication analogously to the volume kernel differentiation action. The surface kernel requires $\max(N_f \times N_{fp}, N_p)$ threads per element per thread block to perform all computations.    
\item {\bf Cubature Kernel:} The nonlinear relaxation term and all $\sigma$ terms in the semi-discrete form \eqref{eq:spatial_9} are evaluated using an appropriately high-order cubature rule in the cubature kernel. Because $\sigma^x$ and $\sigma^y$ are constant in time, we store these variables at cubature integration nodes in global device memory in order to prevent unnecessary interpolation operations. The first operation in this kernel is to copy the elemental field variables, $q^x, q^y$ and $q$ from global memory to shared memory using $N_p$ threads. These shared memory variables are interpolated to the cubature integration points using the interpolation matrix, $\mathcal{I}$ and  $N_c$ threads. Then, nonlinear term and PML terms are computed on cubature node values and stored on the shared memory arrays of size, $N_c$ each. Finally, the cubature kernel performs one more matrix-vector multiplication action for each field to project the results to the interpolation nodes using the operator, $\mathcal{P}$ and $N_p$ threads. To perform all computations, the cubature kernel requests a total of $N_c$ threads per element per thread block.     
\item {\bf Update Kernel:} The update kernel performs time integration updates, which involves global vector operations using the right hand side vectors and necessary amount of solution history depending on the time discretization method and its order. $N_p$ threads per element per thread block are requested by this kernel.
\end{itemize}

Performance of the kernels is highly dependent on the processor architecture, memory usage, and tuning parameters. Here, only the basic performance improvement techniques such as coalescing, loop unrolling and multiple elements per thread block are used. It is possible to use more advance optimization strategies such as utilizing hardware dependent padding, matrix blocking etc., but the study of these is out of the scope of this study. Because, all time stepping techniques studied here use similar kernels with the same level of optimization, we believe that performance results presented in the next section are  independent from the implementation details of computing kernels.   

\section{Numerical Tests}
\label{Sec:NumericalTests}
In this section, we demonstrate the convergence properties, accuracy, and performance of the developed flow solver on distinct PML and non-PML numerical test cases including Couette flow, isothermal vortex problem, and flow around a square cylinder and wall mounted square cylinder.  

\subsection{Unsteady Couette Flow}
As a first computational test we consider shear flow between two parallel plates. The horizontal velocity of the upper plate is specified by $\mathbf{u} = (U,0)$ and a stationary wall boundary condition is enforced at the bottom plate. Periodic flow boundary conditions are enforced at left- and right-hand side of the channel. The Reynolds number, $Re$ of the Couette flow is given by $Re = UL/\nu$ where $L$ is the length of the square channel and $\nu = \tau RT$ is the kinematic viscosity. The analytic solution of the $u$ velocity in the incompressible Navier-Stokes equations is given by  
 \begin{eqnarray*}
u = U\frac{y}{L} + \sum_{n=1}^{\infty}\frac{2U\left(-1\right)^{n}}{\lambda_n L}  e^{ -\nu \lambda_n^2 t}\sin\left( {\lambda_n y}\right),
\end{eqnarray*}
 where $\lambda_n = \frac{n\pi}{L}$.

We derived an analogous Couette flow solution for the Galerkin-Boltzmann equation (the details are listed in the Appendix) which is given in component-wise form as follows,

\begin{eqnarray}
q_2(y,t) &=& \frac{U}{\sqrt{RT}L} y + \frac{1}{\sqrt{RT}}\sum_{n=1}^\infty \frac{2(-1)^n U}{\lambda_n L} \sin(\lambda_n y)e^{\sigma_n t},\nonumber \\ 
q_4(y,t) &=& \frac{1}{RT}\sum_{n=1}^\infty \frac{2(-1)^n U\sigma_n}{\lambda_n^2 L},\nonumber\\ 
q_5(y,t) &=& \frac{U^2}{\sqrt{2}RTL^2}y^2 + \sum_{n=0}^\infty c_n \sin(\lambda_n y)e^{-\frac{t}{\tau}} +  \frac{1}{\sqrt{2}RT}\sum_{n=1}^\infty \frac{2(-1)^n U^2}{\lambda_n L^2(\sigma_n\tau+1)}\sin(\lambda_n y)e^{\sigma_n t} \nonumber \\ &+& \frac{1}{\sqrt{2}RT}\sum_{n=1}^\infty\sum_{m=1}^\infty \frac{4(-1)^n U^2}{\lambda_n\lambda_m L^2(\sigma_n\tau+\sigma_m\tau+1)}\sin(\lambda_n y)\sin(\lambda_m y)e^{(\sigma_n+\sigma_m) t}, \nonumber
\end{eqnarray}
where $q_1=1$, $q_3 = q_6 = 0$, $c_n$ are coefficients chosen to satisfy the initial condition for $q_5$, and $\sigma_n = -\frac{1}{2\tau} + \frac{1}{2\tau}\sqrt{1-4\tau^2RT\lambda_n^2}$, assuming $1\geq 4\tau^2RT\lambda_n^2$.

\begin{figure}[htb!]
\begin{center}
 \begin{subfigure}{0.35\textwidth}
  \includegraphics[width=0.99\textwidth]{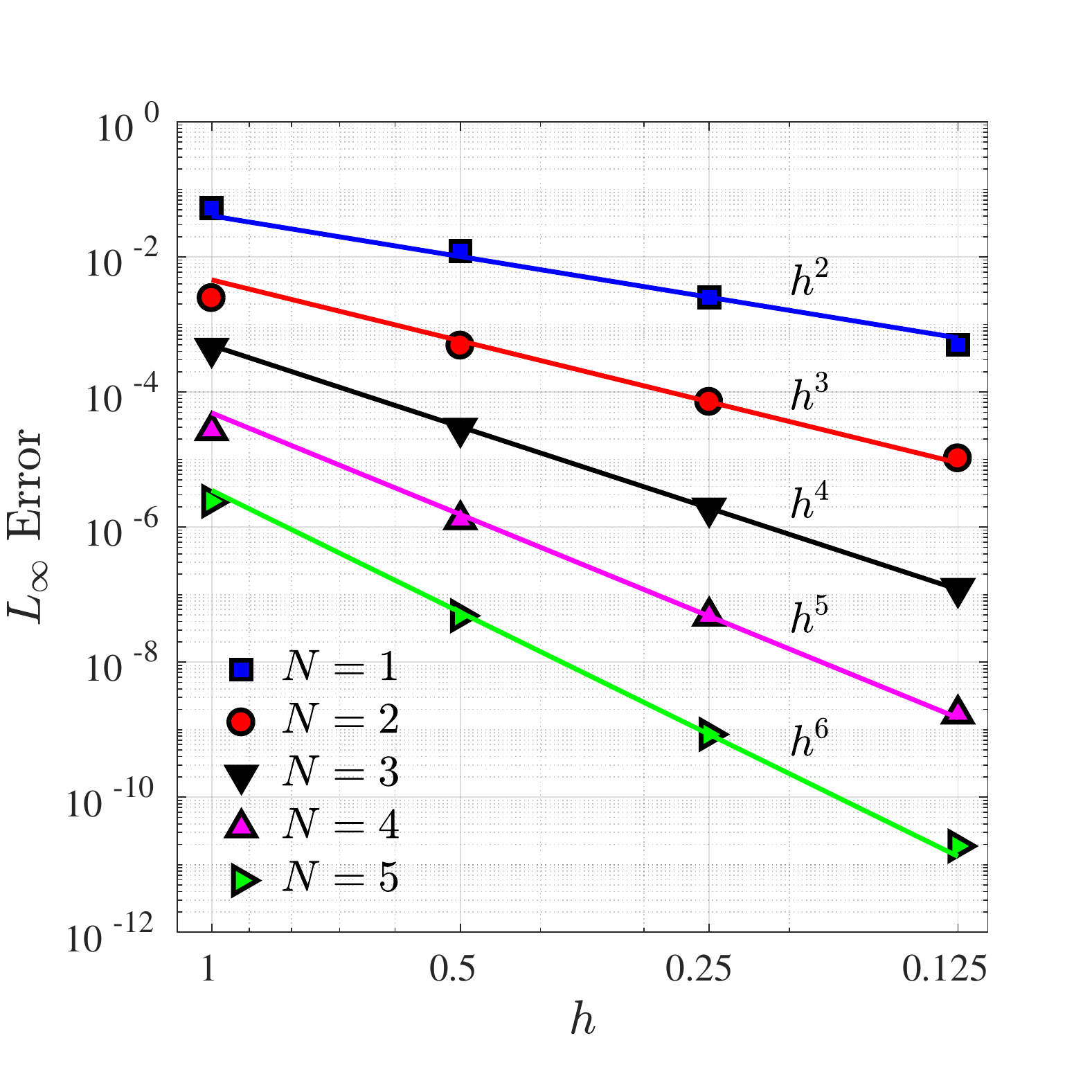}
 \end{subfigure}
 \begin{subfigure}{0.35\textwidth}
  \includegraphics[width=0.99\textwidth]{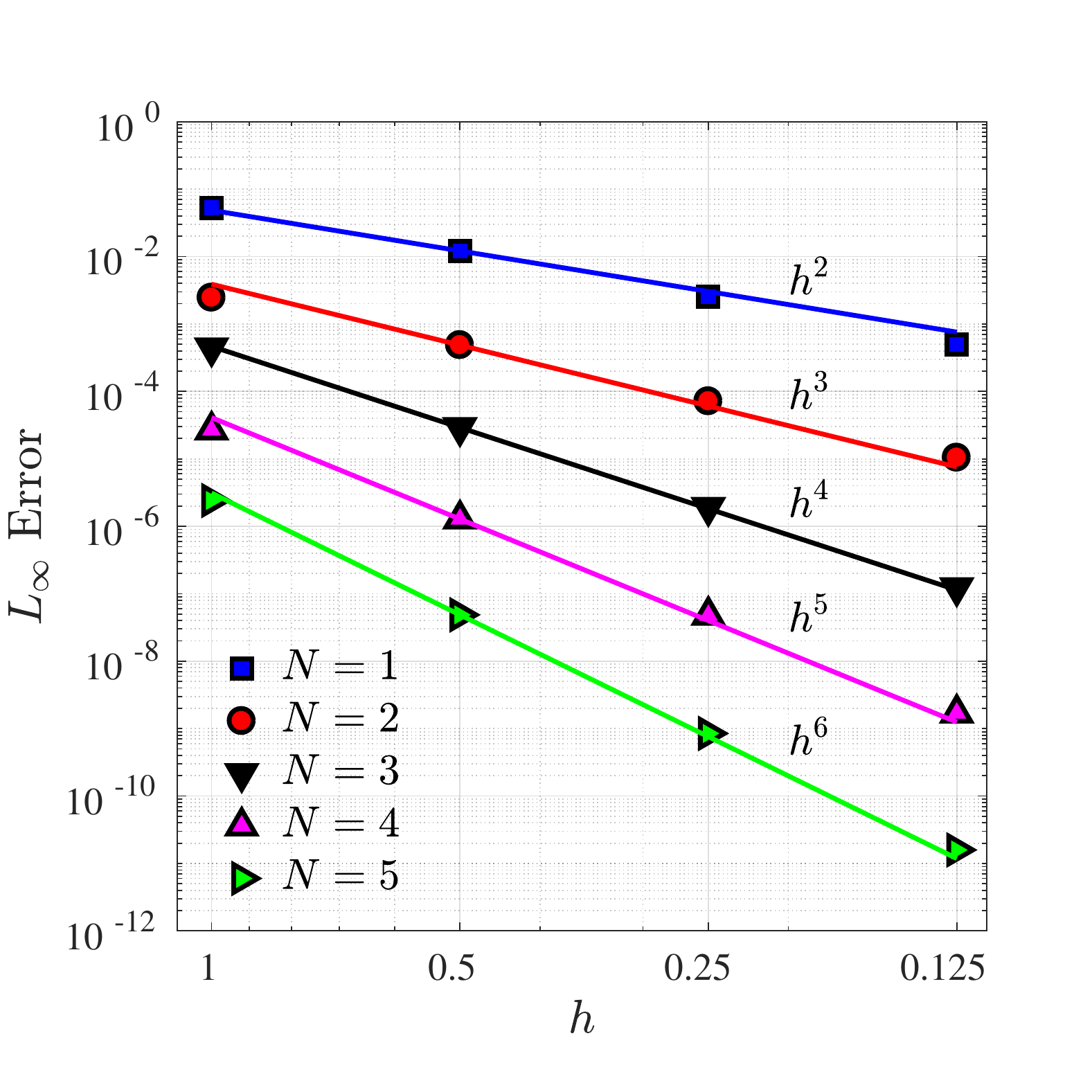}
 \end{subfigure}
  \begin{subfigure}{0.35\textwidth}
  \includegraphics[width=0.99\textwidth]{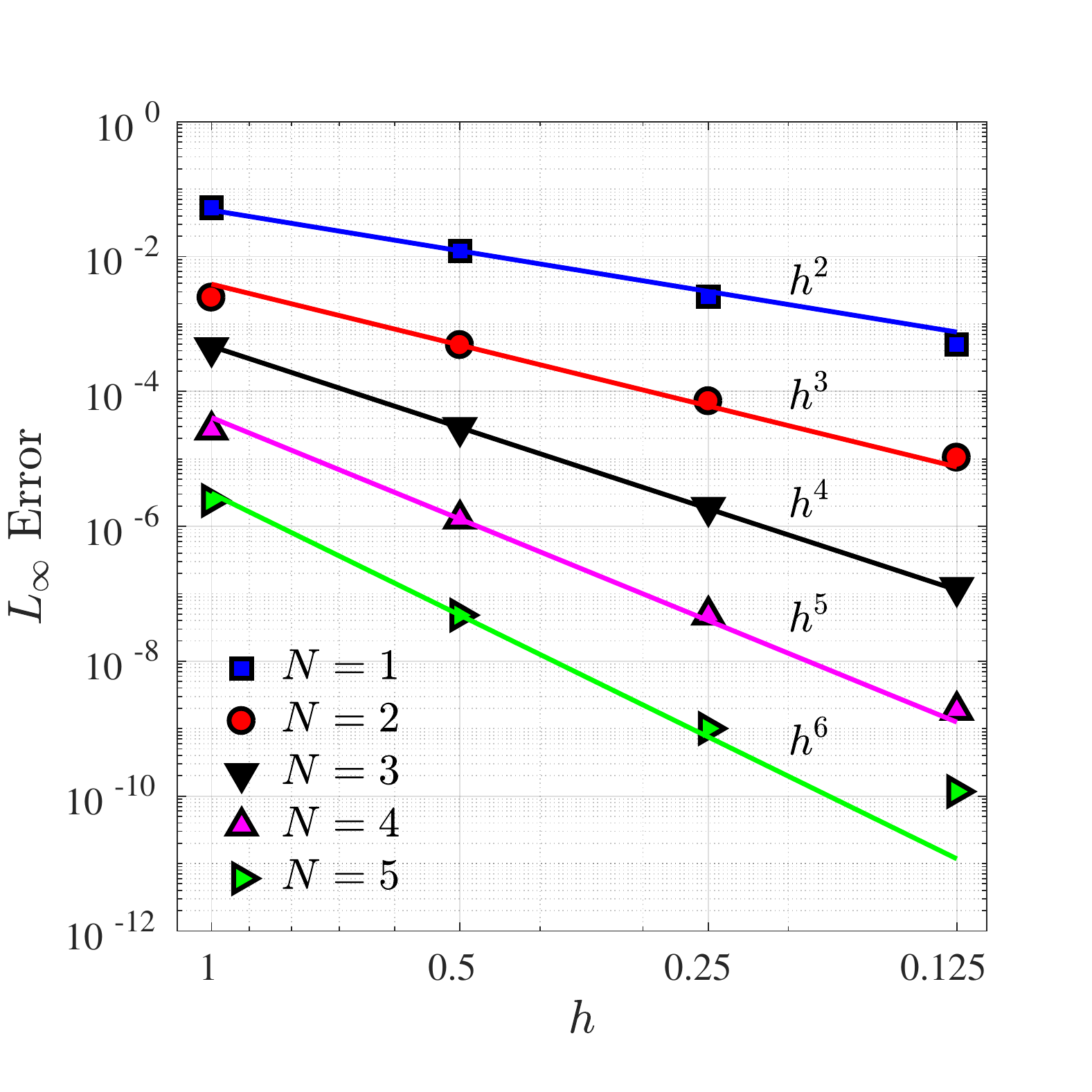}
 \end{subfigure}
 \begin{subfigure}{0.35\textwidth}
  \includegraphics[width=0.99\textwidth]{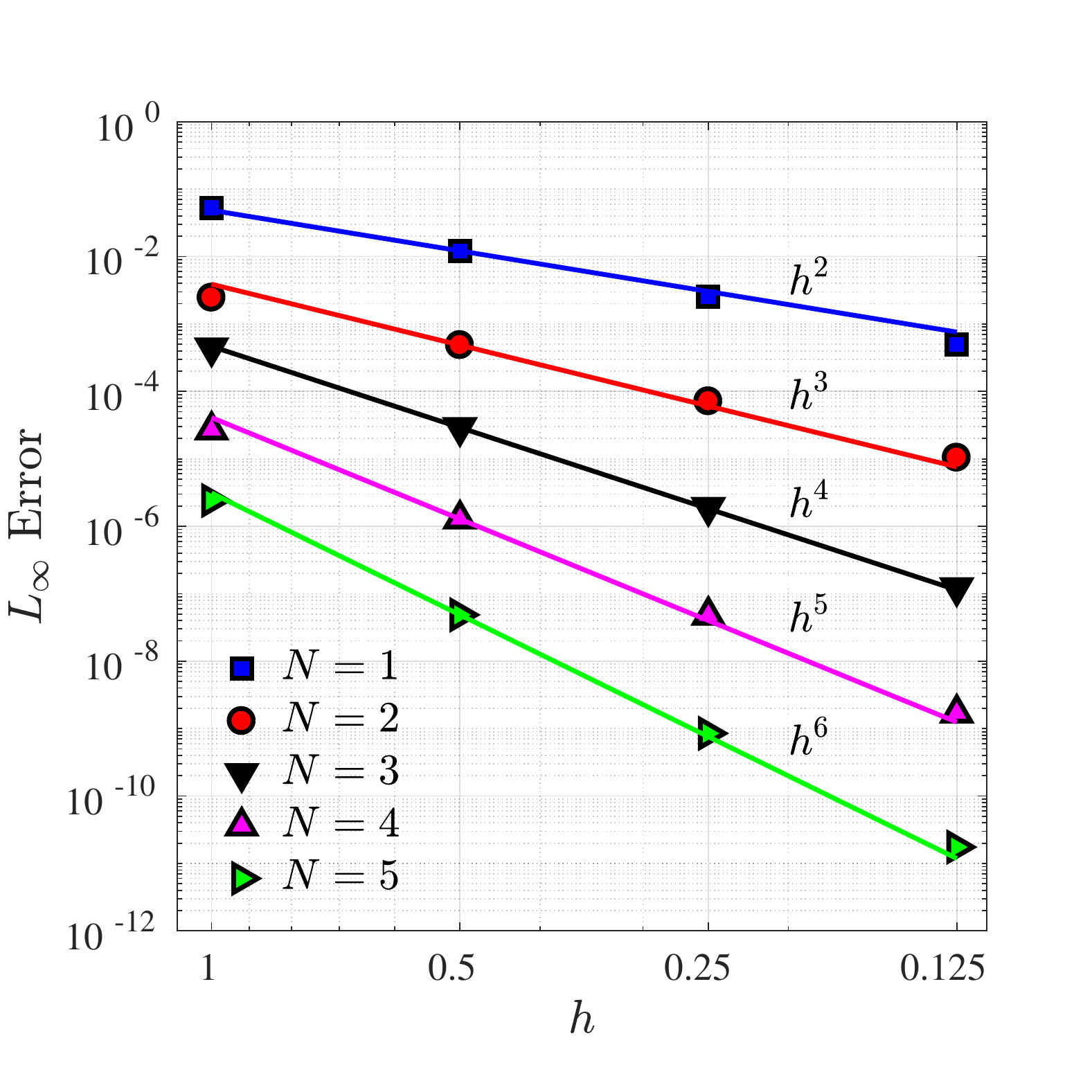}
 \end{subfigure}
\caption{Spatial accuracy test for the unsteady Boltzmann analogy of a Couette flow test problem using $L_\infty$ relative errors for $x$-velocity on successively refined triangular elements. Error plots with reference convergence rate lines are shown for LSERK (top left), SAAB (top right), SARK (bottom left), IMEX (bottom right).}
\label{fig:SpatialConvergence2D}
\end{center}
\end{figure}
We form an exact solution of the Boltzmann equations from the first $10$ modes in the expansions above and solve this problem with the additional initial condition $q_5(y,0)=0$ with \emph{Ma}$=0.1$, $U=1$m/s, $L=1$m, and $\nu = 10^{-2}$ m\textsuperscript{2}/s. Figure \ref{fig:SpatialConvergence2D} shows the computed $L_\infty$ norm of the numerical error for the $x$ component of velocity at the final time $T=1.5$s. We begin with an unstructured mesh of $K=16$ elements and carry out a convergence study with successive uniform mesh refinements and polynomial degree enrichment. The figure demonstrates an $h^{N+1}$ spectral convergence in the numerical error for reference low storage explicit Runge-Kutta (LSERK) and low storage implicit-explicit (LSIMEX) as well as developed semi-analytic Adams-Bashforth (SAAB) and semi-analytic Runge-Kutta (SARK) time integration methods.  
\subsection{Isothermal Vortex Advection}
As a second validation test, we solve an isothermal vortex problem to show the efficacy of our proposed PML formulation. The two-dimensional Euler equations support an advecting vortex solution of the following form \citep{hu_absorbing_2008}, 
\begin{align*}
\rho(\mathbf{x},t)  &=  \rho_r(r)\\
 u(\mathbf{x},t)    &=  U_0 -  u_r(r)\sin{\theta}\\
 v(\mathbf{x},t)    &=  V_0 -  u_r(r)\cos{\theta}
\end{align*}
where $(U_0,V_0)$ is the constant advective velocity, $u_r$ is the given radial velocity, and $r =\sqrt{(x-U_0 t)^2 + (y- V_0t)^2}$.  The radial velocity, density and pressure distribution satisfy the conservation of momentum in the following form,
\begin{align*}
    \frac{d p_r}{dr} = \rho_r \frac{u_r(r)^2}{r}.
\end{align*}
Considering the isothermal flow satisfying  $p_r = \rho_r RT$, density and velocity are related as, 
\begin{align}
    \label{Eq:Vortex_density}
    \frac{d p_r}{\rho_r} = \frac{1}{RT}\frac{u_r(r)^2}{r}.
\end{align}
We consider a radial velocity distribution in the form of
\begin{align}
\label{Eq:Vortex_velocity}
    u_r(r) = \frac{U_{\max}}{b} r e^{  \frac{1}{2} \left( 1-\frac{r^2}{b^2}\right) },
\end{align}
where $U_{\max}$ is the maximum velocity at $r=b$. Density is obtained by integrating \eqref{Eq:Vortex_density} from infinity to $r$ as, 
\begin{align}
\label{Eq:Vortex_density2}
    \ln{\frac{\rho_r}{\rho_\infty}} = -\frac{U_{\max}^2}{2RT}e^{  \frac{1}{2} \left( 1-\frac{r^2}{b^2}\right) }.
\end{align}
Equations \eqref{Eq:Vortex_density2} and \eqref{Eq:Vortex_velocity} are used to initialize the solution with parameters $U_0 =0.5$, $V_0 =0$ and $b=0.2$. The computational domain is set to be the bi-unit square including a surrounding PML domain of width $w$, i.e. $[-1-w,1+w]^2$. The domain discretized with uniform unstructured triangular elements with characteristic length, $h=0.1$. 

Unless stated explicitly otherwise, the PML absorption coefficient is taken to be
\begin{align}
\label{Eq:PML_coefficient}
\hat{\sigma}^{x} = \sigma_{\max} \left|\frac{x-x_0}{w}\right|^\alpha,
\end{align}
where $x_0 = \pm 1$ is the location of interface between physical domain and the PML region. $\hat{\sigma}^{y}$ is computed using an analogous expression.

\begin{figure}[htb!]
\begin{center}
 \begin{subfigure}{0.35\textwidth}
  \includegraphics[width=0.99\textwidth]{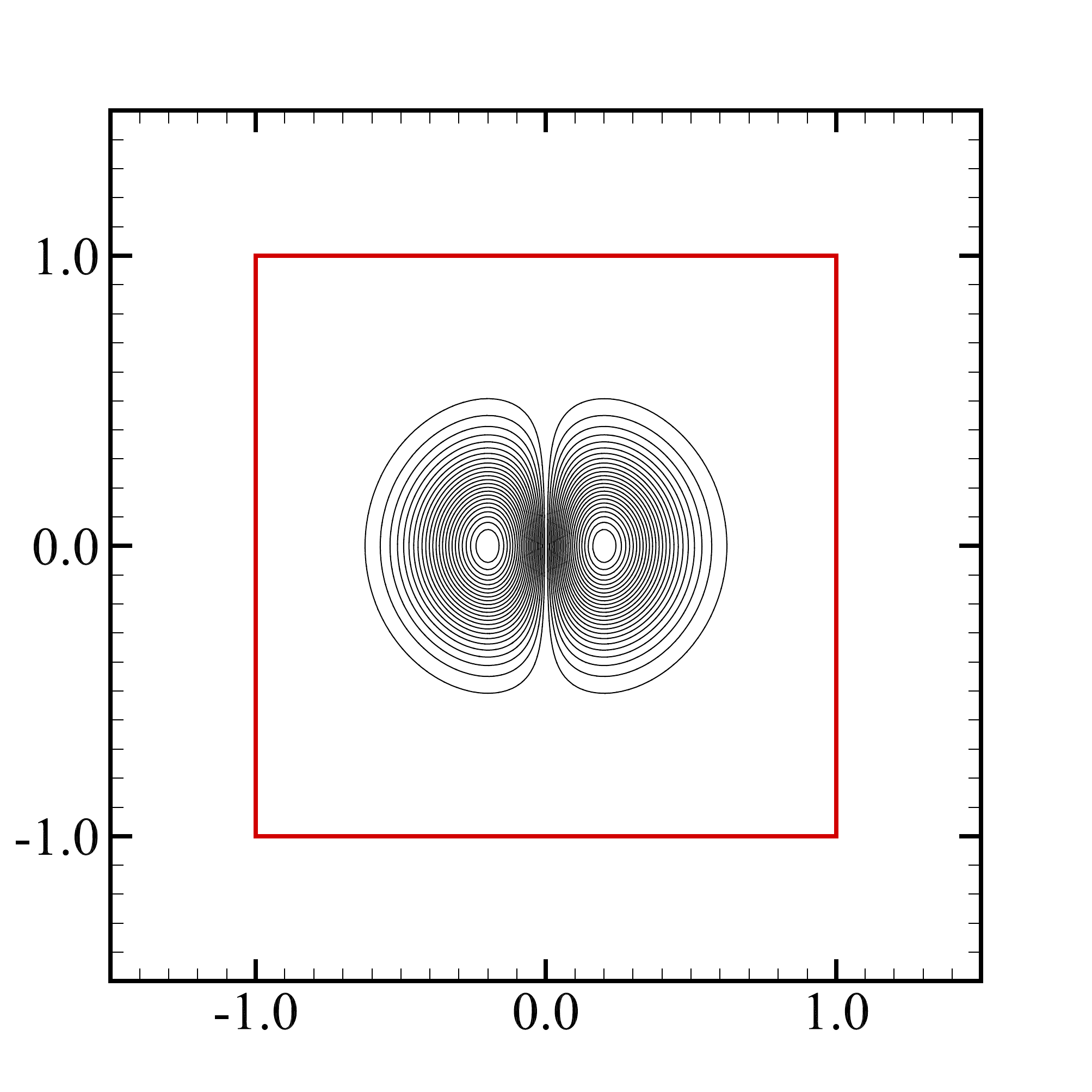}
  \caption{$t=0$}
 \end{subfigure}
 \begin{subfigure}{0.35\textwidth}
  \includegraphics[width=0.99\textwidth]{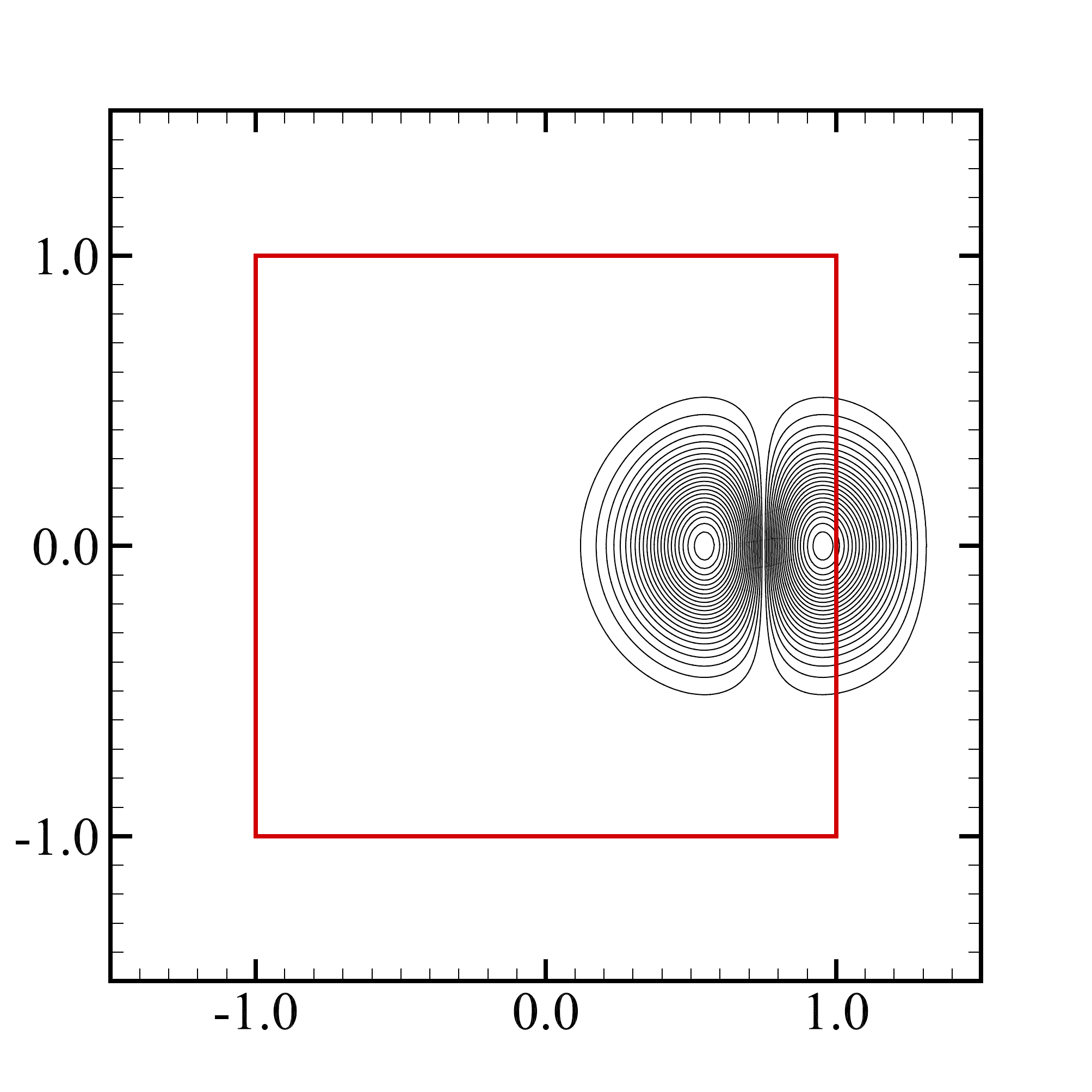}
  \caption{$t=1.5$}
  \end{subfigure}
  \begin{subfigure}{0.35\textwidth}
  \includegraphics[width=0.99\textwidth]{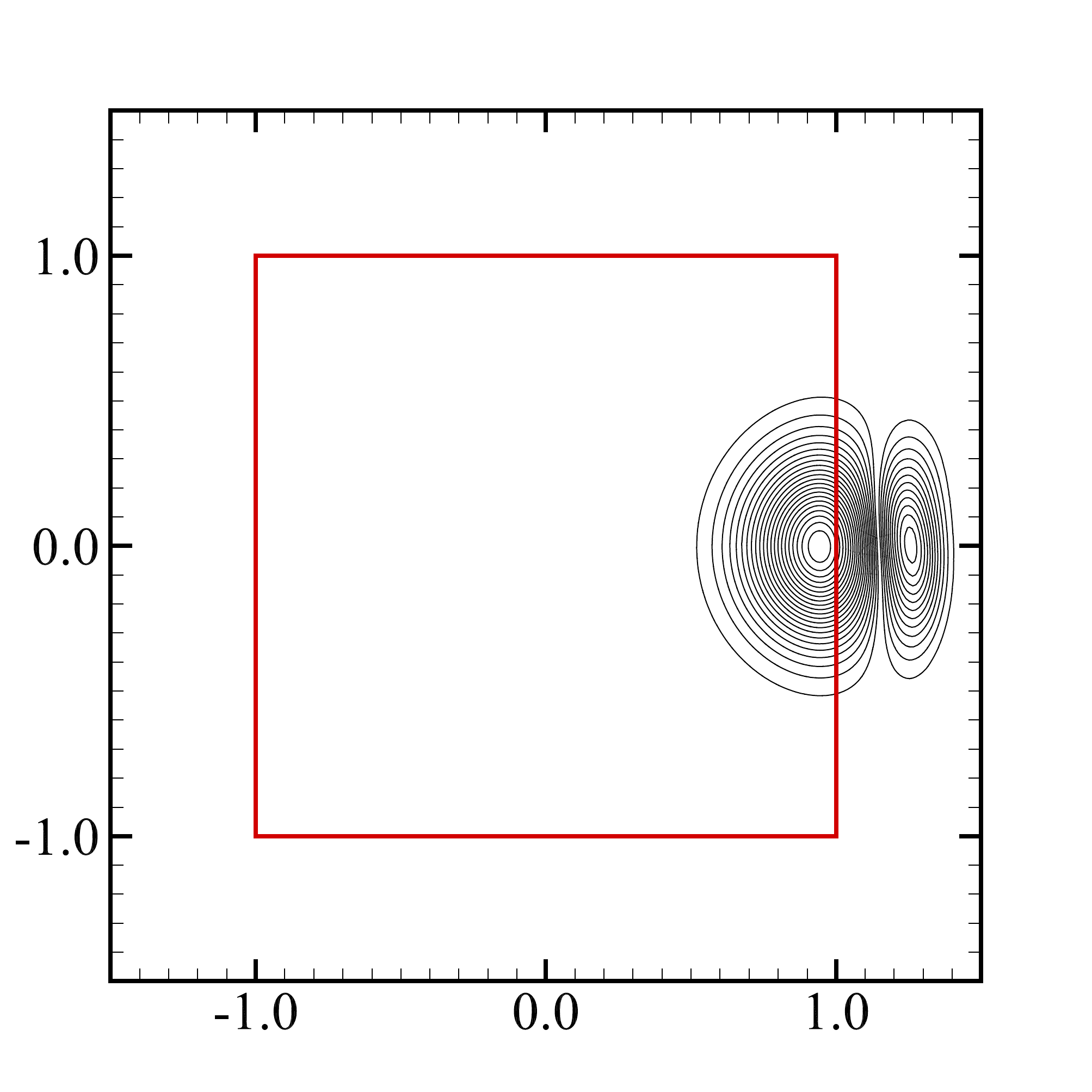}
  \caption{$t=2.3$}
  \end{subfigure}
  \begin{subfigure}{0.35\textwidth}
  \includegraphics[width=0.99\textwidth]{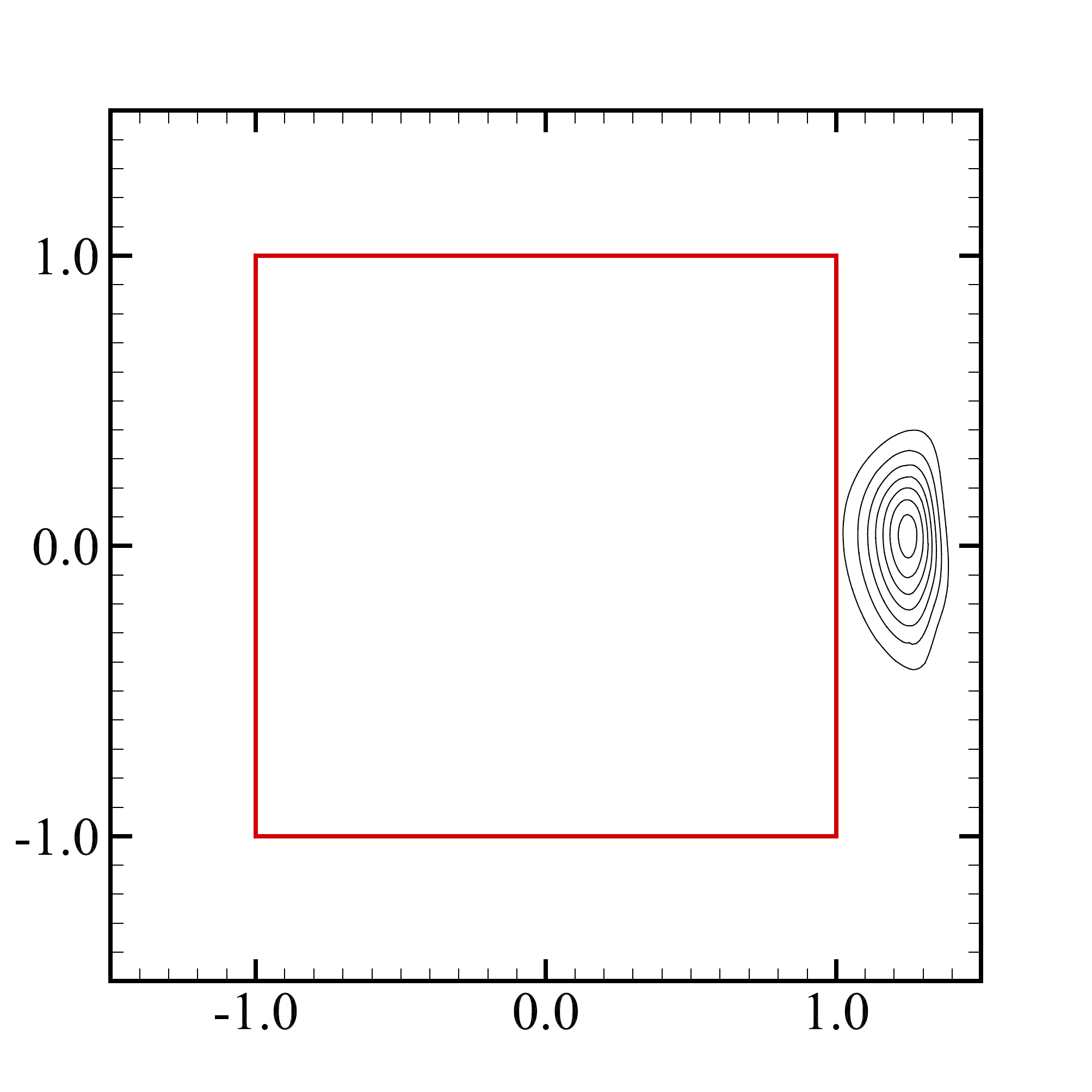}
  \caption{$t=3.3$}
 \end{subfigure}
\caption{Isothermal vortex propagation test for $Re=1000$ and $N=5$ on the domain $[-1.5\times1.5]^2$ with $w=0.5$ PML width. Contours show the $y$-velocity from $-0.25$ to $0.25$ with the increment of $0.0125$ excluding the zero level.}
\label{fig:VortexContour}
\end{center}
\end{figure}
Figure \ref{fig:VortexContour} shows the $v-$velocity contours of a numerical solution at time $t=0, 1.5, 2.3$ and $3.3$, respectively, for $U_{\max} = 0.5 U_0$ and PML width $w=0.5$. The solution is obtained for $N=5$ and Reynolds number $Re=1000$ to preserve the vortex strength at PML and physical domain interface. In the PML region, we select maximum damping coefficient, $\sigma_{\max} = 20$, a fourth-order profile, $\alpha =4$, and multidimensional coefficients $\alpha^x = \alpha^y = 0.1$  for this particular numerical solution. The vortex preserves symmetry while entering the absorbing layer, indicating minimal reflections at the interface. Also, absorption of the vortex in the PML region can be clearly observed in Figure \ref{fig:VortexContour}.

\begin{figure}[htb!]
\begin{center}
  \begin{subfigure}{0.35\textwidth}
  \includegraphics[width=0.99\textwidth]{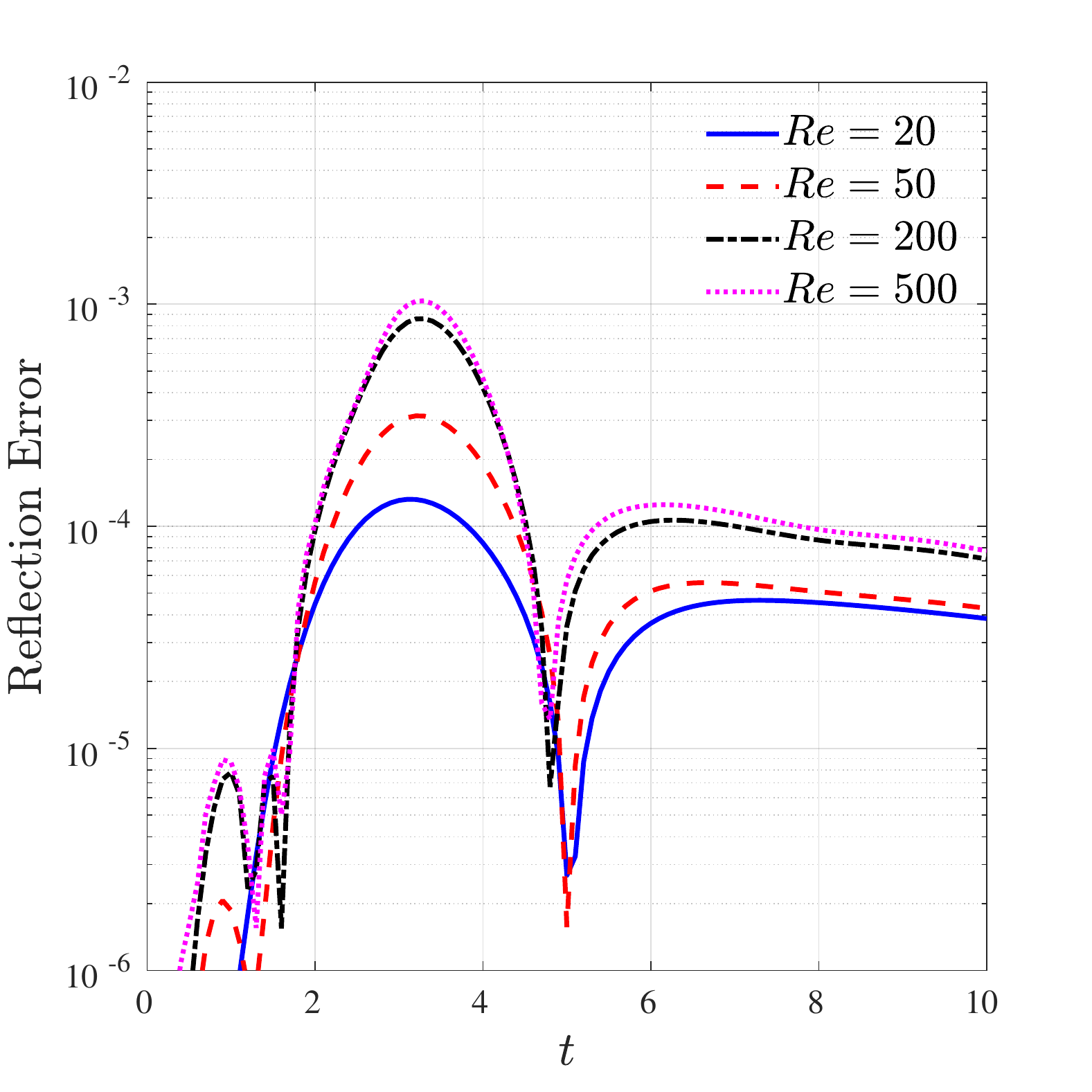}
  \caption{}
 \end{subfigure}
  \begin{subfigure}{0.35\textwidth}
  \includegraphics[width=0.99\textwidth]{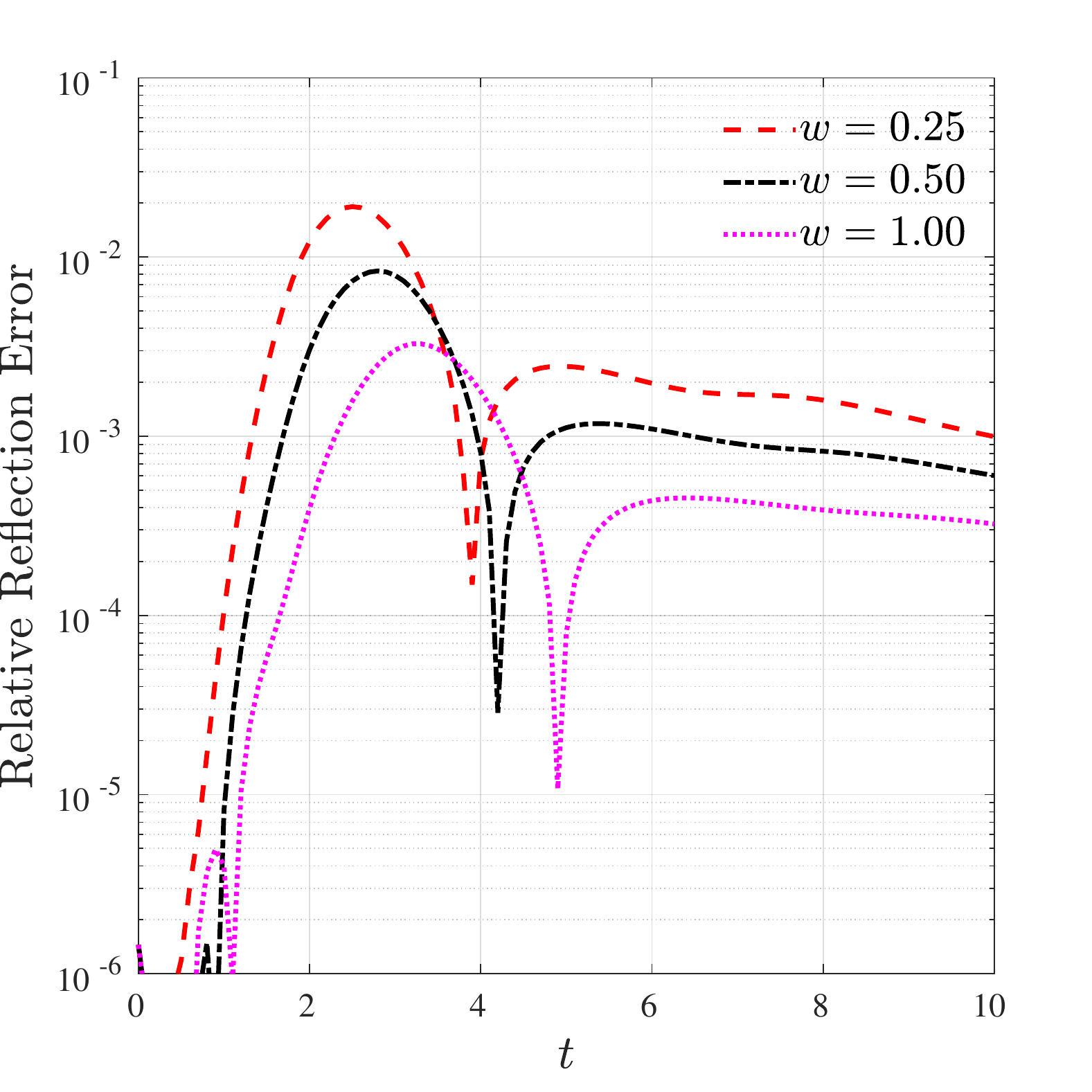}
  \caption{}
 \end{subfigure}
 \begin{subfigure}{0.35\textwidth}
  \includegraphics[width=0.99\textwidth]{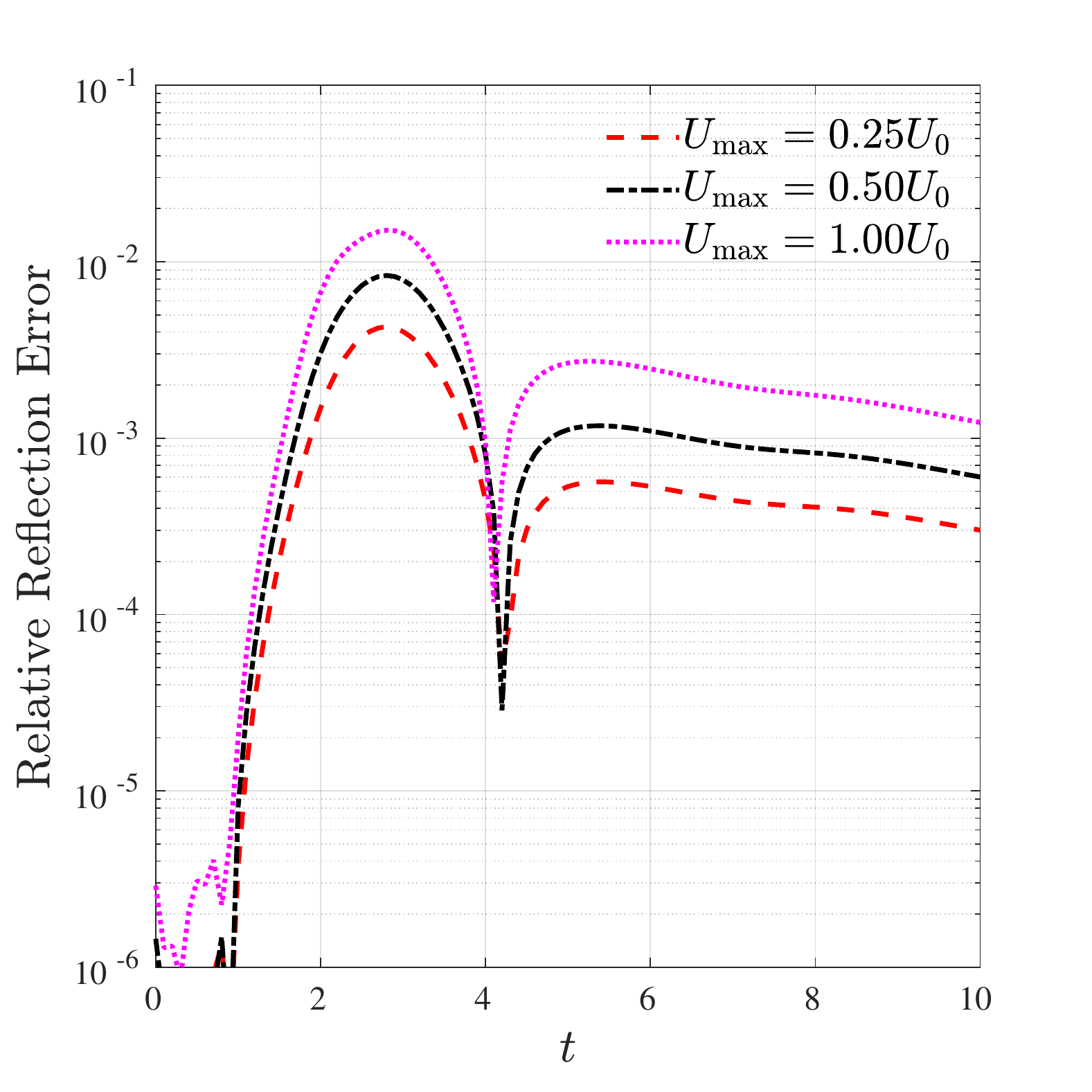}
  \caption{}
 \end{subfigure}
\caption{Isothermal vortex propagation test for $N=3$. Reflection error relative to $U_{\max}$ on the v-velocity component is computed at $(0.9, 0.0)$ (a) for various $Re$ numbers (b) for various  PML widths with  $U_{\max} = 0.5U_0$,$Re=1000$ and (c) for various vortex strengths with  $w= 0.5$, $Re=1000$. }
\label{fig:VortexErrorParametric}
\end{center}
\end{figure}
Figure \ref{fig:VortexErrorParametric} (a) shows the maximum difference between the numerical solution and a reference solution obtained using computational domain which is large enough so that reflections of initial pressure waves do not pollute the solution in the domain of interest. Reflection error is computed at the point, $(0.9,0.0)$ for $N=3$ and for various Reynolds numbers as a function of time. The relative error begins quite small and peaks around $10^{-3}$ for high Reynolds numbers. The relative error decreases with Reynolds number due to weaker vortex strengths with increasing viscous dissipation. To further investigate the maximum reflection error, Figure \ref{fig:VortexErrorParametric} (b-c) show the difference between PML and reference solutions in $v$-velocity component for $Re=1000$ in relation to the PML widths and vortex strengths as a function of time. For the fixed vortex strength of $U_{\max} = 0.5U_0$,  relative reflection error decreases with the increase of the PML width and decreases in time, as expected. In the test case illustrating the effect of the vortex strength on the relative reflection error, a background uniform flow is taken as $U_0 = 0.5$ and the maximum velocity of the vortex is increased from $U_{\max} = 0.25U_0$ to  $U_{\max} = 1.0 U_0$. Although the error increases with the strength of vortex, hence the non-linearity in the system equation, the relative error of less than $1\%$ is achieved for PML width $w=0.5$.

\subsection{Flow Around Square Cylinder}
In our next test we study the accuracy of the Galerkin-Boltzmann approximation and the performance of PML formulation through solving vortex shedding behind a square cylinder test problem. The uniform incoming flow has a Mach number \emph{Ma}$= U_{\infty}/a_{\infty} = 0.2$ where $U_{\infty}$ and $a_{\infty}$ denote the velocity of uniform flow and the speed of sound, respectively. For the present computations, Reynolds number $Re= U_{\infty} d/ v_{\infty} $ is taken as $150$, where $d$ is the characteristic length of the domain and $v_{\infty}$ is the reference kinematic viscosity.  

We solve the problem on a rectangular internal domain $[-5.4, 9.4] \times [-5.4, 5.4]$  surrounded by a PML region of constant width $w$ in all directions. The square cylinder is located at $(0.0,0.0)$ with a unit edge length. The computational domain is discretized with $K=25,000$ unstructured triangular elements with increased resolution near the square cylinder to resolve the boundary layer. 

\begin{figure}[t]
\begin{center}
 \begin{subfigure}{0.45\textwidth}
  \includegraphics[width=0.99\textwidth]{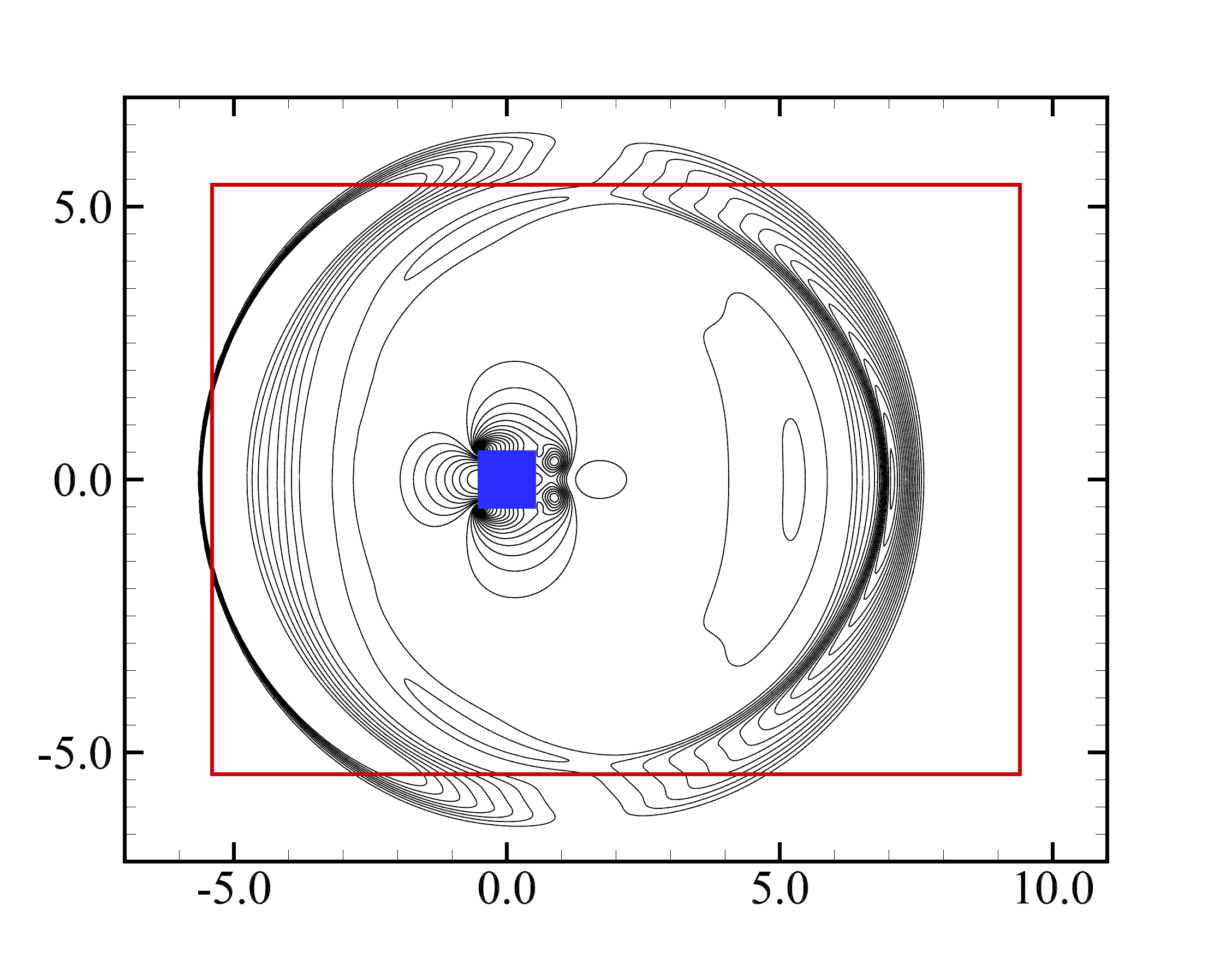}
  \caption{$t=1.1$}
 \end{subfigure}
~
 \begin{subfigure}{0.45\textwidth}
  \includegraphics[width=0.99\textwidth]{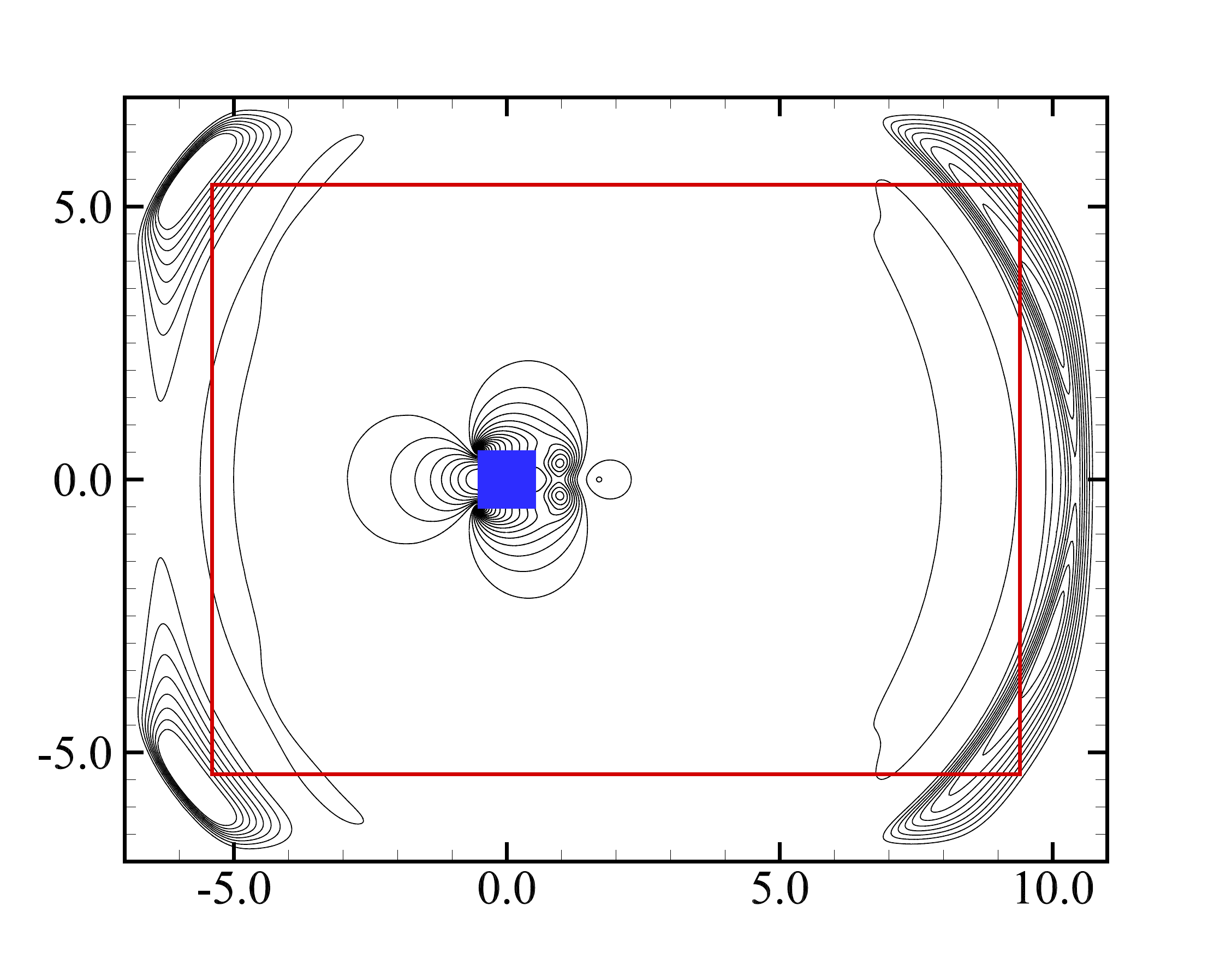}
  \caption{$t=1.8$}
  \end{subfigure}
 ~ 
  \begin{subfigure}{0.45\textwidth}
  \includegraphics[width=0.99\textwidth]{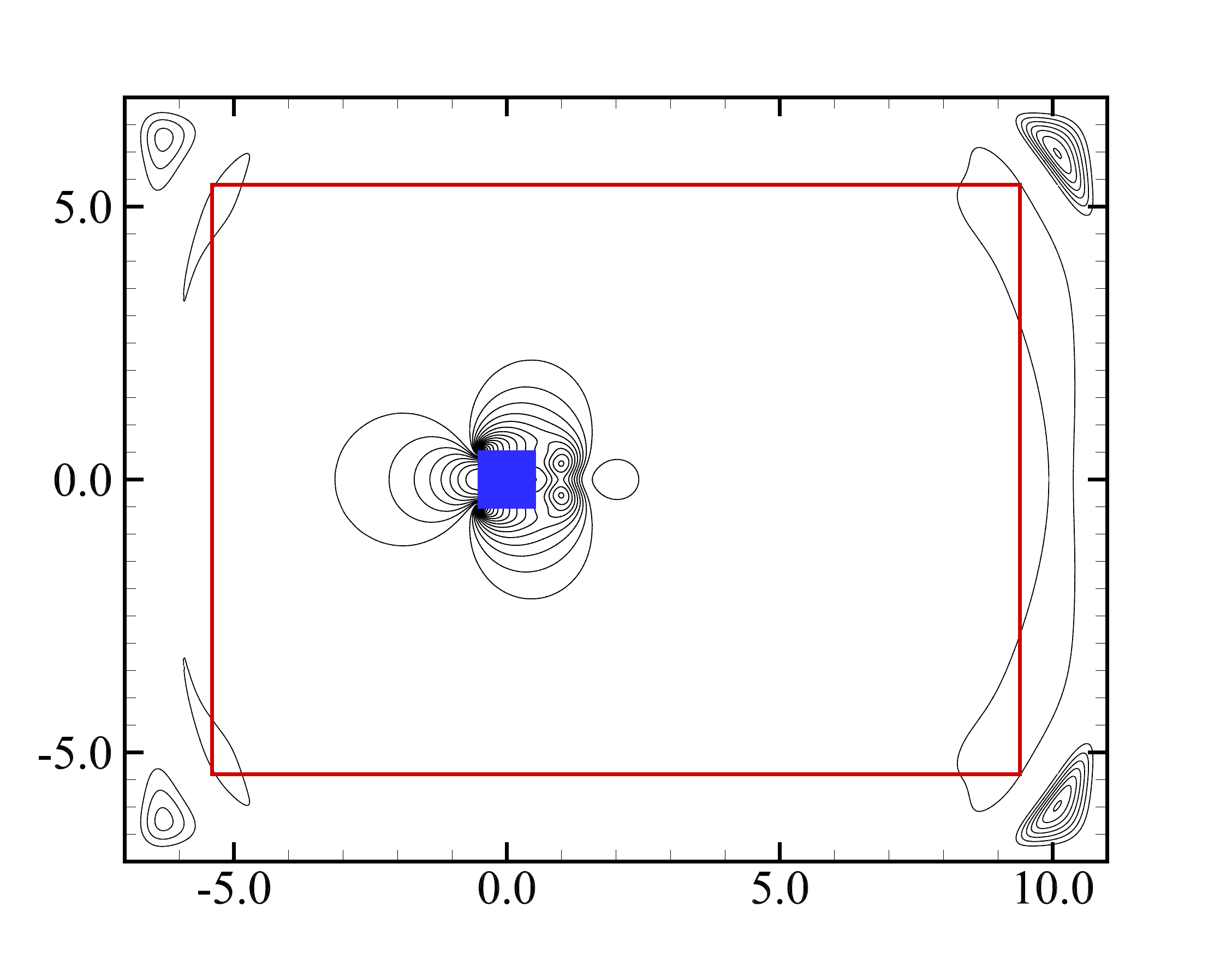}
  \caption{$t=2.1$}
 \end{subfigure}
 ~
 \begin{subfigure}{0.45\textwidth}
  \includegraphics[width=0.99\textwidth]{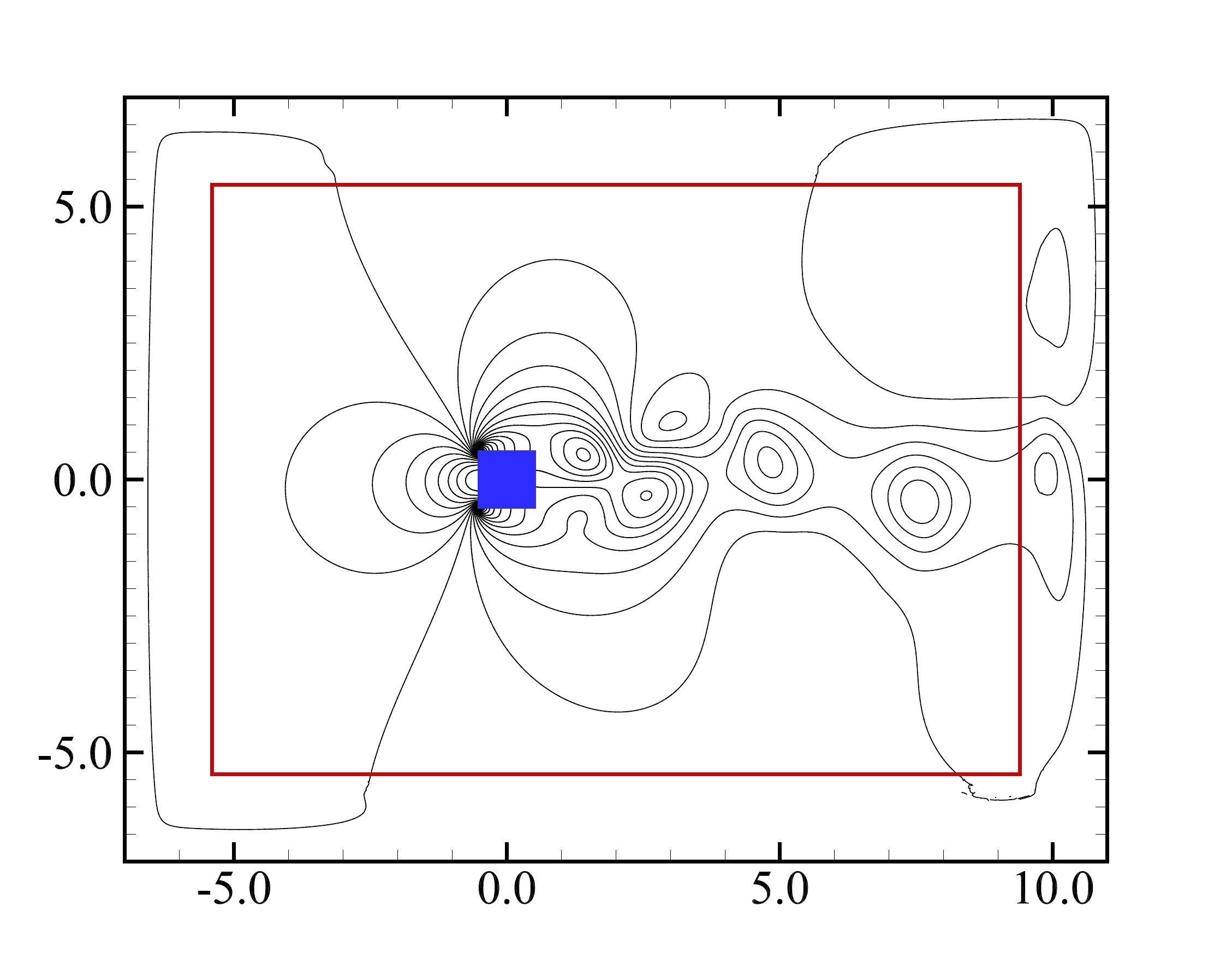}
  \caption{$t=100$}
 \end{subfigure}
\caption{Flow around square cylinder test problem for $Re=150$, $N=3$ on the domain of $[-7,11] \times [-7, 7]$ with a PML region of width, $w=1.6$. Contours show the pressure field from $22.5$ to $25.5$ with the increment of $0.125$. }
\label{fig:SQCylinderPressure}
\end{center}
\end{figure}
Figure \ref{fig:SQCylinderPressure} shows the instantaneous pressure contours for $N=3$ and $w=1.6$ at different solution times.  For this test, the PML parameters are selected similar to the vortex propagation test, i.e. $\sigma_{\max} = 20$, $\alpha = 4$ and $\alpha^x = \alpha^y = 0.1$. Figure \ref{fig:SQCylinderPressure} (a)-(c) clearly demonstrate that initial transient pressure waves are damped out efficiently in the PML region without any noticeable reflections between PML interface  and domain boundaries. Figure \ref{fig:SQCylinderPressure} (d) shows a snapshot of  pressure field after vortex shedding starts and the shear waves dominate the flow. The PML region also performs well in this regime where no visible reflections are observed in the pressure field for this long time simulation.

\begin{figure}[htb!]
\begin{center}
 \begin{subfigure}{0.48\textwidth}
  \includegraphics[width=0.99\textwidth]{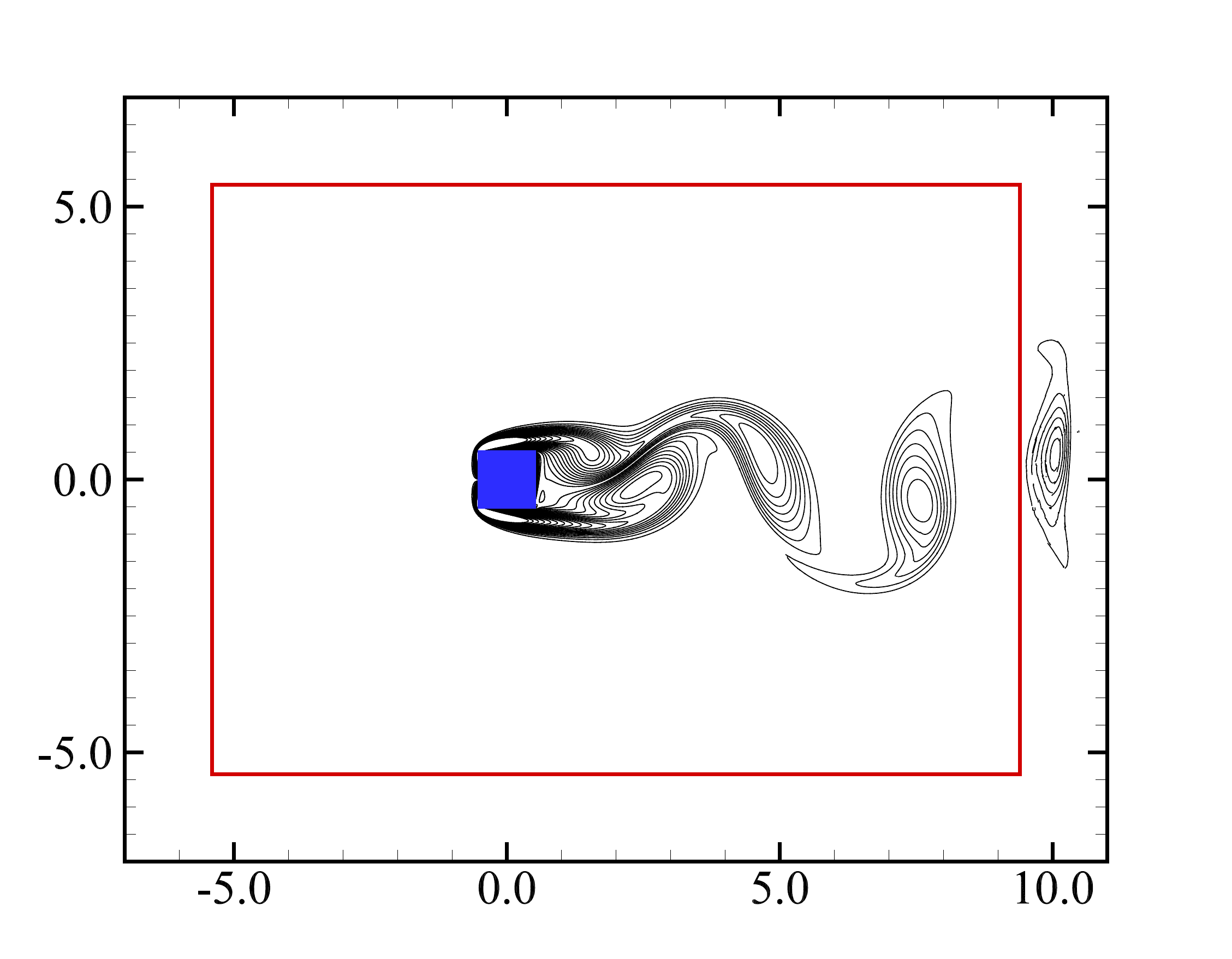}
 \end{subfigure}
~
 \begin{subfigure}{0.48\textwidth}
  \includegraphics[width=0.99\textwidth]{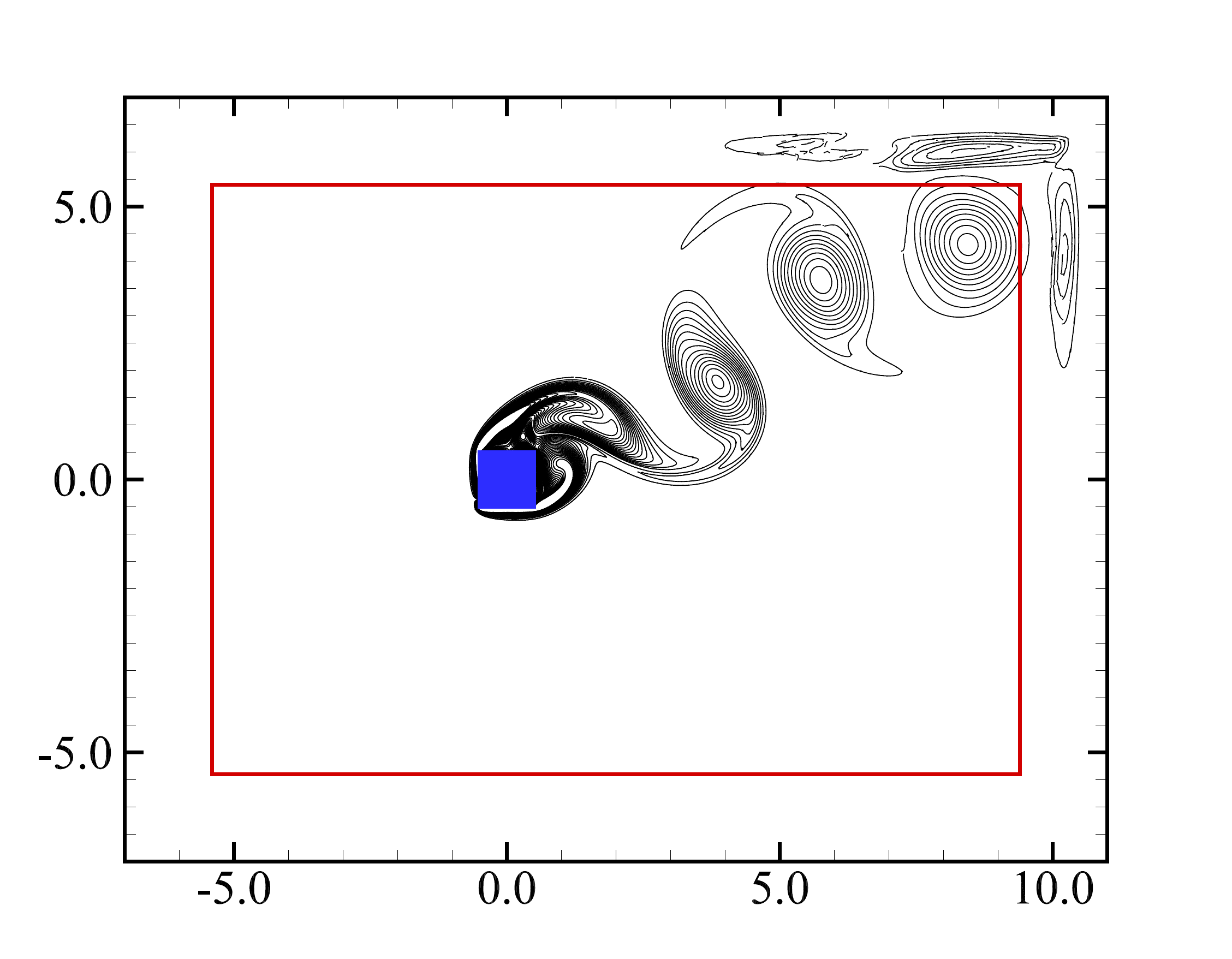}
  \end{subfigure}
\caption{Flow around square cylinder test problem for $Re=150$, $N=5$ on the domain of $[-7,11] \times [-7, 7]$ with a PML region of width, $w=1.6$. Contours show the vorticity field from $-5.0$ to $5.0$ with the increment of $0.25$ excluding the zero level.}
\label{fig:SQCylinderVorticity}
\end{center}
\end{figure}
\begin{figure}[htb!]
\begin{center}
 \begin{subfigure}{0.45\textwidth}
  \includegraphics[width=0.99\textwidth]{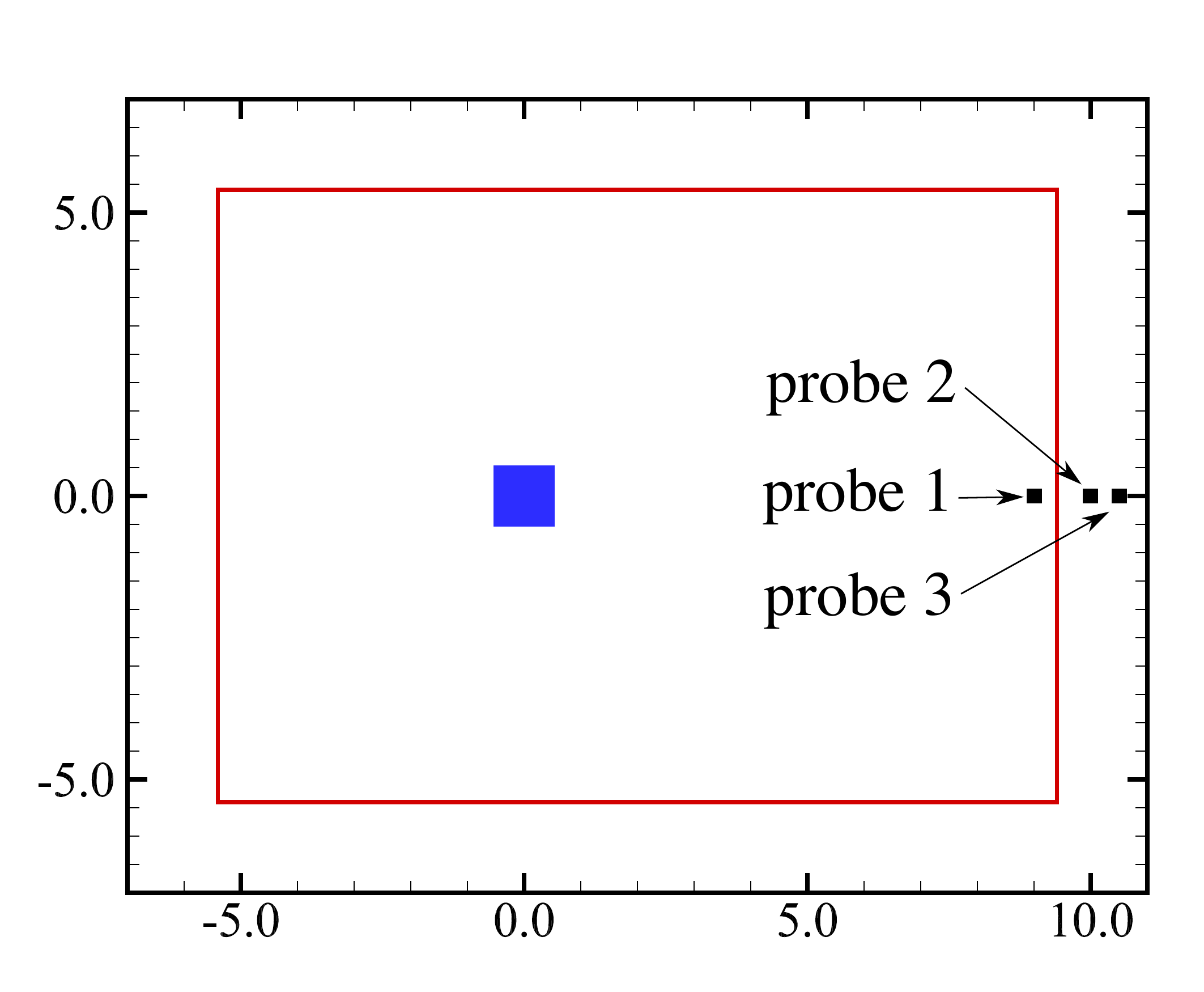}
 \end{subfigure}
 \begin{subfigure}{0.45\textwidth}
  \includegraphics[width=0.99\textwidth]{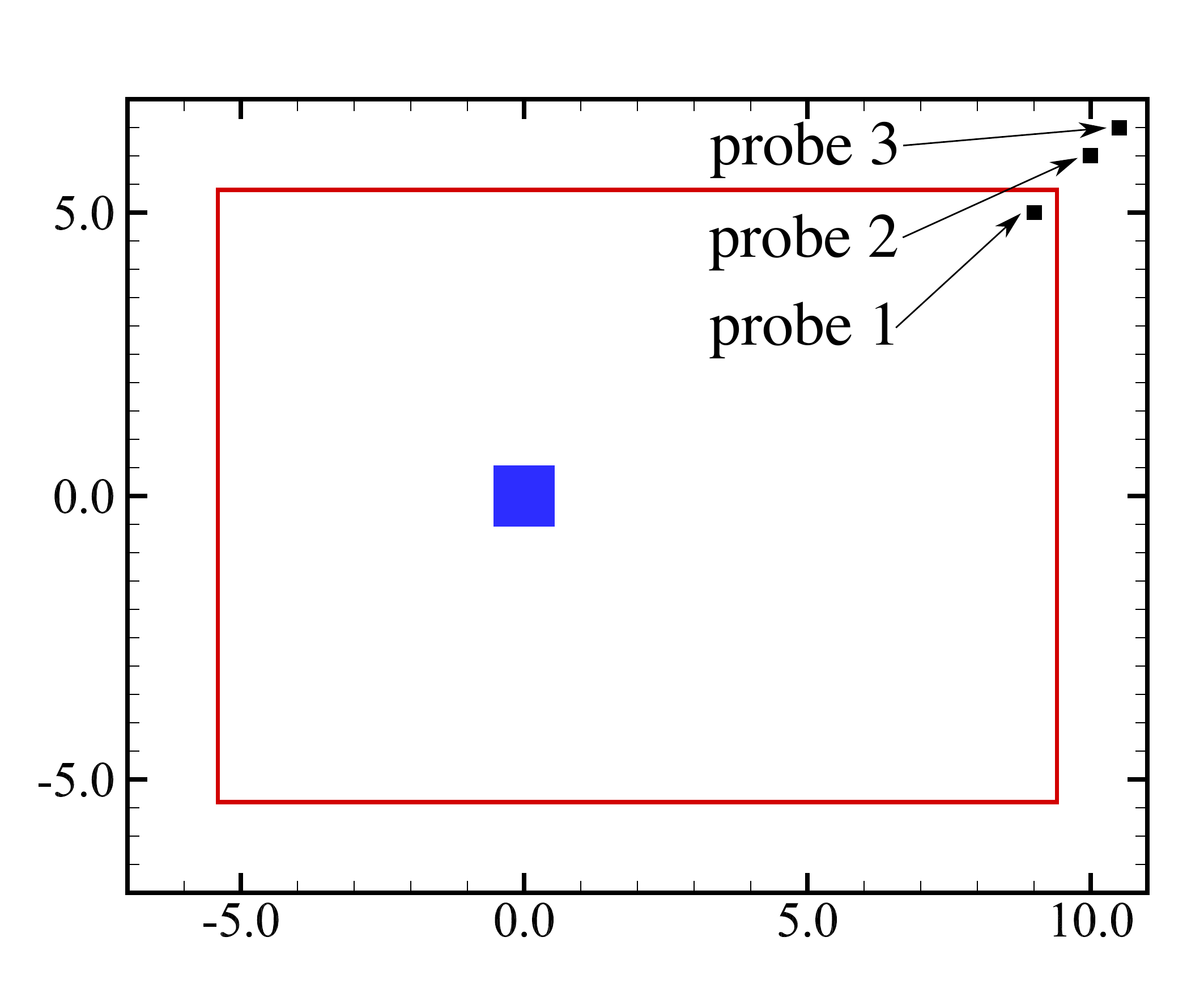}
 \end{subfigure}
 \begin{subfigure}{0.43\textwidth}
  \includegraphics[width=0.99\textwidth]{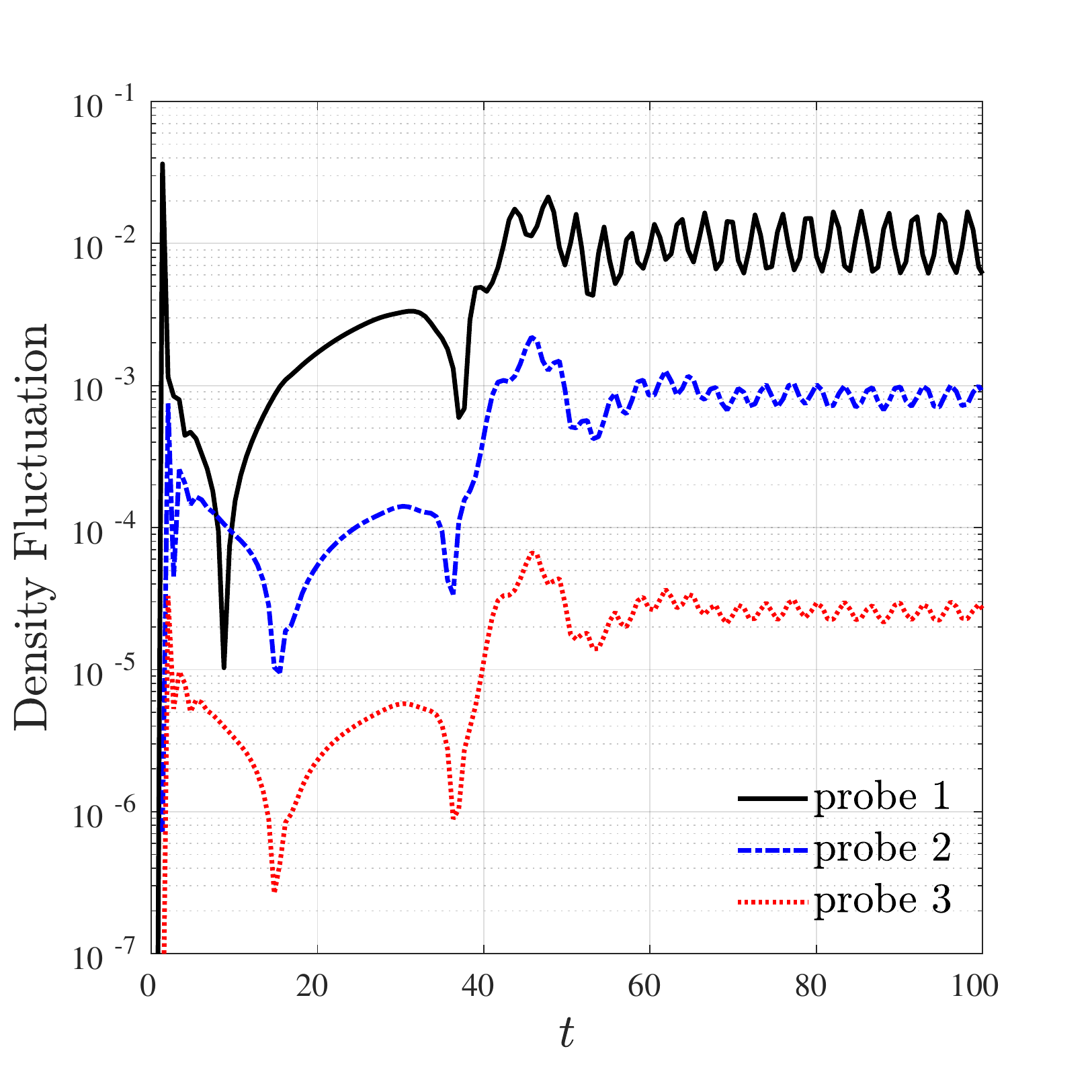}
  \caption{}
 \end{subfigure}
  \begin{subfigure}{0.43\textwidth}
  \includegraphics[width=0.99\textwidth]{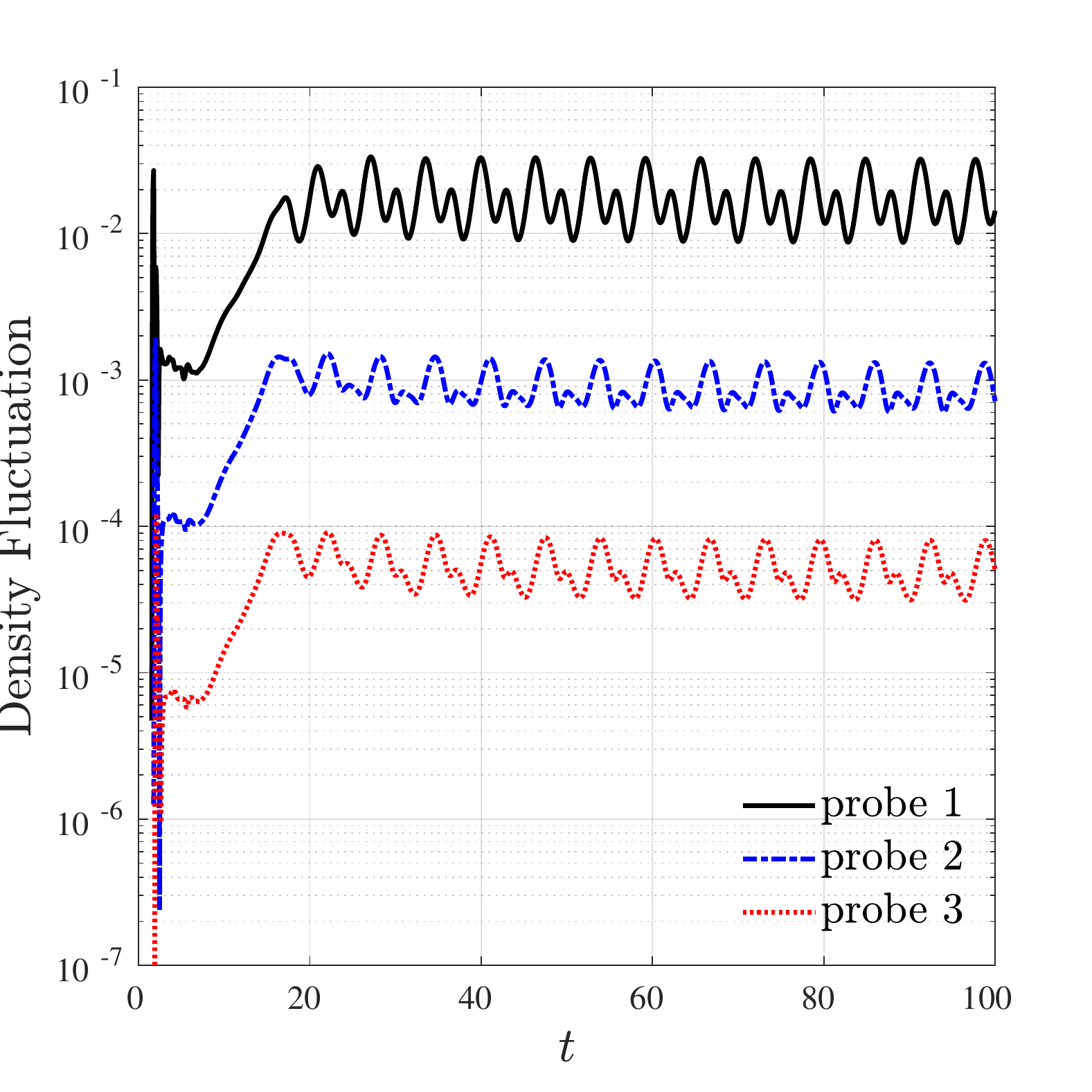}
  \caption{}
 \end{subfigure}
\caption{Flow around square cylinder test problem for $Re=150$, $N=5$ on the domain $[-7,11] \times [-7, 7]$ with a PML region of width $w=1.6$. Density fluctuation history for different probe locations for (a) zero angle of attack and (b) 30 (deg) angle of attack. }
\label{fig:SQCylinderDensity}
\end{center}
\end{figure}

Figure \ref{fig:SQCylinderVorticity} gives two snapshots of the vorticity contours for $Re=150$, $N=5$ and the same PML settings as in the previous test. Absorption of the nonlinear vortices by the PML is clearly seen in the Figure \ref{fig:SQCylinderVorticity} (a) for zero angle of attack.  To demonstrate stability and effectiveness of the present PML formulation for different mean flow directions, the same problem is solved for an angle attack of 30 (deg). As seen in Figure \ref{fig:SQCylinderVorticity} (b), the PML absorbs the incoming vortices almost completely. For this test, the damping efficiency of the PML is shown in \ref{fig:SQCylinderDensity} in terms of density fluctuations, $|\rho/\rho_\infty -1|$ for both zero and $30$ (deg) angle of attack problems on the probes located in three different locations. For zero angle of attack case, probe 1, probe 2, and probe 3 are located at $(9.0,0.0)$, $(10.0, 0.0)$, and $(10.5, 0.0)$, respectively. Similarly, probes are located at $(9.0,5.0)$, $(10.0, 6.0)$ and $(10.5, 6.5)$ in incidence angle of 30 (deg). The density field reaches almost the mean flow value towards the end of the PML region with very small oscillations for both tests, indicating the exponential damping of the PML formulation.      

\begin{figure}[htb]
\begin{center}
 \begin{subfigure}{0.45\textwidth}
  \includegraphics[width=0.99\textwidth]{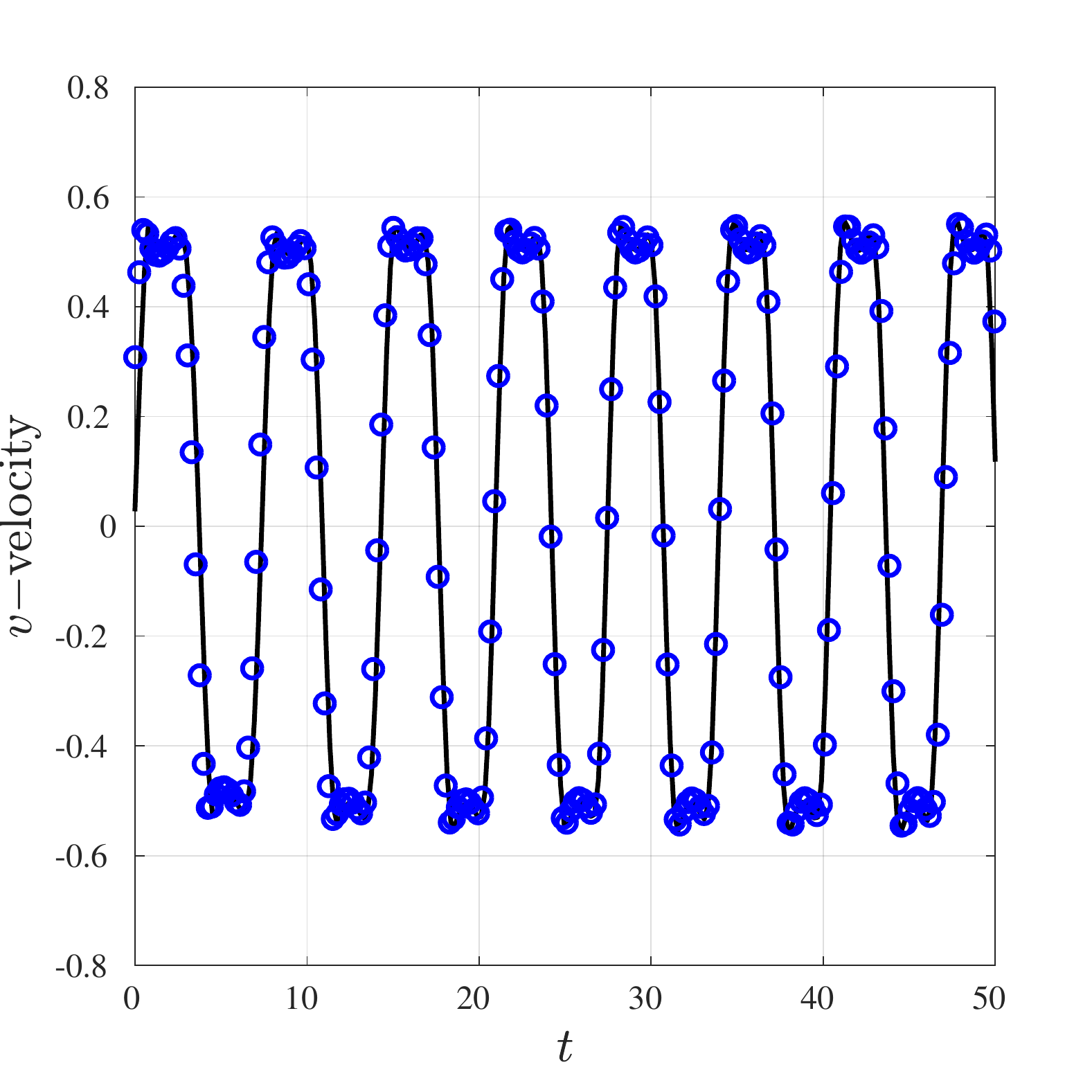}
  \caption{}
 \end{subfigure}
~
 \begin{subfigure}{0.45\textwidth}
  \includegraphics[width=0.99\textwidth]{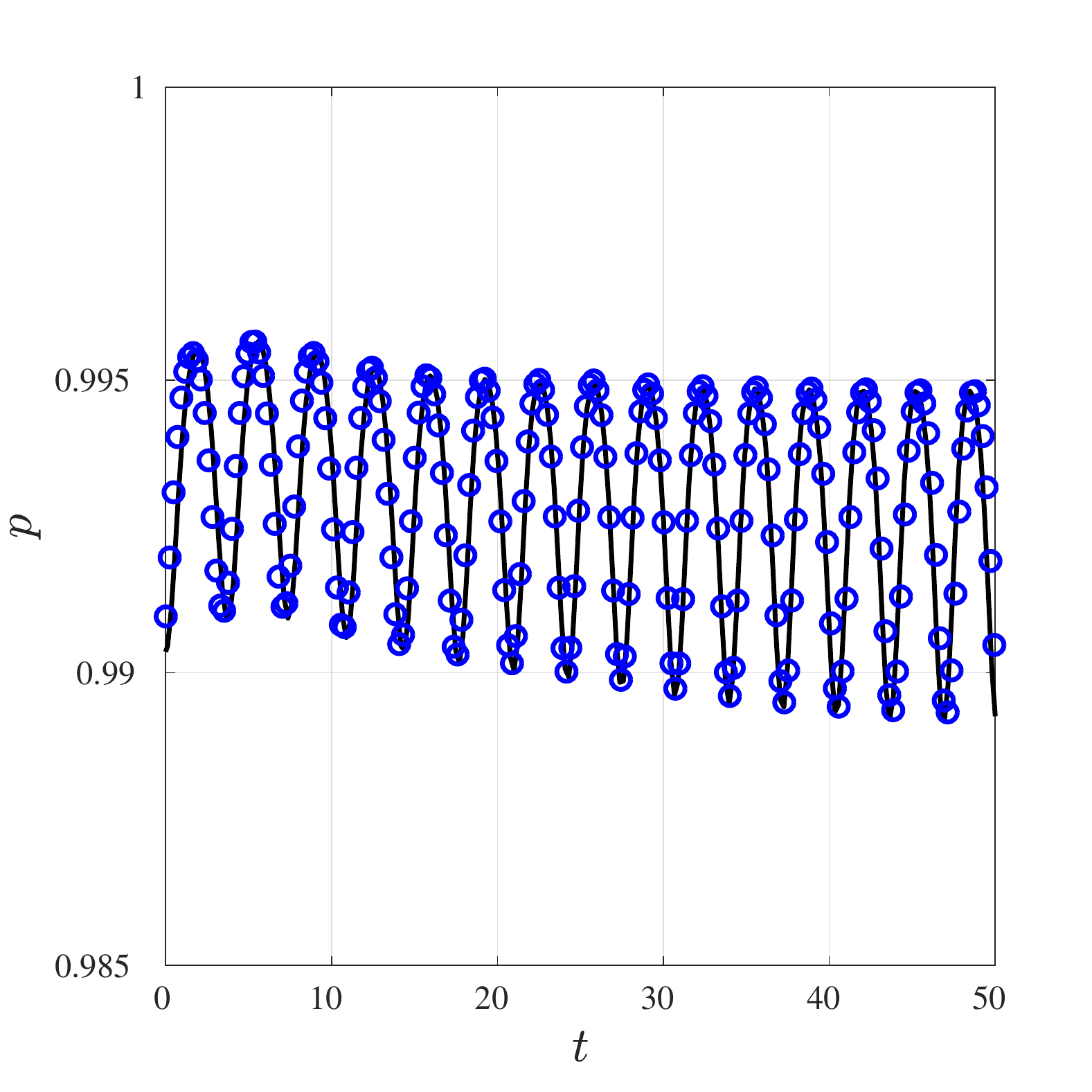}
  \caption{}
  \end{subfigure}
\caption{Flow around square cylinder test problem for $Re=150$, $N=3$ on the domain of $[-7,11] \times [-7, 7]$ with a PML region of width, $w=3.2$. Time history of (a) $y$-velocity and (b) pressure for the PML solution (solid line) and the reference solution (circle).}
\label{fig:SQCylinderProbe}
\end{center}
\end{figure}
In Figure \ref{fig:SQCylinderProbe}, $y$-velocity and pressure are shown at a point $(9.0, 0.0)$ on the outflow side of the computational domain for $Re=150$ and $N=3$. Also, the reference solution is plotted in symbols. The reference solution is obtained on a very large computational domain to ensure that any reflected waves do not pollute the solution in the probe location. The instability mechanism that triggers vortex shedding is extremely sensitive to infinitesimal noise \citep{sohankar_low-reynolds-number_1998}. Because changes in the mesh resolution, blockage, upstream/downstream extent, time step size etc. change the onset time, it is difficult to obtain the same shedding profile between reference and PML solutions. In the figure, we present the results after fully periodic pattern is observed in the $v$-velocity at the probe location from $t=0$ to $t=50$. A very good agreement in the time history of the periodically shed vortices is observed both in $y$-velocity and pressure.

As a final PML accuracy test, we compare the PML solution with the large domain solution in terms of physical parameters of vortes shedding namely Strouhal number, $St$, drag coefficient, $C_d$, and lift coefficient, $C_l$. The Strouhal number is given as $St = f d/U_\infty $ where $f$ is the frequency of shedding computed via a spectral analysis of lift coefficient history sampled over $t=100$ to $t=300$. $C_d$ and $C_l$ are computed using total force on acting on the cylinder surface, $\Gamma$,
\[\mathbf{F}_t = \int_{\Gamma}\left(-\mathbf{\sigma}\cdot\mathbf{n} + \mathbf{n}p\right)d\Gamma.\]
where pressure $p$ is recovered from the equation of state and $\sigma$ is the stress tensor with components,
\[\mathbf{\sigma}_{11} = -RT\left(\sqrt{2}q_5 - \frac{q_2^2}{q_1}\right), \quad \mathbf{\sigma}_{22} = -RT\left(\sqrt{2}q_6 - \frac{q_3^2}{q_1}\right),\]
\[\mathbf{\sigma}_{12} = \mathbf{\sigma}_{21} = -RT\left(q_4 - \frac{q_2 q_3}{q_1}\right).\]
Then, $C_d$ and $C_l$ are computed as follows,
\[ C_d = \frac{\mathbf{F}_t \cdot \mathbf{i}}{\frac{1}{2}\rho_\infty U_\infty d}, \quad C_l = \frac{\mathbf{F}_t \cdot \mathbf{j}}{\frac{1}{2}\rho_\infty U_\infty d}, \]
where $\mathbf{i}$ and $\mathbf{j}$ are the unit normal vectors in the $x$ and $y$ directions, respectively. In Table \ref{Table:SQCylinder}, we show the results for various PML widths and PML strengths with corresponding relative errors computed according to reference solution. Increasing the PML width for fixed a PML strength of $\sigma_{\max} = 20$ decreases the error in the Strouhal number $St$ where the result is obtained to be contain around $1\%$ error for the smallest PML width. On the other hand, increasing the PML damping strength does not effect the $St$ number where very accurate results obtained for $w=3.2$ for various $\sigma$ values. Increasing PML width improves the solution in $C_d$ and $C_l$. However, using a more aggressive PML damping on fixed PML width increases the error due to the need of resolving higher gradients and higher damping near the PML interface. Table \ref{Table:SQCylinder} also shows that the damping performance of our PML formulation is not strongly dependent on the PML width and strength and less than $2\%$ error is achieved even with small PML widths.  

\begin{table}[t]
	\caption{Square cylinder test problem for $Re=150$, \emph{Ma}$=0.2$ and $N=3$. Accuracy of PML formulation in terms of physical averaged quantities, $St$, $C_d$ and $C_l$. }
	\centering
	\label{Table:SQCylinder}
	\small\addtolength{\tabcolsep}{5pt}
	\begin{tabular}{c c c c c c c}
		\hline \hline
		PML Parameter   &$St$ &\%&$C_d$ & \%& $C_l$ & \% \\ 
		\hline \hline
  $w=1.6$ & $0.155$ &$1.113$&$1.462$ & $1.358$ & $0.400$ &$0.688$ \\
  $w=2.4$ & $0.154$ &$0.462$&$1.453$ & $0.755$ & $0.399$ & $0.540$ \\
  $w=3.2$ & $0.153$ &$0.069$&$1.448$ & $0.416$ & $0.398$ & $ 0.305$ \\ \hline
  $\sigma_{\max}=100$  & $0.153$ & $0.069$   & $1.451$ & $0.590$ & $0.399$ & $ 0.448$ \\
  $\sigma_{\max}=200$ & $0.153$ & $0.069$   & $1.454$ & $0.812$ & $0.400$ & $0.564$ \\
	\hline
	\multicolumn{1}{c}{Reference Solution} & $0.153$ & - & $ 1.442$ & - & $0.398$ & - \\    
	\hline \hline
	\end{tabular}
\end{table}

\subsection{Wall Mounted Fence Problem}
In this test, the efficiency of the MRSAAB time stepping method on different flow conditions is studied through solving a two-dimensional wall mounted square cylinder test problem. We select a geometric configuration with a fixed aspect ratio , i.e. ratio of cylinder height to cylinder width, is $5$.  We select \emph{Ma}$ = U_{\infty}/a_{\infty}$ and $Re= U_{\infty} d/ v_{\infty} $ where $U_{\infty}$, $a_{\infty}$, $v_{\infty}$ and $d$ are velocity of the uniform flow, speed of sound, kinematic viscosity, and cylinder height as the characteristic length, respectively.

\begin{figure}[htb!]
\begin{center}
 \begin{subfigure}{0.75\textwidth}
  \includegraphics[width=0.99\textwidth]{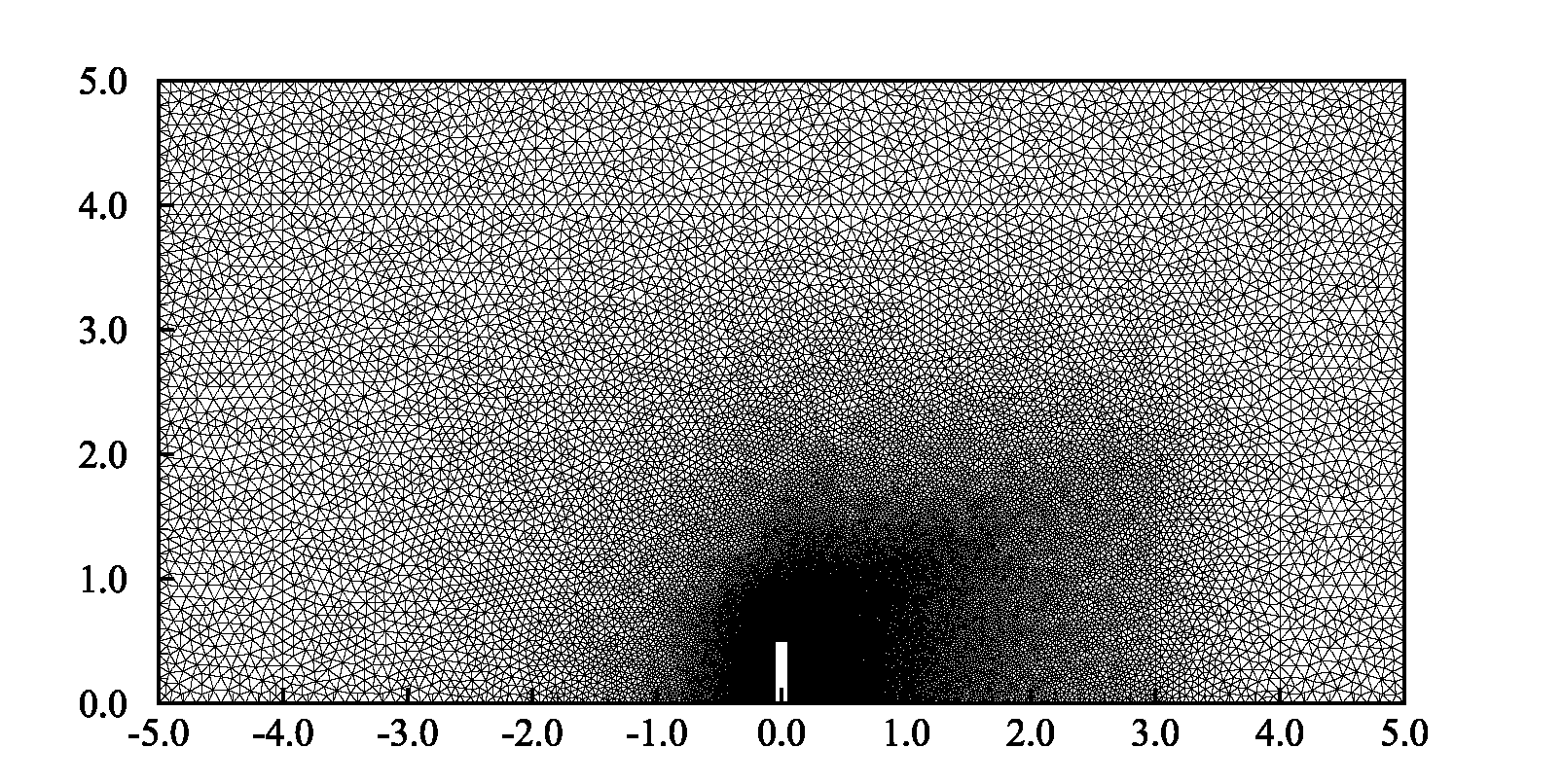}
 \end{subfigure}
\caption{Flow around wall mounted fence test problem. Discretization of the domain in the size of $[-5.0,5.0] \times [0.0, 5.0]$ with a PML region of width, $w=1.0$ into a mesh consisting of 50553 triangular elements. Elements are clustered around the fence and the wake side and uniform element size is used in the PML. }
\label{fig:FenceMesh}
\end{center}
\end{figure}
The speedup in overall runtime that can be achieve with a multirate time stepping method strongly depends on
the distribution of the characteristic stable time step sizes among the elements of the mesh. The gap between the minimum and the maximum stable time steps, as well as the number of elements present in each multirate group, has a significant influence on the computational efficiency. To show the performance of the proposed MRSAAB time discretization method for the different flow configurations in Galerkin-Boltzmann formulation, we fix the mesh and its element organization strategy for this test problem. The computational domain is chosen to be $[-5.0,5.0] \times [0, 5.0]$ with a surrounding PML region of width $w=1.0$. The domain is discretized with $50336$ unstructured triangular elements. To resolve the complex flow structure accurately, resolution is concentrated around the cylinder and in the wake region. The maximum element length is $5$ times the minimum characteristic element length. The mesh structure used in all tests is plotted in Figure \ref{fig:FenceMesh}. Figure \ref{fig:SQCylinderVorticity} gives the instantaneous vorticity contours at time $t=20$ and $t=25$ for $Re=1000$, \emph{Ma}$=0.05$, approximation order $N=5$, and with the same PML settings with used with the square cylinder test problem above.  
\begin{figure}[htb!]
\begin{center}
 \begin{subfigure}{0.6\textwidth}
  \includegraphics[width=0.99\textwidth]{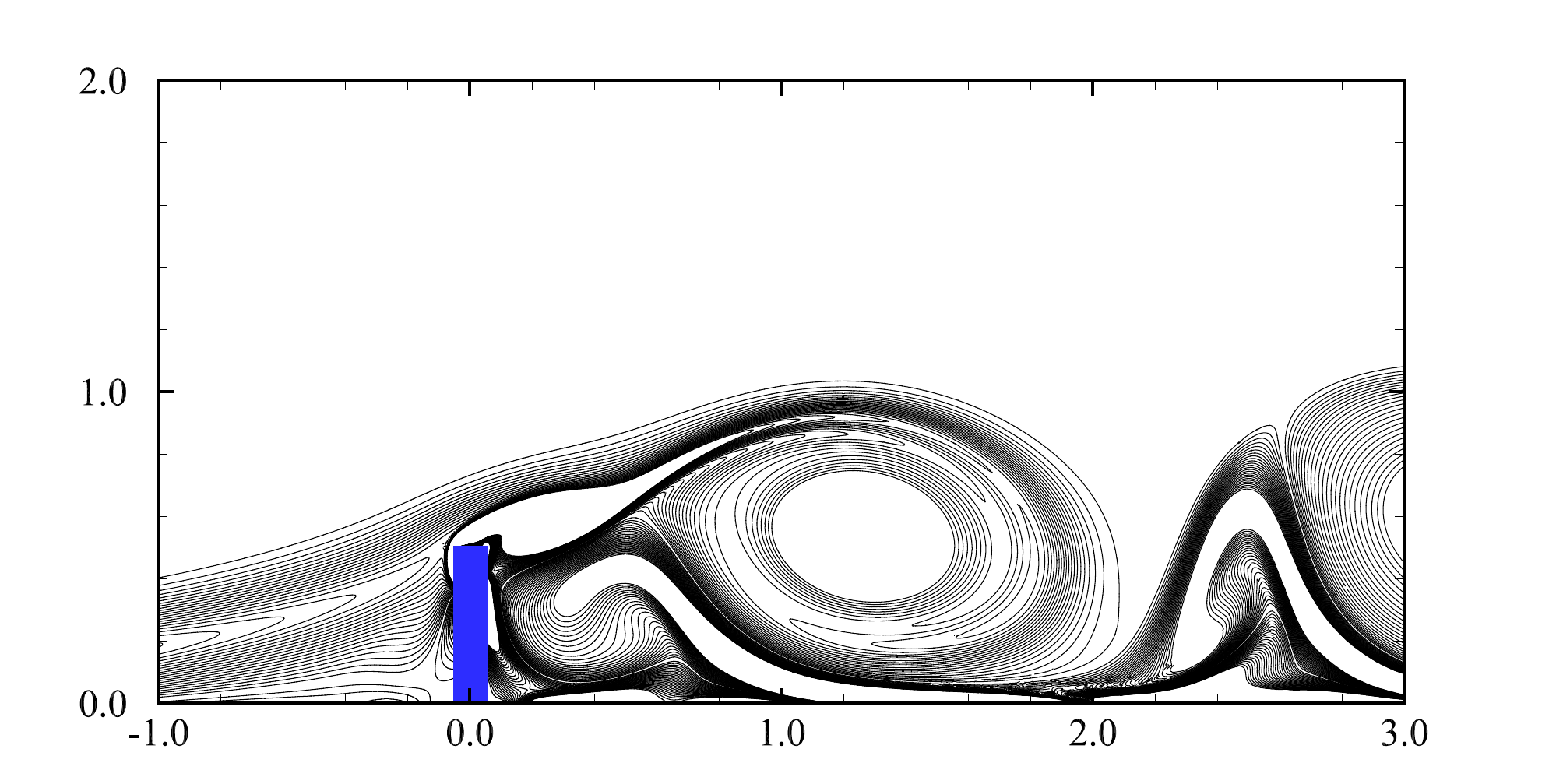}
  \caption{}
 \end{subfigure}
~
 \begin{subfigure}{0.6\textwidth}
  \includegraphics[width=0.99\textwidth]{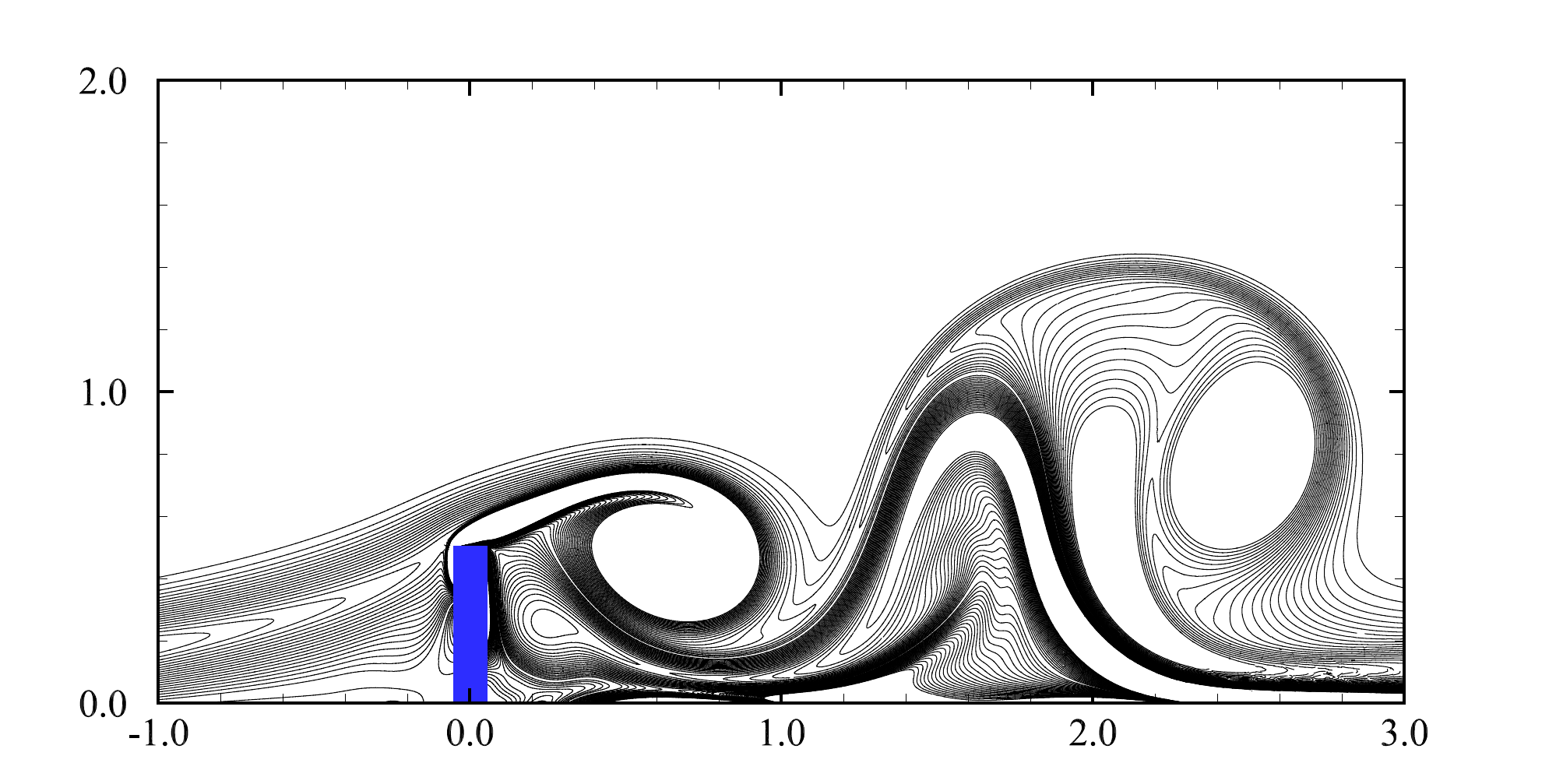}
  \caption{}
  \end{subfigure}
\caption{Flow around wall mounted square cylinder test problem for $Re=1000$, $Ma = 0.05$, $N=3$ on the domain $[-5,5] \times [0, 5]$ with a PML region of width $w=1.0$. The contours show the vorticity field from $-5.0$ to $5.0$ with the increment of $0.2$ excluding the zero level for (a) $t=20$ and (b) $t=25$. Only a part of the domain is shown.}
\label{fig:FenceVorticity}
\end{center}
\end{figure}

Table \ref{table:Fence2D} shows the number of groups and the number of elements each group for the MRSAAB and the standard MRAB time discretization approaches. As noted earlier, the stiff parameter $\frac{1}{\tau}$ depends only on the $Re$ and \emph{Ma} numbers and scales as $Re/$\emph{Ma}$^2$. Consequently, if a fully explicit time discretization is employed in the stiff regime $\frac{1}{\tau} >> 1$ the maximal time step size will be restricted in all elements of the mesh. In this case, a multirate partitioning will not create multiple levels. In contrast, the time step size of MRSAAB discretization is determined solely by the time scale of the wave transport operator and a multirate partitioning strategy will successfully create multirate levels, independent of the flow conditions. 
\begin{table}[t]
	\caption{Number of groups and element numbers in each group for the MRSAAB and MRAB time discretizations for the wall mounted cylinder problem. }
	\centering
	\label{table:Fence2D}
	\small\addtolength{\tabcolsep}{5pt}
	\begin{tabular}{c c c c c }
		\hline \hline
		Method   & $Re$ &$Ma$ & $N_{l}$ & $\#$ Elements in groups\\ 
		\hline \hline
        MRSAAB & -  &  -  & $5$ & $8396$, $13926$, $11926$, $10893$, $5412$ \\ \hline
        \multirow{5}{*}{MRAB} &$200$ & $0.05$ & $1$ & $50553$ \\ 
                              &$200$ & $0.1$ & $2$ & $8396$, $42157$ \\
                              &$200$ & $0.2$ & $3$ & $8396$, $13926$, $28231$ \\ \cline{2-5}
                              &$100$ & $0.1$ & $3$ & $8396$, $13926$, $28231$ \\
                              &$1000$ & $0.1$ & $1$ & $50553$ \\
	\hline \hline
	\end{tabular}
\end{table}

\begin{figure}[htb!]
\begin{center}
 \begin{subfigure}{0.45\textwidth}
  \includegraphics[width=0.99\textwidth]{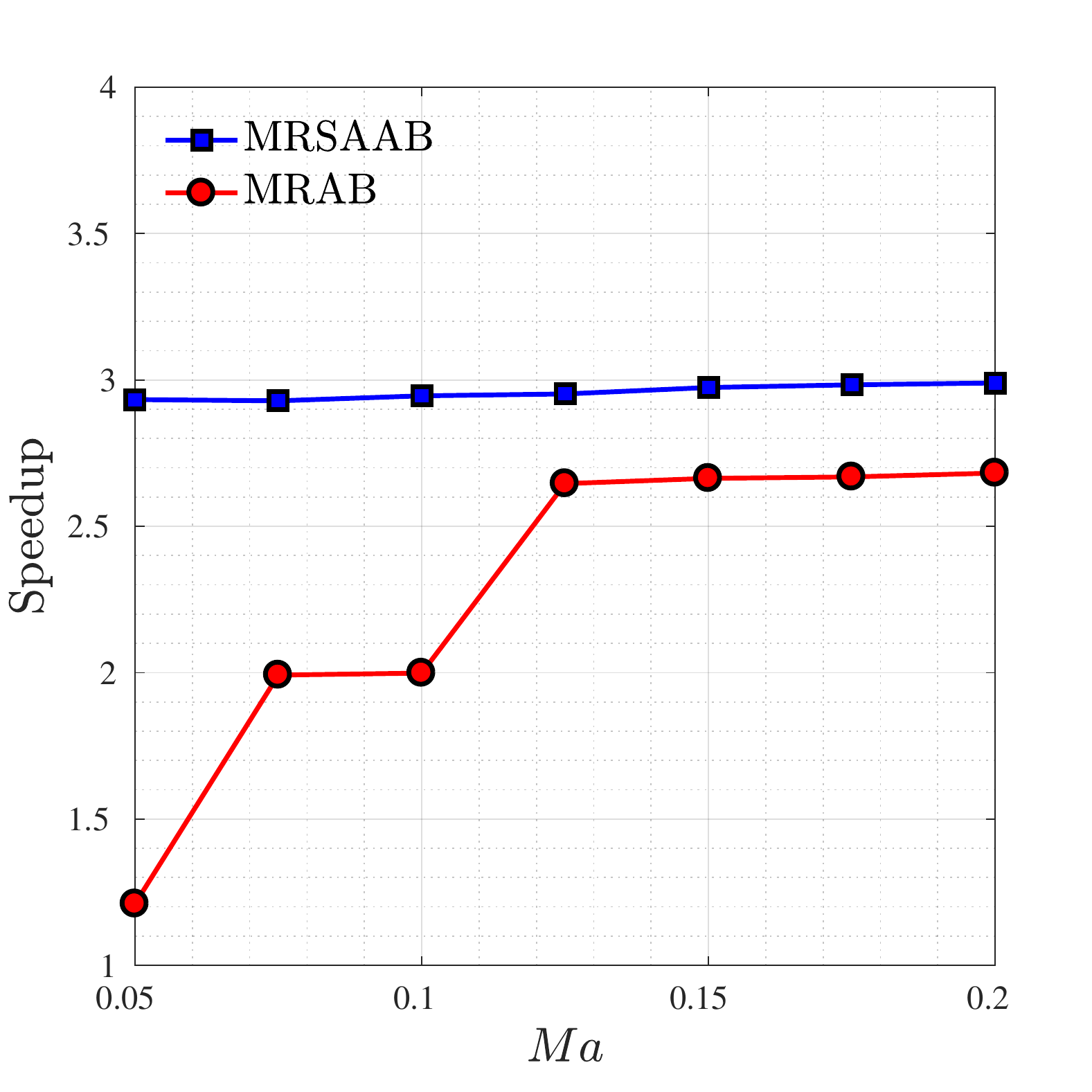}
  \caption{}
 \end{subfigure}
~
 \begin{subfigure}{0.45\textwidth}
  \includegraphics[width=0.99\textwidth]{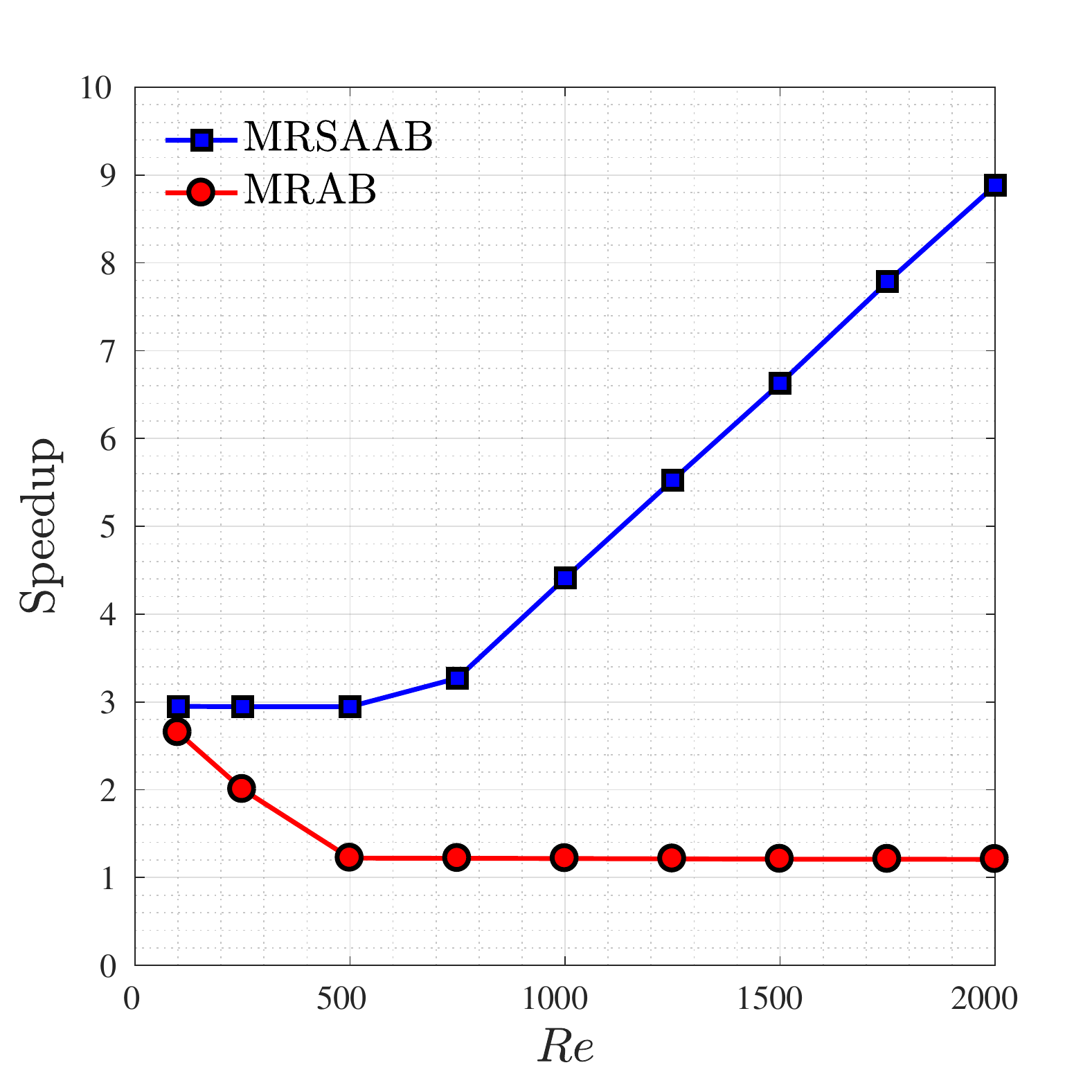}
  \caption{}
  \end{subfigure}
\caption{Multirate speedups for the flow around wall mounted square cylinder test problem on the domain of $[-5,5] \times [0, 5]$ with a PML region of width, $w=1.0$. Effect of (a) various $Ma$ numbers for $Re=200$ and (b) various $Re$ numbers for $Ma=0.1$  on the performance of the MRSAAB and MRAB time steppers. Speedups are computed relative to the LSERK time stepping scheme.}
\label{fig:FenceSpeedup}
\end{center}
\end{figure}
Figure \ref{fig:FenceSpeedup} shows the speedups achieved with MRSAAB and MRAB time discretizations for various flow configurations. The effective speedup values are computed taking the ratio of solution times of MRSAAB and MRAB schemes to the corresponding LSERK scheme. All the time stepping methods are advanced with their maximum time step sizes and LSERK uses roughly $3$ times larger CFL numbers due to its larger stability region. Figure \ref{fig:FenceSpeedup}(a) gives the speedups for $Re=200$ and \emph{Ma} numbers from $0.05$ to $0.2$. In this regime, the time step size restricted by the advective time scale where around a $3$ fold speedup of the MRSAAB scheme originates from building several multirate groups. For the Mach number \emph{Ma}$= 0.05$ the MRAB scheme creates only one multirate group resulting in only a $1.2\times$ speedup. Increasing the Mach number, and hence decreasing the value of the stiff term, results in the MRAB scheme becoming more efficient as it creates $2$ or $3$ groups, gaining $2.7$ fold speedups. Similarly, the effect of varying the Reynolds number $Re$ for a fixed Mach \emph{Ma}$=0.1$ on the performance multirate time steppers is given in Figure \ref{fig:FenceSpeedup} (b). When the time step size is restricted by the advective time scale i.e., $Re<500$, the MRSAAB method gives around a $3$ fold relative speedup. For higher Reynolds numbers, the stiffness resulting from larger $1/\tau$ factors also increases which makes the pure explicit schemes inefficient. The MRSAAB method reaches around a $9$ fold speedup in this regime. In contrast, using a fully explicit multirate approaches loses its efficiency and the solution process is accelerated only $1.2$ times for $Re>500$.

\begin{figure}[htb!]
\begin{center}
 \begin{subfigure}{0.48\textwidth}
  \includegraphics[width=0.99\textwidth]{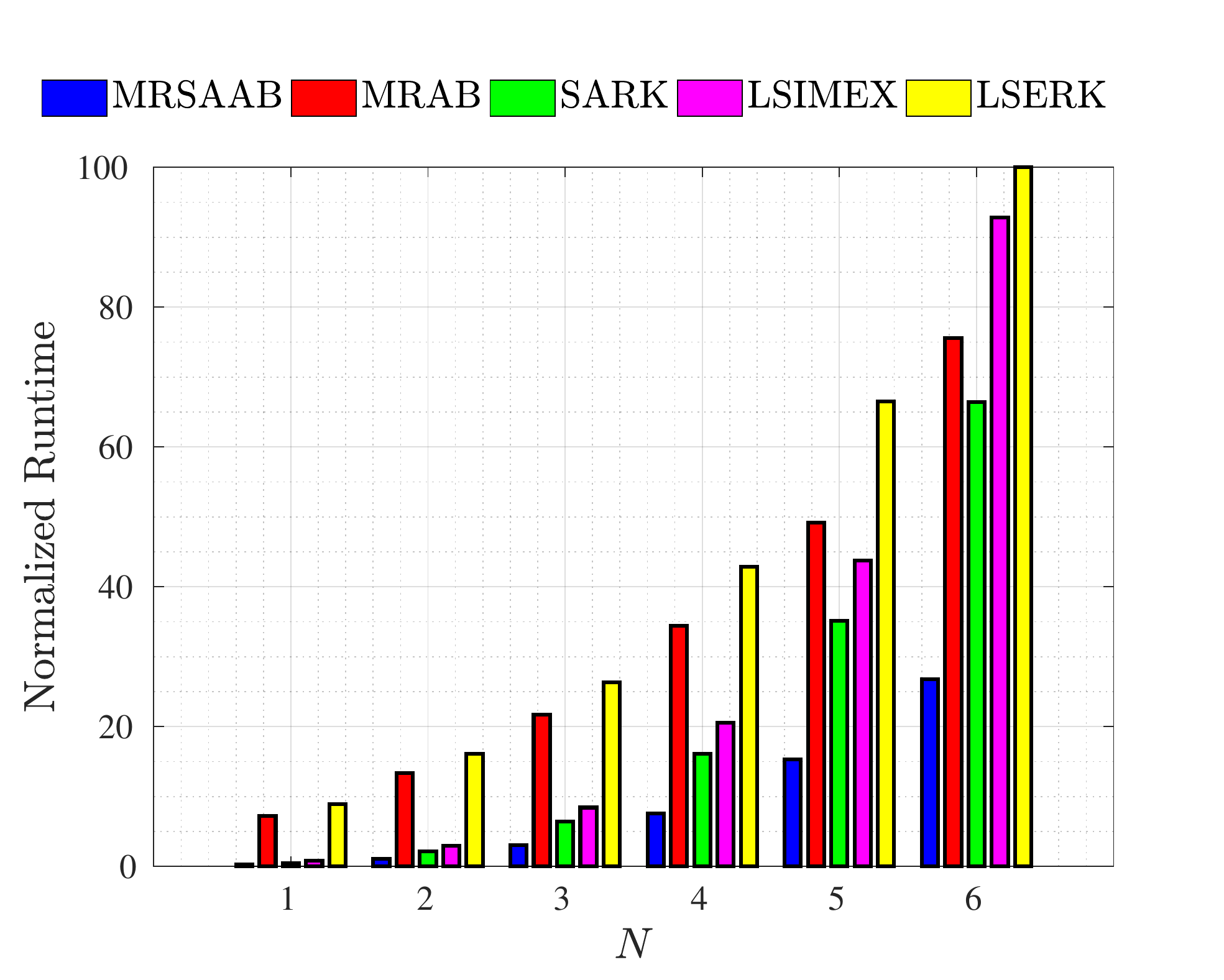}
 \end{subfigure}
 ~
\begin{subfigure}{0.48\textwidth}
  \includegraphics[width=0.99\textwidth]{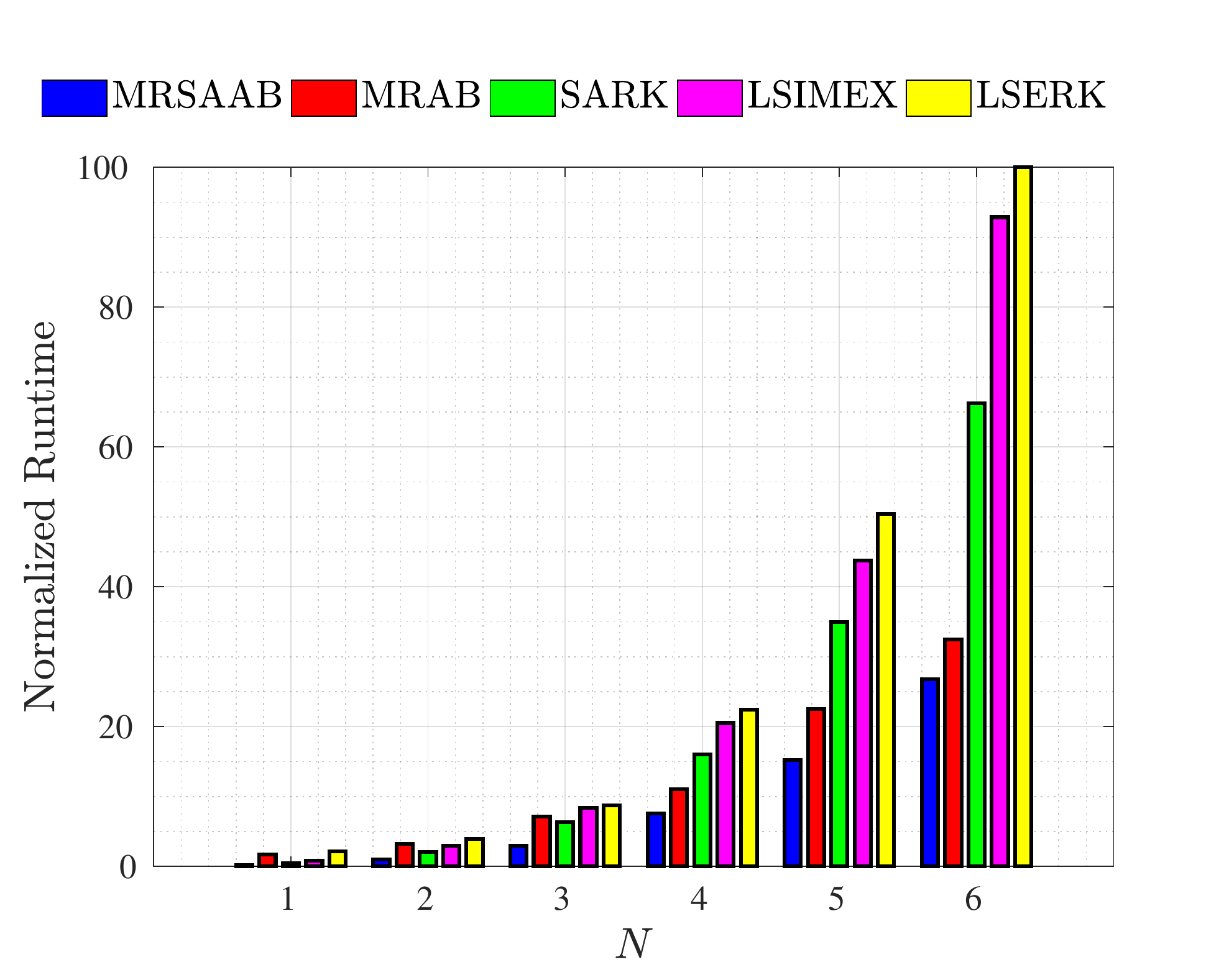}
 \end{subfigure}
\caption{Multirate speedups for the flow around a wall mounted square cylinder test problem on the domain $[-5,5] \times [0, 5]$ with a PML region of width $w=1.0$. Normalized runtimes for (top) $Ma=0.05$ and $Re=1000$  (bottom) $Ma=0.2$ and $Re=200$. Runtimes are normalized according to the runtime of the LSERK method at $N=6$.}
\label{fig:SQRunTime}
\end{center}
\end{figure}

Figure \ref{fig:SQRunTime} illustrates the breakdown of the normalized runtimes for different time discretization methods for orders $N=1\ldots6$. For $Re=1000$ and $Ma=0.05$, the Boltzmann equations become stiff and MRSAAB method results in considerably lower runtimes as given in \ \ref{fig:SQRunTime}(a). In this regime, the multirate method without semi-analytic integration is not effective due to lack of sufficient grouping where even SARK scheme outperforms MRAB. The LSIMEX method produces slightly longer runtimes in all orders compared to the SARK method due to additional operations and data movement. The fully explicit scheme gives the highest computational times as expected. Figure \ref{fig:SQRunTime} presents the same results for $Re=200$ and $Ma=0.2$. The Boltzmann equations are not as stiff in this regime and the MRAB scheme creates $3$ groups for $N=5$. The MRSAAB scheme gives the smallest run times but this time MRAB outperforms SARK and LSIMEX schemes for all orders. For these two different flow regimes, the MRSAAB method is the fastest method and its performance is independent of the flow regime.   

\section{Conclusion}
\label{Sec:Conclusion}
We presented a high-order nodal discontinuous Galerkin method for the Boltzmann equations discretized with Hermite polynomials in velocity space, and used it to simulate nearly incompressible flows. We also introduced a stabilized unsplit perfectly matching layer (PML) formulation for the resulting nonlinear flow equations. The equations are advanced in time with developed semi-analytic and multirate schemes. Numerical tests show that the proposed M-PML formulation exponentially damps the difference between the nonlinear fluctuation and a prescribed mean flow. Because of the non-linearity in the equations, the absorbing layers are not formally perfectly matched with the governing equations, in contrast to their linear counterparts. However, numerical examples give satisfactory results even in the severely truncated domains. We tested the performance of developed semi-analytic time integration in terms of accuracy and efficiency and compared the performance with  implicit-explicit and fully explicit Runge-Kutta methods. Numerical result indicate that the performance of the multirate semi-analytic method combined with the Galerkin-Boltzmann equations is very promising for modeling physically relevant flow problems requiring spatially and temporally varying scales.      

\section*{Acknowledgements}
This research was supported in part by the Exascale Computing Project (17-SC-20-SC), a collaborative effort of two U.S. Department of Energy organizations (Office of Science and the National Nuclear Security Administration) responsible for the planning and preparation of a capable exascale ecosystem, including software, applications, hardware, advanced system engineering, and early testbed platforms, in support of the nations exascale computing imperative.

The computational experiments reported in this paper were performed on systems provided by the Advance Research Computing group at Virginia Tech. Finally, this research was supported in part by the John K. Costain Faculty Chair in Science at Virginia Tech.

\printbibliography

\section*{Appendix: Derivation of Couette Flow for Galerkin-Boltzmann Equations}
\label{Sec.Appendix}
To derive an analogous Couette flow solution for the Galerkin-Boltzmann equation, we start from a shear flow assumption i.e., we assume that that ${\bf q}={\bf q}(y,t)$. Simplifying Equation \ref{eq:BoltzmannSystem}, we obtain
 \begin{eqnarray*}
 \frac{\partial \mathbf{\tilde{q}}}{\partial t} & = & A_y \frac{\partial \mathbf{\tilde{q}}}{\partial y} + \mathcal{N}(q).
 \end{eqnarray*}
We further assume uniform density, namely $q_1=1$, and horizontal flow $q_3 = q_6 =0$ to obtain the following system,
  \begin{eqnarray*}
  \frac{\partial q_2}{\partial t} & = & -\sqrt{RT} \frac{\partial q_4}{\partial y}, \\
    \frac{\partial q_4}{\partial t} & = & -\sqrt{RT} \frac{\partial q_2}{\partial y} - \frac{1}{\tau} q_4,\\
     \frac{\partial q_5}{\partial t} & = & - \frac{1}{\tau} \left( q_5 - \frac{q_2^2}{\sqrt{2}}\right).
 \end{eqnarray*}
We are interested in the $y$-velocity profiles of this shear flow and therefore focus our attention on obtaining a $q_2$ which solves this system. To begin, we note that the first two equations correspond to the telegrapher's equation. We eliminate $q_4$ by differentiating the first equation in respect to $t$. Substituting the second equation, we obtain
 \begin{eqnarray*}
  \frac{\partial^2 q_2}{\partial t^2} & = & {RT} \frac{\partial^2 q_2}{\partial y^2} - \frac{1}{\tau}\frac{\partial q_2}{\partial t}.
  \end{eqnarray*}
This equation admits the trivial solutions $q_2 = 1$ and $q_2 = y$ which are sufficient to satisfy the shear boundary conditions $u|_{y=0} = 0$ and $u|_{y=L} = U$. We find additional general solutions to this equation by finding solutions that satisfy homogeneous boundary conditions. These solutions can be found by assuming the separable form 
\[
q_2^{n} = \sin\left(\lambda_n y\right) e^{\sigma_n t},
\]
where $\lambda_n = \frac{n\pi}{L}$ and $\sigma_n$ satisfies $\sigma_n^2 = -RT\lambda_n^2 - \frac{1}{\tau}\sigma_n$. Hence,
\begin{eqnarray*}
\sigma_n & = & -\frac{1}{2\tau} \pm \sqrt{\frac{1}{4\tau^2}-RT\lambda_n^2},\\
& = & -\frac{1}{2\tau} \pm \frac{1}{2\tau}\sqrt{1-4\tau^2RT\lambda_n^2}.
\end{eqnarray*}
Assuming that $\tau^2RT\lambda_n^2$ is small, we can use a small parameter estimate to write
\begin{eqnarray*}
\sigma_n & = & -\frac{1}{2\tau} \pm \left(\frac{1}{2\tau}-\tau RT\lambda_n^2 + \mathcal{O}(\tau^3 R^2T^2\lambda_n^4)\right).
\end{eqnarray*}
Therefore, to select the branch that best approximates the incompressible Navier-Stokes shear mode, we take the positive branch, i.e.,
\begin{eqnarray*}
\sigma_n & = &  -\frac{1}{2\tau} + \frac{1}{2\tau}\sqrt{1-4RT\lambda^2_n\tau^2}.
\end{eqnarray*}
so that 
\begin{eqnarray*}
\sigma_n & = &  -\tau RT\lambda_n^2 + \mathcal{O}(\tau^3 R^2T^2\lambda_n^4), \\
      & = & -\nu\lambda_n^2 + \mathcal{O}\left(\frac{\nu^3}{a^2}\lambda_n^4\right),
\end{eqnarray*}
by the definition of fluid viscosity $\nu$ and the speed of sound $a$. We note that the relative discrepancy between the incompressible Navier-Stokes and Galerkin-Boltzmann shear mode is
\[
\epsilon_n := \frac{\sigma_n + \nu \lambda^2_n}{\nu \lambda^2_n} = \mathcal{O}\left(\frac{\nu^2}{a^2}\lambda_n^2\right). 
\]
Thus, we know a priori that the decay rates of linear shear for these flow models diverge for sufficiently high-order modes. On the other hand, the decay rates will agree well for shear flows with low viscosity, large Mach number, and low mode numbers.

Using the equation for $q_4$, we find that each homogeneous solution $q_2^n$ has an associated solution $q_4^n$ component,
\[
q_4^n = \frac{\sigma_n}{\sqrt{RT}\lambda_n} \cos(\lambda_n y)e^{\sigma_n t}.
\]
Using the homogeneous and non-homogeneous solutions for $q_2$, we write the exact solution of the horizontal momentum of the Galerkin-Boltzmann shear flow as an expansion, which satisfies the initial condition $q_2(y,0) = 0$ and the shear boundary conditions $q_2(0,t) = 0$ and $q_2(L,t) = \frac{U}{\sqrt{RT}}$ as follows
\begin{eqnarray}
q_2(y,t) &=& \frac{U}{\sqrt{RT}L} y + \frac{1}{\sqrt{RT}}\sum_{n=1}^\infty \frac{2(-1)^n U}{\lambda_n L} \sin(\lambda_n y)e^{\sigma_n t}.
\end{eqnarray}
The associated solution for $q_4$ is written as
\begin{eqnarray}
q_4(y,t) &=& \frac{1}{RT}\sum_{n=1}^\infty \frac{2(-1)^n U\sigma_n}{\lambda_n^2 L} \cos(\lambda_n y)e^{\sigma_n t}. 
\end{eqnarray}
Finally, once we return to the equation for $q_5$, we find that $q_5$ satisfies 
\begin{align*}
\frac{\partial q_5}{\partial t} &= -\frac{1}{\tau}q_5 + \frac{1}{\sqrt{2}\tau}q_2^2, \\
& = -\frac{1}{\tau}q_5 + \frac{1}{\sqrt{2}\tau}\left( \frac{U}{\sqrt{RT}L} y + \frac{1}{\sqrt{RT}}\sum_{n=1}^\infty \frac{2(-1)^n U}{\lambda_n L} \sin(\lambda_n y)e^{\sigma_n t}\right)^2.
\end{align*}
This equation is solved using the homogeneous boundary conditions $q_5(0,t) =0$ and $q_5(L,t) = \frac{U^2}{\sqrt{2}RT}$ to obtain
\begin{multline*}
q_5(y,t) = \frac{U^2}{\sqrt{2}RTL^2}y^2 + \sum_{n=0}^\infty c_n \sin(\lambda_n y)e^{-\frac{t}{\tau}} \\ +  \frac{1}{\sqrt{2}RT}\sum_{n=1}^\infty \frac{2(-1)^n U^2}{\lambda_n L^2(\sigma_n\tau+1)}\sin(\lambda_n y)e^{\sigma_n t} \\ + \frac{1}{\sqrt{2}RT}\sum_{n=1}^\infty\sum_{m=1}^\infty \frac{4(-1)^n U^2}{\lambda_n\lambda_m L^2(\sigma_n\tau+\sigma_m\tau+1)}\sin(\lambda_n y)\sin(\lambda_m y)e^{(\sigma_n+\sigma_m) t},
\end{multline*}
where $c_n$ are coefficients chosen to satisfy the initial condition for $q_5$. 

\end{document}